\documentclass[a4paper,11pt,reqno]{amsart}
\usepackage[T1]{fontenc}
\usepackage[utf8]{inputenc}
\usepackage{lmodern}
\usepackage{amsmath,amssymb,amsthm,mathtools,mathrsfs,enumerate}
\usepackage{microtype}
\microtypesetup{expansion=false}
\usepackage[hidelinks]{hyperref}
\hypersetup{pdftitle={Eisenstein series twisted Shintani zeta functions over number fields},pdfauthor={Eun Hye Lee and Ramin Takloo-Bighash}}
\newtheorem{theorem}{Theorem}[section]
\newtheorem{proposition}[theorem]{Proposition}
\newtheorem{lemma}[theorem]{Lemma}
\newtheorem{corollary}[theorem]{Corollary}
\theoremstyle{definition}

\theoremstyle{remark}
\newtheorem{remark}[theorem]{Remark}
\newcommand{\A}{\mathbb A}
\newcommand{\C}{\mathbb C}
\newcommand{\R}{\mathbb R}

\newcommand{\OO}{\mathcal O}
\newcommand{\cS}{\mathscr S}

\newcommand{\one}{\mathbf 1}
\newcommand{\diag}{\operatorname{diag}}
\newcommand{\vol}{\operatorname{vol}}

\newcommand{\GL}{\mathrm{GL}}

\allowdisplaybreaks[2]

\title[Eisenstein twists of Shintani zeta functions]{Automorphic Form Twisted Shintani Zeta Functions over Number Fields II: Eisenstein Series}
\author{Eun Hye Lee}
\address{Department of Mathematics, Texas Christian University, Fort Worth, TX 76129}
\email{eun.hye.lee@tcu.edu}
\author{Ramin Takloo-Bighash}
\address{Department of Mathematics, Statistics, and Computer Science, University of Illinois at Chicago, Chicago, IL 60607}
\email{rtakloo@uic.edu}
\date{26 November 2026}
\begin{document}
\begin{abstract}
We study Shintani zeta integrals twisted by Eisenstein series over an arbitrary number field. 
\end{abstract}
\maketitle
\tableofcontents

\section{Introduction and main theorem}\label{sec:statement}

In this paper, following \cite{LTcuspidal}, we consider the twist of the zeta function considered by Shintani \cite{Shintani} and Wright \cite{Wright} by Eisenstein series. The cuspidal case is treated in the companion preprint \cite{LTcuspidal}. The nonzero constant term of the Eisenstein series prevents the singular theta pairing from being an ordinary convergent integral, even on a fixed determinant fiber. The principal task is to regularize that pairing compatibly with unfolding, Poisson summation, and spectral specialization.

Hough and Lee \cite{HoughLee} established Eisenstein-twisted continuation and pole formulas over $\mathbb Q$ in the spherical trivial central character case.  We generalize their result to a general number field and allow arbitrary inducing Hecke characters, compact types, and Schwartz--Bruhat tests.

The regularization follows Zagier's subtraction of integrated cusp powers \cite{ZagierRS}; related automorphic truncation and period constructions appear in \cite{Arthur,ArthurII,JLR} and \cite[Section~4.3]{MV}. The zero vector and both nonzero singular orbits use one geometric cutoff. Their intermediate cutoff powers cancel, leaving the same primitive on both sides of Poisson summation. This equality is essential: separately chosen orbit regularizations would not by themselves justify the Poisson identity. In the remaining paragraphs of this introduction we explain exactly what we prove and outline the argument. 

\subsection{Binary cubic forms and the zeta integral}
Let $k$ be a number field, let $n=[k:\mathbb Q]$, and write $\A=\A_k$. Let $M_k=M_0\sqcup M_\infty$ be its places, partitioned into finite and archimedean places, and put $r=|M_\infty|=r_1+r_2$; thus $n=r_1+2r_2$. For an algebraic group or vector space $X$ and a $k$-algebra $R$, write $X_R=X(R)$. Set $G=\GL_2$ and identify $V=k^4$ with binary cubic forms by
\[
 F_x(u,v)=x_1u^3+x_2u^2v+x_3uv^2+x_4v^3.
\]
We use the determinant twist in Wright's adelic framework \cite{Wright}. Our left action $\rho$ is defined by
\begin{equation}\label{eq:action}
 F_{g\cdot x}(u,v)=\det(g)^{-1}F_x\bigl((u,v)g\bigr).
\end{equation}
Thus $\rho$ is the determinant twist of the symmetric-cube representation; a scalar matrix $aI_2$ acts on $V$ by multiplication by $a$. The discriminant is the degree-four polynomial
\[
 \operatorname{Disc}(x)=x_2^2x_3^2+18x_1x_2x_3x_4
 -4x_2^3x_4-4x_1x_3^3-27x_1^2x_4^2.
\]
It satisfies
\begin{equation}\label{eq:discriminant-character}
 \operatorname{Disc}(gx)=(\det g)^2\operatorname{Disc}(x).
\end{equation}

Let $S$ be the discriminant-zero hypersurface and put $V'=V\setminus S$.

Write $U_\A=\prod_v U_v$ for the product of the standard maximal
compact subgroups of $G(k_v)$; $K$-finiteness below means finiteness
under right translation by $U_\A$. Put
\[
 n(c)=\begin{pmatrix}1&0\\c&1\end{pmatrix},\quad
 a(a,d)=\operatorname{diag}(a,d),\quad d(z)=zI_2,\quad
 g^\iota=(\det g)^{-1}g.
\]

The lower Borel subgroup is denoted by $B$; its unipotent radical is $N=\{n(c)\}$.
We also write $a(t)=a(t,t^{-1})$, $w_0=\left(\begin{smallmatrix}0&1\\1&0\end{smallmatrix}\right)$, and $n^+(c)=w_0n(c)w_0$.
For a smooth automorphic form $\varphi$ on $G_k\backslash G_\A$ with unitary central character $\omega$, and $\Phi\in\cS(V_\A)$, define initially in a right half-plane
\begin{equation}\label{eq:Z}
Z(s,\Phi,\varphi)=\int_{G_\A/G_k}|\det g|_\A^s\varphi(g^{-1})\sum_{x\in V'_k}\Phi(gx)\,dg.
\end{equation}
The measures are fixed in Section~\ref{subsec:measures}. Singular analogues require the finite part defined in Section~\ref{subsec:eis-finite-part}; their Mellin transforms are defined in \eqref{eq:76-IJ-definition}.
For $z>0$, let $\underline z$ be the idele whose finite components are $1$ and whose infinite components are $z^{1/n}$. Thus $|\underline z|_\A=z$ and $d(\underline z)$ acts on $V_\A$ by scalar multiplication by $\underline z$.

\begin{lemma}\label{lem:central-norm-normalization}
There is a unique $\kappa\in\R$ such that $\omega(\underline z)=z^{i\kappa}$.
Set $\omega_0=\omega|\cdot|_\A^{-i\kappa}$ and $\varphi_0(h)=|\det h|_\A^{-i\kappa/2}\varphi(h)$. Then $\omega_0(\underline z)=1$ and
\begin{equation}\label{eq:unitary-central-shift}
Z(s,\Phi,\varphi)=Z(s-i\kappa/2,\Phi,\varphi_0).
\end{equation}
The same translation applies to the singular Mellin transforms.
\end{lemma}
\begin{proof}
This is \cite[Lemma~1.1]{LTcuspidal}. For the finite parts used here, the additional observation is that the norm character is constant on each determinant fiber. It can therefore be moved outside both the truncated integral and its subtraction before passing to the limit.
\end{proof}
Until the norm component is restored explicitly, all central characters are normalized by $\omega(\underline z)=1$.

\subsection{The Eisenstein family and the main theorem} For Eisenstein twists fix unitary Hecke characters
$\eta_1,\eta_2$ trivial on the chosen norm section and set
$\omega=\eta_1\eta_2$. Let $f_\tau$ be a smooth right-$U_\A$-finite
flat section satisfying
\begin{equation}\label{eq:intro-inducing-section}
 f_\tau(a(a,d)n(c)h)
 =\eta_1(a)\eta_2(d)|d/a|_\A^{1/2+\tau}f_\tau(h).
\end{equation}
Flatness means $f_\tau|_{U_\A}$ is independent of $\tau$; finite
holomorphic linear combinations of such sections are allowed. Put
\begin{equation}\label{eq:intro-eisenstein}
 E_\tau(h)=\sum_{\gamma\in B_k\backslash G_k}f_\tau(\gamma h)
                           \qquad(\Re\tau>1/2),
\end{equation}
and continue meromorphically. The central character is $\omega$.
Every smooth $K$-finite Eisenstein series with unitary central
character is covered after the common norm twist of
Lemma~\ref{lem:central-norm-normalization} and absorption of the
opposite norm characters into $\tau$.

Put $(\eta_{+,1},\eta_{+,2})=(\eta_1,\eta_2)$ and
$(\eta_{-,1},\eta_{-,2})=(\eta_2,\eta_1)$. Its constant term is
$f_\tau+M(\tau)f_\tau$ and, for
$g=u a(\underline\lambda)a(t_1,t_2)n(c)$, $t_i\in\A^1$, equals
\begin{equation}\label{eq:intro-eisenstein-constant}
 E_{\mathrm c}(\tau,g^{-1})
 =\sum_{\epsilon=\pm1}\lambda^{1+2\epsilon\tau}
    \eta_{\epsilon,1}(t_1)^{-1}\eta_{\epsilon,2}(t_2)^{-1}
    c_\epsilon(\tau,u),
\end{equation}
where $c_+(\tau,u)=f_\tau(u^{-1})$ and
$c_-(\tau,u)=(M(\tau)f_\tau)(u^{-1})$.
Write $\delta(\xi)=1$ if the norm-normalized Hecke character
$\xi$ is trivial, and zero otherwise.
A parameter is regular when the Eisenstein family is holomorphic
there; regular parameters outside
$\{\pm1/2,\pm1,\pm3/2\}$ are called generic.

\begin{theorem}\label{thm:eisenstein-intro}
At every regular spectral parameter, the nonsingular integral
$Z(s,\Phi,E_\tau)$ converges absolutely in the sufficient region
\begin{equation}\label{eq:intro-eisenstein-convergence}
 \Re s>5/2+|\Re\tau|.
\end{equation}
It admits joint meromorphic continuation in $(s,\tau)$. For generic
$\tau$, its only possible poles in $\Re s>1/2+|\Re\tau|$ are
simple poles at
\begin{equation}\label{eq:intro-eisenstein-poles}
 s=11/6\pm\tau/3,\qquad s=5/2\pm\tau,
\end{equation}
when these points lie in that half-plane. In the plus branch, the
triple-root and double-root residues can be nonzero only if
$\eta_1\eta_2^2=1$ and $\eta_2=1$, respectively; in the minus
branch the corresponding conditions are $\eta_1^2\eta_2=1$
and $\eta_1=1$. 

The singular integrals admit a common cusp regularization, defined
in \eqref{eq:eis-reg-subtraction} and \eqref{eq:76-IJ-definition},
and continue jointly meromorphically. Their nonconstant contribution
has the Euler quotient
\begin{equation}\label{eq:intro-general-euler}
 \frac{L^\Sigma(3s-\tfrac12,\pi_\tau\otimes\omega^{-2})
       L^\Sigma(s+\tfrac12,\pi_\tau\otimes\omega^{-1})}
      {L^\Sigma(4s,\omega^{-2})}\,H_{\Phi,\Sigma}(s,\tau),
\end{equation}
up to the factor $1/2$ in \eqref{eq:eis-euler}. Here $\Sigma$
contains the infinite and nonstandard places, $H_{\Phi,\Sigma}$
is the integral over the places in $\Sigma$ defined in \eqref{eq:75-H-definition}, and
\begin{equation}\label{eq:intro-principal-L}
 L^\Sigma(w,\pi_\tau)
 =L^\Sigma(w+\tau,\eta_2)L^\Sigma(w-\tau,\eta_1).
\end{equation}

If all four characters
\begin{equation}\label{eq:intro-nonresonance}
 \eta_1,\quad\eta_2,\quad\eta_1^2\eta_2,\quad\eta_1\eta_2^2
\end{equation}
are nontrivial, the regularized dual singular integral is entire,
and the regularized $J$ and ordinary $Z$ are holomorphic, at every
regular spectral parameter, in
\begin{equation}\label{eq:intro-nontrivial-holomorphy}
 \Re s>\beta(\tau):=
 \max\{1/4,(1/2+|\Re\tau|)/3,|\Re\tau|-1/2\}.
\end{equation}
On the unitary axis this is $\Re s>1/4$.
If $\Phi|_{S(\A)}=0$, then for generic $\tau$ the locations
\eqref{eq:intro-eisenstein-poles} are the only possible poles on
the whole plane; under \eqref{eq:intro-nonresonance}, $Z$ is
entire at every regular parameter. Exceptional values and spectral
Laurent coefficients are described in \S\ref{subsec:eis-exceptional}.
A common norm twist translates all $s$-locations by $i\kappa/2$
as in \eqref{eq:unitary-central-shift}.
\end{theorem}

\subsection{Structure of proof and relation to \cite{LTcuspidal}}
Section~\ref{sec:radial} fixes measures and recalls the lattice bounds used in the companion paper. Section~\ref{subsec:eis-family} proves absolute convergence of the ordinary zeta integral. Section~\ref{subsec:eis-finite-part} constructs a common cusp finite part and proves its unfolding at finite cutoff. Section~\ref{subsec:eis-constant} evaluates the triple-root and double-root constant terms; the additive Jacobian in the latter is retained explicitly. Section~\ref{subsec:eis-nonconstant} completes the nonconstant radial Mellin transform and unfolds it into Whittaker integrals. Section~\ref{subsec:eis-local} continues those distributions and computes the Euler quotient. Poisson summation in Section~\ref{subsec:eis-continuation} proves the main theorem, and Section~\ref{subsec:eis-exceptional} treats exceptional values and spectral Laurent coefficients. Appendix~\ref{subsec:eis-smoothed} identifies the finite part with the residue of a smoothed pairing.

The shared analytic lemmas from \cite{LTcuspidal} are stated at their points of use, with precise references for their proofs. The spectral-family extensions needed here are justified separately. The method of radial splitting is Tate's \cite{Tate}; the singular orbit geometry is from Wright \cite{Wright}; and the Whittaker calculation is related to Hough's double Hecke series \cite[Lemma~6]{Hough}. The broader prehomogeneous-vector-space and automorphic-period framework is developed in \cite{SS,Sato94,Sato06}. For the arithmetic local theory and adelic smoothing see \cite{DW2,DW3,Yukie}.

\subsection{Zagier's method}
The principal difficulty encountered here, compared with a cuspidal
twist \cite{LTcuspidal}, arises from the constant term of the Eisenstein series.
Although the original nonsingular zeta integral converges in a
right half-plane, Poisson summation introduces singular theta
pairings that need not converge, even at a fixed radial parameter.
We treat these pairings by subtracting their explicit cusp
asymptotics. To explain the idea, we first recall the classical
renormalization procedure of Zagier
\cite[pp.~424--426]{ZagierRS}.

Let $\Gamma=\mathrm{SL}_2(\mathbb Z)$, and let $\mathcal D$ be
its standard fundamental domain in the upper half-plane.
For this classical discussion, write $w=x+iy$, with hyperbolic
measure
\begin{equation}\label{eq:ch7-intro-measure}
 d\mu(w)=\frac{dx\,dy}{y^2}.
\end{equation}
Put $\mathcal D_T=\{w\in\mathcal D:\operatorname{Im}w\le T\}$.
For $T>1$, the complementary cusp is the strip
\begin{equation}\label{eq:ch7-intro-cusp}
 \mathcal D\setminus\mathcal D_T
 =
 \left\{x+iy:-\frac12\le x\le\frac12,\ y>T\right\},
\end{equation}
up to boundary identifications. Suppose that $F$ is a continuous
$\Gamma$-invariant function with an expansion
\begin{equation}\label{eq:asymp}
 F(x+iy)=\varphi(y)+r(x,y),
 \qquad
 \varphi(y)=\sum_{j=1}^{\ell}
                  c_j y^{\alpha_j}(\log y)^{n_j},
\end{equation}
where $\alpha_j\in\mathbb C$, $n_j$ are nonnegative integers,
and $r(x,y)=O_N(y^{-N})$ for every $N$, uniformly in $x$.

Define $\widehat\varphi$ by
\begin{equation}\label{eq:ch7-intro-primitive}
 \widehat\varphi'(y)=\frac{\varphi(y)}{y^2},
\end{equation}
choosing the unique antiderivative that is a finite linear
combination of power--logarithm functions with no standalone
constant term. Here the hat denotes this normalized
antiderivative, not a Fourier transform. The factor $y^{-2}$
comes from the hyperbolic measure. For example,
\begin{equation}\label{eq:ch7-intro-primitive-example}
 \varphi(y)=y^\alpha
 \quad\Longrightarrow\quad
 \widehat\varphi(y)=
 \begin{cases}
 \displaystyle\frac{y^{\alpha-1}}{\alpha-1},
      &\alpha\ne1,\\[6pt]
 \log y,&\alpha=1.
 \end{cases}
\end{equation}
Fixing the additive constant is essential, since changing it
would change the renormalized integral by the opposite constant.

\begin{proposition}\label{prop:classical-cusp-finite-part}
Suppose $F$ has the expansion \eqref{eq:asymp}, and normalize
$\widehat\varphi$ by \eqref{eq:ch7-intro-primitive}. Then the limit
\[
 \lim_{U\to\infty}\left(\int_{\mathcal D_U}F\,d\mu
                                  -\widehat\varphi(U)\right)
\]
exists. For every $T>1$ it equals
\[
 \int_{\mathcal D_T}F\,d\mu
 +\int_{\mathcal D\setminus\mathcal D_T}(F-\varphi)\,d\mu
 -\widehat\varphi(T).
\]
This expression is independent of $T$. If all
$\operatorname{Re}\alpha_j<1$, it agrees with the ordinary integral
of $F$.
\end{proposition}
\begin{proof}
For $U>T>1$, split the truncated domain at height $T$ and write
$F=(F-\varphi)+\varphi$ in the intervening cusp strip. This gives
the exact identity
\begin{equation}\label{eq:ch7-intro-finite-strip}
 \begin{aligned}
 \int_{\mathcal D_U}F\,d\mu
 ={}&
 \int_{\mathcal D_T}F\,d\mu
 +\int_{\mathcal D_U\setminus\mathcal D_T}
                         (F-\varphi)\,d\mu\\
 &+\int_{\mathcal D_U\setminus\mathcal D_T}
                         \varphi(y)\,d\mu.
 \end{aligned}
\end{equation}
Since the strip has width one, the last integral equals
\begin{equation}\label{eq:ch7-intro-endpoints}
 \int_{\mathcal D_U\setminus\mathcal D_T}\varphi(y)\,d\mu
 =
 \int_T^U\frac{\varphi(y)}{y^2}\,dy
 =
 \widehat\varphi(U)-\widehat\varphi(T).
\end{equation}
Consequently,
\begin{equation}\label{eq:ch7-intro-cutoff-identity}
 \begin{aligned}
 \int_{\mathcal D_U}F\,d\mu-\widehat\varphi(U)
 ={}&
 \int_{\mathcal D_T}F\,d\mu
 +\int_{\mathcal D_U\setminus\mathcal D_T}
                         (F-\varphi)\,d\mu\\
 &-\widehat\varphi(T).
 \end{aligned}
\end{equation}
The rapid decay of $F-\varphi$ makes its integral over the
infinite cusp absolutely convergent. We may therefore let
$U\to\infty$ and define
\begin{equation}\label{eq:ch7-intro-renormalized}
 \begin{aligned}
 \operatorname{R.N.}\!\int_{\Gamma\backslash\mathbb H}F\,d\mu
 &:=
 \lim_{U\to\infty}
 \left(
 \int_{\mathcal D_U}F\,d\mu-\widehat\varphi(U)
 \right)\\
 &=
 \int_{\mathcal D_T}F\,d\mu
 +\int_{\mathcal D\setminus\mathcal D_T}
                         (F-\varphi)\,d\mu
 -\widehat\varphi(T).
 \end{aligned}
\end{equation}
The finite-cutoff identity proves that the last expression is
independent of $T$. It also explains the minus sign:
$\widehat\varphi(T)$ is the lower-endpoint contribution in
\eqref{eq:ch7-intro-endpoints}. If
$\widehat\varphi(U)\to0$, as happens when all
$\operatorname{Re}\alpha_j<1$, this prescription agrees with
the ordinary convergent integral.

\end{proof}

The usefulness of this construction is that it retains the
unfolding underlying the Rankin--Selberg method. Let
$E(w,s)$ be the standard uncompleted Eisenstein series for
$\Gamma$, and put
\begin{equation}\label{eq:ch7-intro-constant-term}
 a_0(y)=\int_0^1F(x+iy)\,dx.
\end{equation}
For rapidly decreasing $F$, ordinary unfolding gives, initially
in a right half-plane,
\begin{equation}\label{eq:ch7-intro-ordinary-unfolding}
 \int_{\Gamma\backslash\mathbb H}
          F(w)E(w,s)\,d\mu(w)
 =
 \int_0^\infty a_0(y)y^{s-2}\,dy.
\end{equation}
When $F$ has the more general asymptotic expansion
\eqref{eq:asymp}, these integrals may diverge. The theorem on
page~426 of \cite{ZagierRS} supplies the corresponding
renormalized identity
\begin{equation}\label{eq:ch7-intro-zagier-unfolding}
 R(F;s)
 :=
 \operatorname{R.N.}\!\int_0^\infty
                    a_0(y)y^{s-2}\,dy
 =
 \operatorname{R.N.}\!\int_{\Gamma\backslash\mathbb H}
                    F(w)E(w,s)\,d\mu(w),
\end{equation}
initially for $\operatorname{Re}s$ sufficiently large, and
subsequently by meromorphic continuation. The renormalization
on the right uses the cusp expansion of the product
$F(w)E(w,s)$. The Eisenstein-series expression yields the
meromorphic continuation of $R(F;s)$ and, after the usual
completion, its symmetry under $s\mapsto1-s$. Residues recover
the renormalized average of $F$, with correction terms at
certain exceptional cusp exponents. In the present adelic setting, the analogous construction is
carried out at each fixed radial parameter $z>0$.

\subsection*{Acknowledgment}
The authors thank Robert Hough, Takashi Taniguchi, and Frank Thorne for extensive discussions. The first author was partially supported by a research grant from Texas Christian University, and the second by a collaboration grant from the Simons Foundation. The treatment of general compact types, the joint meromorphic continuation arguments, and the preparation and review of this article were developed with assistance from ChatGPT. Access to ChatGPT 6 Astra was provided through OpenAI's ChatGPT for Academic Researchers program via Julius Ross, whom the authors thank.

\section{Measures and lattice estimates}\label{sec:radial}
\subsection{Measures}\label{subsec:measures}
We use normalized local absolute values, so that $|z|_{\C}=z\bar z$. Fix a nontrivial continuous additive character $\psi_k$ of $\A/k$; the letter $\psi$ without a subscript is reserved for the smoothing function in $\mathscr E(\psi,w;h)$. Additive Haar measures on $k_v$ are self-dual for $\psi_{k,v}$, and their product is denoted $du$; in particular $\vol(\A/k)=1$. At a finite place put $\vol(\OO_v^\times,d^\times t)=1$. At an archimedean place $K$, write $t=r\omega$, with $r>0$ and $\omega\in\Omega_K$, where $\Omega_{\R}=\{\pm1\}$ and $\Omega_{\C}=S^1$. We take
\[
 d^\times t=e_K\frac{dr}{r}\,d\omega,
 \qquad e_{\R}=1,\quad e_{\C}=2,
\]
where $d\omega$ has total mass one. The residue constants below use precisely these normalizations.

We use Tate's decomposition of the idele norm and Haar measure \cite[Section~4]{Tate}, with the following choice of section. Write $\A^1=\ker(|\cdot|_{\A}:\A^\times\to\R_{>0})$. For $z>0$, let $\underline z$ be the idele with finite components $1$ and every archimedean component $z^{1/n}$. Then
\begin{equation}\label{eq:section}
 |\underline z|_{\A}=z,\qquad
 \underline{z_1z_2}=\underline z_1\underline z_2.
\end{equation}
We define a multiplicative Haar measure by
\[
 d^\times t=\frac{dz}{z}\,d^\times t^1,
 \qquad t=\underline z\,t^1,\quad t^1\in\A^1, \quad z > 0.
\]
The quotient measures by $k$, $k^\times$, and $G_k$ are always those obtained from counting measure on these discrete groups. We do not separately renormalize the inherited measure on $\A^1/k^\times$.

Put
\[
 n(c)=\begin{pmatrix}1&0\\c&1\end{pmatrix},\quad
 a(t_1,t_2)=\diag(t_1,t_2),\quad a(t)=a(t,t^{-1}),\quad
 U_\A=\prod_v U_v,
\]
where $U_v$ is a standard maximal compact subgroup. Equip $U_\A$ with its probability Haar measure $du$.  We specify measure normalizations explicitly as in Wright \cite{Wright}. We fix the Haar measure on $G_\A$ by the right-Iwasawa formula
\begin{equation}\label{eq:haar}
 \int_{G_\A} f(g)\,dg
 =\frac12\int_{U_\A}\int_{(\A^\times)^2}\int_{\A}
 f\bigl(u n(c)a(a,d)\bigr)\,dc\,d^\times a\,d^\times d\,du.
\end{equation}

Let $A$ be the diagonal group, $N$ the lower unipotent group, and $B=AN$. Put $G_\A^1=\{g:|\det g|_{\A}=1\}$, $B_\A^1=B_\A\cap G_\A^1$, and $B_\A^0=\{a(t_1,t_2)n(c):t_1,t_2\in\A^1\}$. Write $d(t)=a(t,t)$. In right Iwasawa coordinates $g=u a(t_1,t_2)n(c)$, set $t(g)=|t_1/t_2|_\A^{1/2}$ and $\chi_P(g)=t(g)^2$. The factor $1/2$ in \eqref{eq:haar} gives
\begin{equation}\label{eq:radial-measure}
 dg=\frac{dz}{z}\,dg^1,\qquad g=d(\underline z)g^1.
\end{equation}
Indeed, write $a=\underline z\,\underline\lambda\,t_1$ and $d=\underline z\,\underline\lambda^{-1}t_2$, with $t_i\in\A^1$. The determinant of the map
\[
 (\log z,\log\lambda)\longmapsto
 (\log|a|_{\A},\log|d|_{\A})
\]
is $-2$. Thus one half of the product measure on the two positive torus coordinates is $dz/z\,d\lambda/\lambda$, proving \eqref{eq:radial-measure}. In particular $|\det d(\underline z)|_{\A}=z^2$.

\subsection{Lattice bounds}
For $z>0$ put
\begin{equation}\label{eq:scalar-Schwartz-dilation}
 \Phi^z(x)=\Phi(d(\underline z)x).
\end{equation}
Here we record the lattice
estimates used throughout the paper.  All norms in the next two lemmas
are Euclidean norms.

\begin{lemma}\label{lem:uniform-lattice}
Let $L\subset\R^d$ be a full lattice.  If $Q>d$, then
\begin{equation}\label{eq:uniform-translated-lattice}
 \sup_{b\in\R^d}\sum_{\ell\in L}
       (1+R\|\ell+b\|)^{-Q}\ll_{L,Q}1+R^{-d}
       \qquad(R>0).
\end{equation}
If $A\ge0$ and $Q>d+A$, then
\begin{equation}\label{eq:uniform-nonzero-lattice}
 \sum_{0\ne\ell\in L}(1+R\|\ell\|)^{-Q}
       \ll_{L,Q,A} R^{-d}(1+R)^{-A}.
\end{equation}
Consequently, for every $f\in\mathscr S(\R^d)$ these bounds apply to
$\sum_{\ell\in L}|f(R(\ell+b))|$ and
$\sum_{0\ne\ell\in L}|f(R\ell)|$, respectively, with the right-hand
side multiplied by the Schwartz seminorm
\begin{equation}\label{eq:weighted-Schwartz-seminorm}
 p_Q(f)=\sup_y(1+\|y\|)^Q|f(y)|.
\end{equation}
They hold uniformly for bounded subsets of Schwartz space and after
replacing $f$ by any fixed continuous differential operator with
polynomial coefficients applied to $f$.
\end{lemma}

\begin{proof}
See \cite[Lemma~2.1]{LTcuspidal}. The statement is recalled here with the same lattice and Schwartz-seminorm conventions.
\end{proof}

\begin{lemma}\label{lem:triangular-lattice}
Let $L\subset E=\R^d$ be a full lattice and let $m\ge1$.  For
$R_1,\ldots,R_m>0$ and arbitrary real linear maps $B_{ih}:E\to E$
with $h<i$, put
\begin{equation}\label{eq:triangular-scaled-variables}
 y_i=R_i\left(x_i+\sum_{h<i}B_{ih}x_h\right),
       \qquad x=(x_1,\ldots,x_m)\in L^m.
\end{equation}
For $f\in\mathscr S(E^m)$ one has
\begin{equation}\label{eq:triangular-full-lattice}
 \sum_{x\in L^m}|f(y_1,\ldots,y_m)|
 \ll p(f)\prod_{i=1}^m(1+R_i^{-d}).
\end{equation}
For $1\le j\le m$ and $A\ge0$, one also has
\begin{equation}\label{eq:triangular-nonzero-stratum}
 \sum_{\substack{x\in L^m\\x_1=\cdots=x_{j-1}=0,\ x_j\ne0}}
       |f(y_1,\ldots,y_m)|
 \ll p_A(f) R_j^{-d}(1+R_j)^{-A}
                  \prod_{i>j}(1+R_i^{-d}).
\end{equation}
Here $p$ and $p_A$ are continuous Schwartz seminorms, and the constants
are independent of all the maps $B_{ih}$ and all the $R_i$.
The same assertions hold uniformly after replacing $f(y)$ by $f(Cy)$,
where $C$ ranges over a compact subset of $\GL(E^m)$.
\end{lemma}

\begin{proof}
Use the product Euclidean norm
$\|y\|^2=\sum_{i=1}^m\|y_i\|^2$ on $E^m$.  For $Q>0$ we have
\[
 \prod_{i=1}^m(1+\|y_i\|)^Q
 \le (1+\|y\|)^{mQ}.
\]
It follows that
\begin{equation}\label{eq:Schwartz-product-majorant}
 |f(y_1,\ldots,y_m)|
 \le C_{m,Q}p_{mQ}(f)
             \prod_{i=1}^m(1+\|y_i\|)^{-Q}.
\end{equation}
For the displayed product Euclidean norm one may take $C_{m,Q}=1$;
any other fixed choice of Euclidean norm changes only this constant.
In particular, the seminorm on the right is a seminorm of $f$ itself.

Choose $Q>d$, and, for $1\le i\le m$, put
\[
 \begin{split}
 b_i(x_1,\ldots,x_{i-1})
     &=\sum_{h<i}B_{ih}x_h,\\
 w_i(x_1,\ldots,x_i)
     &=\bigl(1+R_i\|x_i+b_i(x_1,\ldots,x_{i-1})\|\bigr)^{-Q}.
 \end{split}
\]
The empty sum defining $b_1$ is zero.  For every fixed choice of the
earlier variables, \eqref{eq:uniform-translated-lattice} gives
\[
 \sum_{x_i\in L}w_i(x_1,\ldots,x_i)
 \le C_{L,Q}(1+R_i^{-d}).
\]
This bound holds for every value of the translate $b_i$ and therefore
does not involve any norm of any $B_{ih}$.

The triangular order permits us to use these bounds successively.
Since $w_i$ depends only on $x_1,\ldots,x_i$, summing first in $x_m$
gives
\[
 \sum_{x\in L^m}\prod_{i=1}^m w_i
 \le C_{L,Q}(1+R_m^{-d})
       \sum_{(x_1,\ldots,x_{m-1})\in L^{m-1}}
                    \prod_{i=1}^{m-1}w_i.
\]
Iterating in the order $x_{m-1},\ldots,x_1$ yields
\[
 \sum_{x\in L^m}\prod_{i=1}^m w_i
 \le C_{L,Q}^{\,m}\prod_{i=1}^m(1+R_i^{-d}).
\]
Tonelli's theorem justifies this order of summation because every
term is nonnegative; the resulting finite upper bound also proves
convergence.  Combining it with
\eqref{eq:Schwartz-product-majorant} proves
\eqref{eq:triangular-full-lattice}, with $p(f)=p_{mQ}(f)$.

For the second assertion, fix $j$ and $A\ge0$, and now choose
$Q>d+A$.  On the specified stratum,
$x_1=\cdots=x_{j-1}=0$, so $y_i=0$ for $i<j$ and
$y_j=R_jx_j$.  Sum first over $x_m,\ldots,x_{j+1}$, exactly as above.
The factors with $i<j$ in the product majorant are all equal to one,
and the remaining sum is bounded by
\[
 C_{m,Q}p_{mQ}(f)C_{L,Q}^{\,m-j}
 \prod_{i>j}(1+R_i^{-d})
 \sum_{0\ne x_j\in L}(1+R_j\|x_j\|)^{-Q}.
\]
The last sum is \emph{untranslated}.  Applying
\eqref{eq:uniform-nonzero-lattice} to it gives
\eqref{eq:triangular-nonzero-stratum}, with
$p_A(f)=p_{mQ}(f)$.  This also covers $j=1$ and $j=m$:
the conditions on earlier variables and the product over later
variables are then interpreted as empty when appropriate.
The vanishing of the earlier variables is precisely what permits
the nonzero-lattice estimate at the $j$th step.

For completeness, let $\mathcal C\subset\GL(E^m)$ be compact and set
\[
 K=\max\left(1,\sup_{C\in\mathcal C}\|C^{-1}\|_{\rm op}\right)<\infty.
\]
For every $s\ge0$ and $C\in\mathcal C$, changing variables $v=Cy$
in the defining supremum gives
\[
 \begin{split}
 p_s(f\circ C)
 &=\sup_v(1+\|C^{-1}v\|)^s|f(v)|\\
 &\le K^s p_s(f).
 \end{split}
\]
Apply the two estimates just proved to $f\circ C$ and use this
inequality with $s=mQ$.  This proves the final uniformity assertion.
The constants may depend on $L,d,m,Q,A$, the fixed norms, and
$\mathcal C$, but not on the $R_i$ or the maps $B_{ih}$.
\end{proof}

\begin{remark}\label{rem:adelic-lattice-uniformity}
In this paper we will be applying these lemmas to adelic sums.  To make this application precise, fix a compact support $K_f$ and a
compact open additive subgroup $L_f$ under which the finite-adelic
functions are invariant.  The corresponding space of finite-adelic
test functions is finite dimensional.  It is enough to treat a fixed
basis of that space, with its finitely many archimedean Schwartz
coefficients.  If the finite group factors range over a compact set,
their inverse images of $K_f$ lie in a single compact subset of
$\A_f^m$, hence in $\mathfrak a^{-1}\widehat{\OO}_k^m$ for some
nonzero integral ideal $\mathfrak a$.  Every contributing rational
point therefore lies in the full lattice
$(\mathfrak a^{-1})^m\subset k_\infty^m$.
The constants are controlled by finitely many seminorms of those
archimedean coefficients.  
\end{remark}
\begin{remark} Here is an elementary observation that will be used later on. We have 
\begin{equation}\label{eq:log-absorption}
 (1+|\log R|)^q\ll_{q,\varepsilon}R^{\varepsilon}+R^{-\varepsilon}
 \qquad(R>0,\ q\ge0,\ \varepsilon>0).
\end{equation}
\end{remark}

\subsection{Reduction theory}
Reduction theory for $\GL_2$ \cite[Section~10]{JL}, with the compactness of $\A^1/k^\times$ and $\A/k$ \cite[Section~4]{Tate}, provides compact representative sets $C_1\subset\A^1$, $C_N\subset\A$, and $\lambda_0>0$ such that
\begin{equation}\label{eq:explicit-Siegel-set}
\mathfrak S=\{u a(\underline\lambda)a(t_1,t_2)n(c):u\in U_\A,\ \lambda\ge\lambda_0,\ t_1,t_2\in C_1,\ c\in C_N\}
\end{equation}
maps onto $G_\A^1/G_k$. The integral of a nonnegative quotient function is bounded by a fixed multiple of its integral on this cover with density
\begin{equation}\label{eq:explicit-Siegel-density}
du\,\lambda^{-2}\frac{d\lambda}{\lambda}\,d^\times t_1d^\times t_2dc.
\end{equation}
Only compact transverse representative sets and this domination property are used before the exact high-cusp coordinates are established in Section~\ref{subsec:eis-finite-part}.

\section{The Eisenstein integral and absolute convergence}\label{subsec:eis-family}\label{sec:eisenstein}
Retain the action of Section~\ref{sec:statement} and the measures and right-Iwasawa conventions of Section~\ref{subsec:measures}.  In this section write $G^1=G_\A^1$ and
$U=U_\A$.  Thus
\begin{equation}\label{eq:Iwasawa}
 g=u a(\underline\lambda)a(t_1,t_2)n(c),\qquad
 t_i\in\A^1,\qquad t(g)=\lambda
\end{equation}
on $G^1$, and $t$ is unchanged by scalar multiplication.  Put
$t_L(h)=t(h^{-1})$. Let $\eta_1,\eta_2$ be unitary Hecke
characters of $k^\times\backslash\A^\times$, trivial on the chosen
positive radial section $\underline r$, and put $\omega=\eta_1\eta_2$.
The restriction to the radial section is a normalization, not a
restriction on the unitary central characters covered by the theorem.
Indeed any unitary central character has a unique expression
$\omega_0|\cdot|_\A^{i\kappa}$ with $\omega_0(\underline r)=1$ and
$\kappa\in\R$. Twisting a normalized family by
$|\det|_\A^{i\kappa/2}$ restores it and replaces $s$ throughout by
$s-i\kappa/2$. The two remaining opposite norm exponents in an
inducing pair are absorbed into $\tau$. Thus it suffices to prove
all assertions with the displayed norm-normalized characters.
No unramifiedness assumption is imposed. Let $f_\tau$ be a smooth
$U$-finite flat section satisfying
\begin{equation}\label{eq:eis-general-induction}
 f_\tau\bigl(a(a,d)n(c)h\bigr)
 =\eta_1(a)\eta_2(d)|d/a|_\A^{1/2+\tau} f_\tau(h).
\end{equation}
Flat means that $f_\tau|_U$ is independent of $\tau$. The section
has fixed finite level and is unramified at almost all places;
finite sums of such sections are allowed. Define
\begin{equation}\label{eq:eis-general-family}
 E(\tau,h)=E(f_\tau,h)
   =\sum_{\gamma\in B_k\backslash G_k}f_\tau(\gamma h),
                    \qquad \Re\tau>\tfrac12.
\end{equation}
This is left $G_k$-invariant and has central character $\omega$.
The two unitary characters are independent. In particular neither is
required to be unramified, or to be the inverse of the other. We
consider arbitrary smooth $U$-finite sections, without requiring a
spherical vector.

Write $M_\tau$ for the standard global intertwining operator,
initially given in its convergence half-plane by
\begin{equation}\label{eq:eis-global-intertwiner}
 (M_\tau f_\tau)(h)
   =\int_\A f_\tau(w n(v)h)\,dv,\qquad
 w=\begin{pmatrix}0&1\\-1&0\end{pmatrix}.
\end{equation}
It carries the section for $(\eta_1,\eta_2,\tau)$ to the section for
$(\eta_2,\eta_1,-\tau)$, and
\begin{equation}\label{eq:eis-constant-definition}
 E_{\mathrm c}(\tau,h)
 :=\int_{\A/k}E(\tau,n(c)h)\,dc
 = f_\tau(h)+(M_\tau f_\tau)(h).
\end{equation}
The meromorphic continuation, functional equation, and constant
term in this generality are the standard $\GL_2$ Eisenstein
results; see \cite[Section~16]{JL} and
\cite[Sections~4.1.6--4.1.8]{MV}, with the opposite Borel
convention changed as in \eqref{eq:eis-general-induction}.
Set
\begin{equation}\label{eq:eis-branch-data}
 \begin{gathered}
 C=\vol(\A^1/k^\times),\qquad p_\epsilon=1+2\epsilon\tau,\\
 (\eta_{+,1},\eta_{+,2})=(\eta_1,\eta_2),\qquad
 (\eta_{-,1},\eta_{-,2})=(\eta_2,\eta_1),\\
 c_+(\tau,u)=f_\tau(u^{-1}),\qquad
 c_-(\tau,u)=(M_\tau f_\tau)(u^{-1}),\\
 b_\epsilon(\tau;u,t_1,t_2)
 =\eta_{\epsilon,1}(t_1)^{-1}\eta_{\epsilon,2}(t_2)^{-1}
                  c_\epsilon(\tau,u).
 \end{gathered}
\end{equation}
Then the right-Iwasawa constant term is
\begin{equation}\label{eq:eis-ct}
 E_{\mathrm c}(\tau,g^{-1})
       =\sum_{\epsilon=\pm1}\lambda^{p_\epsilon}
                            b_\epsilon(\tau;u,t_1,t_2).
\end{equation}
Thus the scattering operator acts on the section: it cannot in
general be replaced by a scalar scattering coefficient.

For the trivial pair and the section equal to $1$ on $U$,
$c_+=1$ and $c_-=m(\tau)$, where
\begin{equation}\label{eq:eis-spherical-scattering}
 \begin{gathered}
 m(\tau)=\frac{\Lambda_k(2\tau)}{\Lambda_k(1+2\tau)},\\
 \Lambda_k(v)=|D_k|^{v/2}
 (\pi^{-v/2}\Gamma(v/2))^{r_1}
 (2(2\pi)^{-v}\Gamma(v))^{r_2}\zeta_k(v).
 \end{gathered}
\end{equation}
We retain $m$ and $\Lambda=\Lambda_k$ for this spherical auxiliary
family, particularly in Appendix~\ref{subsec:eis-smoothed}.
Put $E_{\mathrm n}=E-E_{\mathrm c}$; this difference is defined
on $N_k\backslash G_\A$, and is not asserted to be automorphic
under all of $G_k$. Regularity has the meaning fixed in the introduction.
The compact unipotent integral \eqref{eq:eis-constant-definition}
shows that regularity of $E$ implies regularity of its constant
term and hence of $M_\tau f_\tau$.
For a regular parameter and $\gamma=|\Re\tau|$,
\begin{equation}\label{eq:eis-bound}
 |E(\tau,g^{-1})|\ll_{\tau,f}\lambda^{1+2\gamma}
                         \qquad(\lambda\ge1)
\end{equation}
on a fixed Siegel set. The nonconstant part decreases faster than
any power of $\lambda$. These estimates are uniform on compact
regular spectral sets after a harmless enlargement of the exponent;
local pole-clearing gives the corresponding assertion near a pole.
They follow from the Fourier expansion and the fixed-vector
Whittaker estimates used below.

When $\eta_1,\eta_2$ and the section are fixed, we abbreviate $E_\tau=E(\tau)$ and suppress these data in subscripts: for example, $R_{\Phi,E_\tau}=R_{\Phi,\tau}$. The same convention applies to $D_\epsilon$, $H$, $P_J$, and the singular Mellin transforms.

Values, residues, and higher Laurent coefficients at singular
parameters are addressed in Section~\ref{subsec:eis-exceptional}.

We shall use two measure consequences of \eqref{eq:haar}:
\begin{equation}\label{eq:eis-measures}
 \begin{aligned}
 dg&=\frac{dz}{z}\,dg^1,\qquad g=d(\underline z)g^1,\\
 dg^1&=du\,\lambda^{-2}\frac{d\lambda}{\lambda}
                   d^\times t_1d^\times t_2dc, \quad g \text{ as in \eqref{eq:Iwasawa}}.
 \end{aligned}
\end{equation}
The logarithmic Jacobian of $(z,\lambda)\mapsto(z\lambda,z\lambda^{-1})$
is $2$, canceling the factor $1/2$ in \eqref{eq:haar}; moving the
unipotent factor to the right contributes $\lambda^{-2}$.
On $G^1/B_k$, the variables $t_i$ range over $\A^1/k^\times$ and
$c$ over $\A/k$.  

\begin{lemma}\label{lem:eis-initial}
For $\Phi\in\cS(V_\A)$ and a regular $\tau$, the integral
\[
 Z(s,\Phi,E(\tau))=\int_{G_\A/G_k}|\det g|_\A^s
       E(\tau,g^{-1})\sum_{x\in V'_k}\Phi(gx)\,dg
\]
converges absolutely for $\Re s>5/2+\gamma$.
Its upper radial part, obtained by restricting $z\ge1$, is entire in $s$.
\end{lemma}
\begin{proof}
Write $g=d(\underline z)h$, with $h\in G^1$, and put
\begin{equation}\label{eq:nonneg-theta}
 \Theta_\Phi(z,h)=\sum_{x\in V'_k}|\Phi(d(\underline z)hx)|,
 \qquad
 A_\tau(z)=\int_{G^1/G_k}|E(\tau,h^{-1})|
                                      \Theta_\Phi(z,h)\,dh.
\end{equation}
Since $|\omega|=1$, Tonelli's theorem reduces
absolute convergence to
\[
 \int_0^\infty z^{2\Re s}A_\tau(z)\,\frac{dz}{z}<\infty.
\]
Choose $\Gamma>|\Re\tau|$. On a Siegel set the Eisenstein
factor is $O(\lambda^{1+2\Gamma})$, with transverse variables
in the fixed compact representative sets of Section~\ref{sec:radial}.

For completeness, we record the consequence of the common lattice
estimate that will also be used for the singular theta series.
Set
\begin{equation}\label{eq:eis-cubic-scales}
 A_1=z\lambda^3,\qquad A_2=z\lambda,\qquad
 A_3=z\lambda^{-1},\qquad A_4=z\lambda^{-3}.
\end{equation}
The finite support condition confines all contributing rational
coefficients to a fixed fractional ideal. The lower-unipotent
coordinate map is triangular, with diagonal entries one.
Lemma~\ref{lem:triangular-lattice}, applied with
$R_i=A_i^{1/n}$, therefore gives, for every $M>0$,
\begin{equation}\label{eq:eis-first-two-strata}
\begin{aligned}
 \sum_{x_1\ne0}|\Phi(d(\underline z)hx)|
 &\ll_M\prod_{i=1}^4(1+A_i^{-1})(1+A_1)^{-M},\\
 \sum_{x_1=0,\,x_2\ne0}|\Phi(d(\underline z)hx)|
 &\ll_M\prod_{i=2}^4(1+A_i^{-1})(1+A_2)^{-M}.
\end{aligned}
\end{equation}
Here the sums range over all remaining rational coefficients, and
$n=[k:\mathbb Q]$. Compact invertible factors are absorbed in the
Schwartz seminorm. In passing from the lemma to this display, choose
its decay order at least $nM$ and use
$(1+A^{1/n})^n\ge1+A$.

A nonsingular cubic has $(x_1,x_2)\ne(0,0)$. Since $A_i\le A_1$
and $A_2A_3A_4=z^3\lambda^{-3}$ for $\lambda\ge1$, the two
bounds imply, for every $N>0$,
\begin{equation}\label{eq:eis-theta-bound}
 \Theta_\Phi(z,h)\ll_N
 z^{-4}(1+z\lambda^3)^{-N}
 +z^{-3}\lambda^3(1+z\lambda)^{-N}.
\end{equation}
Indeed, write $1+A_i^{-1}=A_i^{-1}(1+A_i)$, and absorb
the resulting three or four growing factors by increasing $M$.
This estimate holds for all $z>0$ and $\lambda\ge1$.

For $0<z\le1$, multiplication by the Eisenstein bound and the
cusp density gives the two integrals
\begin{equation}\label{eq:eis-radial-majorants}
\begin{aligned}
 B_1(z)&=z^{-4}\int_1^\infty
       \lambda^{-1+2\Gamma}(1+z\lambda^3)^{-N}\frac{d\lambda}{\lambda},\\
 B_2(z)&=z^{-3}\int_1^\infty
       \lambda^{2+2\Gamma}(1+z\lambda)^{-N}\frac{d\lambda}{\lambda}.
\end{aligned}
\end{equation}
For $N>2+2\Gamma$, the substitution $v=z\lambda$ gives
$B_2(z)\ll z^{-5-2\Gamma}$. Using
$\lambda^{-1+2\Gamma}\le\lambda^{2+2\Gamma}$ and then
$v=z\lambda^3$ gives
$B_1(z)\ll z^{-4-(2+2\Gamma)/3}\le z^{-5-2\Gamma}$.
The bounded-height region contributes $O(z^{-4})$ by
Lemma~\ref{lem:uniform-lattice}. Hence
\begin{equation}\label{eq:eis-small-radius-average}
 A_\tau(z)\ll z^{-5-2\Gamma}\qquad(0<z\le1).
\end{equation}

For $z,\lambda\ge1$, the counting factors in
\eqref{eq:eis-first-two-strata} are $O(\lambda^4)$ in total,
whereas the nonzero-coordinate factor gives arbitrary decay in
$z\lambda$. Thus, for every $P,Q>0$,
\begin{equation}\label{eq:eis-large-radius-theta}
 \Theta_\Phi(z,h)\ll_{P,Q}z^{-P}\lambda^{-Q}.
\end{equation}
The bounded-height region has the same arbitrary decay in $z$, by
the nonzero-vector estimate in Lemma~\ref{lem:uniform-lattice}.
Taking $Q>2\Gamma-1$ therefore proves
\begin{equation}\label{eq:eis-large-radius-average}
 A_\tau(z)\ll_P z^{-P}\qquad(z\ge1)
\end{equation}
for every $P>0$. Given $\Re s>5/2+|\Re\tau|$, choose
$|\Re\tau|<\Gamma<\Re s-5/2$ and then $P>2\Re s$.
The two radial integrals converge by
\eqref{eq:eis-small-radius-average} and
\eqref{eq:eis-large-radius-average}.
The latter estimate is uniform on compact sets of $s$ and remains
integrable after multiplication by any fixed power of $\log z$;
this proves the assertion about the upper radial part.
\end{proof}

\begin{corollary}\label{cor:eis-upper-entire}
The integral 
\begin{equation}\label{eq:eis-upper-integral}
 Z_+(s,\Phi,E(\tau))=
 \int_1^\infty\int_{G^1/G_k}z^{2s}E(\tau,h^{-1})
       \sum_{x\in V'_k}\Phi(d(\underline z)hx)\,dh\,\frac{dz}{z}
\end{equation}
is entire in $s$. On compact sets of regular spectral parameters,
the convergence and all $s$-derivatives are locally uniform.
\end{corollary}


\section{A common cusp finite part and its unfolding}\label{subsec:eis-finite-part}
Until \S\ref{subsec:eis-exceptional}, the parameter is generic
in the sense fixed in the introduction.

For $X=S$ or $V'$, use the notation
\[
 \Theta_{X,\Phi}(z,g)=\sum_{x\in X(k)}\Phi(\underline z gx).
\]
In particular, $\Theta_{S,\Phi}$ includes the zero orbit.

Put $V_0=\{x_1=x_2=0\}$ and fix $\Phi\in\cS(V(\A))$.
A nonzero singular cubic is, up to a nonzero scalar, either a cube of
a linear form or a product of a squared linear form and a different
linear form.  Its repeated
root is unique and is defined over $k$: Galois automorphisms preserve
root multiplicity.  The stabilizer of the line defining the repeated
root of $v^2(\alpha u+\beta v)$ is $B_k$.
Consequently distinct translates of $V_0(k)\setminus\{0\}$ under
$G_k/B_k$ are disjoint, and their union is $S_k\setminus\{0\}$.
Thus
\begin{equation}\label{eq:eis-function}
 \Theta_{S,\Phi}(z,g)=\Phi(0)+P f_z(g),\qquad
 P f_z(g)=\sum_{\gamma\in G_k/B_k}f_z(g\gamma),
\end{equation}
where
\[
 f_z(g)=\sum_{(\alpha,\beta)\in k^2\setminus\{(0,0)\}}
          \Phi(\underline z g(0,0,\alpha,\beta)).
\]
For every fixed $z,g$ these sums
converge absolutely.

Use the Euclidean column norm over $\R$, its square over $\C$, and
the maximum norm at finite places.  These norms scale by the normalized
local absolute value.  If $v_\gamma$ is the second column of
$\gamma\in G_k$, then
\[
 t(g\gamma)^{-1}=\prod_v\|g_vv_\gamma\|_v\qquad(g\in G^1).
\]
Indeed, the second column of $u a(A,D)n(c)$ is $D u e_2$,
and $|A D|_\A=1$.
If $\gamma B_k\ne\gamma'B_k$, the two rational second columns
are linearly independent.  At every place,
\[
 |\det(g_vv_\gamma,g_vv_{\gamma'})|_v
       \le\|g_vv_\gamma\|_v\|g_vv_{\gamma'}\|_v.
\]
Multiplying and using $|\det g|_\A=1$ and the product formula for
the nonzero rational determinant gives
\begin{equation}\label{eq:eis-separated-cusps}
 t(g\gamma)t(g\gamma')\le1.
\end{equation}
For $T>1$, let $\mathcal F_T$ consist of the classes for which
$t(g\gamma)\le T$ for every rational cusp $\gamma B_k$.
The Siegel cover shows that $\mathcal F_T$ is compact once
$T>\max(1,\lambda_0)$.

Write
\[
 \begin{aligned}
 A_\Phi(u)&=\int_{\A^2}\Phi(u(0,0,x,v))\,dx\,dv,\\
 \widetilde\Phi(x)&=\int_U\Phi(ux)du,
 \qquad A_\Phi=\int_U A_\Phi(u)du.
 \end{aligned}
\]
For $g=u a(\underline\lambda)a(t_1,t_2)n(c)$, the linear map on
$V_0$ is
\begin{equation}\label{ex:plane-map}
 L_{z,\lambda,t,c}(0, 0, \alpha,\beta)=
 \left(0, 0, \underline z\underline\lambda^{-1}t_2\alpha,
 \underline z\underline\lambda^{-3}t_2^2t_1^{-1}
                                                (c\alpha+\beta)\right).
\end{equation}

\begin{lemma}\label{lem:eis-function-asymptotic}
For $g=u a(\underline\lambda)a(t_1,t_2)n(c)$ and $\lambda\ge1$,
\begin{equation}\label{eq:eis-function-asymp}
 f_z(g)=z^{-2}\lambda^4A_\Phi(u)-\Phi(0)
                           +O_{N,z}(\lambda^{-N})
\end{equation}
for every $N\ge0$. The estimate is uniform in the fixed compact
representatives for $u,t_1,t_2,c$, and for $z$ in compact subsets
of $(0,\infty)$.
\end{lemma}
\begin{proof}
Put $\Psi_u(x,v)=\Phi(u(0,0,x,v))$ and write the plane map
\eqref{ex:plane-map} as
\begin{equation}\label{eq:eis-plane-matrix}
 L=\begin{pmatrix}a&0\\bc&b\end{pmatrix},\qquad
 a=\underline z\underline\lambda^{-1}t_2,\qquad
 b=\underline z\underline\lambda^{-3}t_2^2t_1^{-1}.
\end{equation}
Its additive determinant has norm $z^2\lambda^{-4}$.
Adelic Poisson summation, with the self-dual measures, gives
\begin{equation}\label{eq:eis-function-poisson}
 f_z(g)+\Phi(0)
 =z^{-2}\lambda^4
          \sum_{\eta\in k^2}\widehat\Psi_u(L^{-t}\eta).
\end{equation}
The zero frequency is $A_\Phi(u)$.
The family $\widehat\Psi_u$ is bounded in Schwartz space, with
uniformly compact finite support. Since the finite components of
$L$ and $L^{-1}$ range in compact sets, all contributing rational
frequencies lie in a fixed fractional-ideal lattice
$\Gamma\subset k_\infty^2\simeq\mathbb R^{2n}$.
At infinity, uniformly for $z$ in a fixed compact interval,
\begin{equation}\label{eq:eis-dual-plane-expansion}
 \|L_\infty^t\|_{\mathrm{op}}\ll\lambda^{-1/n},
 \qquad
 \|L_\infty^{-t}\eta\|\gg\lambda^{1/n}\|\eta\|.
\end{equation}
Apply the nonzero-vector estimate in Lemma~\ref{lem:uniform-lattice}
in dimension $2n$ with $R\asymp\lambda^{1/n}$. The sum over
$\eta\ne0$ in \eqref{eq:eis-function-poisson} is
$O_M(\lambda^{-M})$ for every $M$. Its prefactor is
$O_z(\lambda^4)$, which is absorbed by increasing $M$.
This proves \eqref{eq:eis-function-asymp}; the term $-\Phi(0)$
records the omitted zero vector in the function.
\end{proof}

\begin{lemma}\label{lem:eis-other-cusps}
With the same uniformity, for every $N\ge0$,
\begin{equation}\label{ex:other-cusps}
 P f_z(g)-f_z(g)=O_{N,z}(\lambda^{-N})
                       \qquad(\lambda\ge1).
\end{equation}
The estimate also holds for the sum of the absolute values of
all summands on the left.
\end{lemma}
\begin{proof}
The repeated-root decomposition identifies the difference with
\begin{equation}\label{eq:other-cusps-proof-identity}
 P f_z(g)-f_z(g)
 =\sum_{x\in S_k\setminus V_0(k)}\Phi(d(\underline z)gx).
\end{equation}
Every vector in this sum has $x_1\ne0$, or $x_1=0$ and
$x_2\ne0$. Consequently its absolute sum is bounded by the two
strata in \eqref{eq:eis-first-two-strata}; no discriminant
condition was used in deriving that estimate. For $z$ in a fixed
compact interval and $\lambda\ge1$, their counting factors are
$O(\lambda^4)$ and their respective decay factors are
$O_M(\lambda^{-3M})$ and $O_M(\lambda^{-M})$.
Choosing $M>N+4$ proves the assertion.
\end{proof}

\begin{lemma}\label{lem:eis-small-function}
For $0<\lambda\le1$ and every $N\ge0$,
\begin{equation}\label{ex:function-small}
 \sum_{(\alpha,\beta)\in k^2\setminus\{(0,0)\}}
 |\Phi(d(\underline z)g(0,0,\alpha,\beta))|
                               =O_{N,z}(\lambda^N),
\end{equation}
uniformly in the same compact parameters.
\end{lemma}
\begin{proof}
For the matrix in \eqref{eq:eis-plane-matrix},
$L^{-1}=\left(\begin{smallmatrix}a^{-1}&0\\-ca^{-1}&b^{-1}
\end{smallmatrix}\right)$, so
$\|L_\infty^{-1}\|_{\mathrm{op}}\ll_z\lambda^{1/n}$.
Thus $\|L_\infty\xi\|\gg_z\lambda^{-1/n}\|\xi\|$.
Finite support again confines the rational vectors to a fixed
lattice in dimension $2n$. Apply Lemma~\ref{lem:uniform-lattice}
with $R\asymp_z\lambda^{-1/n}$ to the absolute Schwartz
majorant. Its nonzero-vector estimate gives arbitrary positive
powers of $\lambda$.
\end{proof}

The zero term in \eqref{eq:eis-function} cancels the term $-\Phi(0)$ in
\eqref{eq:eis-function-asymp}.  Together with \eqref{ex:other-cusps}, this gives
\begin{proposition}\label{prop:eis-singular-theta-asymptotic}
For every $N\ge0$, uniformly in the same compact parameters,
\[
 \Theta_{S,\Phi}(z,g)=z^{-2}\lambda^4 A_\Phi(u)
                                      +O_{N,z}(\lambda^{-N}).
\]
\end{proposition}

Define the averages
\begin{equation}\label{eq:eis-transverse-coefficients}
 \begin{aligned}
 D_\epsilon(\Phi,\tau)
 &=\int_{U\times(\A^1/k^\times)^2}
 A_\Phi(u)b_\epsilon(\tau;u,t_1,t_2)
                    \,du\,d^\times t_1d^\times t_2,\\
 e_\epsilon(\tau)
 &=\int_{U\times(\A^1/k^\times)^2}
 b_\epsilon(\tau;u,t_1,t_2)
                    \,du\,d^\times t_1d^\times t_2.
 \end{aligned}
\end{equation}
These are the averages with the inherited Iwasawa measure.
Orthogonality of characters on the compact group $\A^1/k^\times$
gives
\begin{equation}\label{eq:eis-transverse-selection}
 \begin{aligned}
 D_\epsilon(\Phi,\tau)
 &=C^2\mathbf1_{\{\eta_1=\eta_2=1\}}
                  \int_U A_\Phi(u)c_\epsilon(\tau,u)\,du,\\
 e_\epsilon(\tau)
 &=C^2\mathbf1_{\{\eta_1=\eta_2=1\}}
                  \int_U c_\epsilon(\tau,u)\,du.
 \end{aligned}
\end{equation}

After multiplying the theta leading term by the Eisenstein
constant term and integrating the transverse variables, the cusp
density is
\[
 z^{-2}\sum_\epsilon D_\epsilon\lambda^{3+2\epsilon\tau}
                                      \frac{d\lambda}{\lambda}.
\]
Its normalized antiderivative is
\begin{equation}\label{eq:finite-part-explicit-subtraction}
 B_T(z,\tau)=z^{-2}\sum_{\epsilon=\pm1}
 D_\epsilon(\Phi,\tau)\frac{T^{q_\epsilon}}{q_\epsilon},
 \qquad q_\epsilon=3+2\epsilon\tau.
\end{equation}
Unless $\eta_1=\eta_2=1$, this subtraction vanishes identically.
The resulting convergence of geometric cusp truncations need not
be absolute: the character cancellation is performed over each
complete transverse quotient. In the trivial inducing case the
section-weighted compact averages in
\eqref{eq:eis-transverse-selection} must still be retained, even
when the section has no spherical component.

\begin{proposition}
\label{prop:eis-finite-part-exists}
Fix $\Phi\in\mathcal S(V_{\mathbb A})$ and $z>0$. Suppose that
the Eisenstein family $E(\tau)$ and its intertwining section
$M_\tau f_\tau$ are regular at the spectral parameter under consideration,
and that
\begin{equation}\label{eq:finite-part-nonresonance}
q_\epsilon:=3+2\epsilon\tau\ne0
\qquad(\epsilon=\pm1).
\end{equation}
Then the limit 
\begin{equation}\label{eq:eis-reg-subtraction}
 R^{\rm reg}_{\Phi,\tau}(z)=
 \lim_{T\to\infty}\left\{
 \int_{\mathcal F_T}E(\tau,g^{-1})\Theta_{S,\Phi}(z,g)\,dg^1
                                      -B_T(z,\tau)\right\}
\end{equation}
exists.

The convergence is locally uniform for $z>0$ and for regular spectral parameters satisfying
\eqref{eq:finite-part-nonresonance}. The resulting expression is holomorphic in $\tau$ when the Eisenstein family and its intertwining section are holomorphic and \(3\pm2\tau\ne0\) and
defines a locally meromorphic family in $\tau$. If $(\eta_1,\eta_2)\ne(1,1)$,
the same assertions hold without \eqref{eq:finite-part-nonresonance},
with $B_T=0$ at all regular parameters.
\end{proposition}

\begin{proof}
Let
\begin{equation}\label{eq:finite-part-truncated-pairing}
 \mathcal A_T(z,\tau)=\int_{\mathcal F_T}
 E(\tau,g^{-1})\Theta_{S,\Phi}(z,g)\,dg^1.
\end{equation}
Choose $T_0>1$ sufficiently large. By
\eqref{eq:eis-separated-cusps}, the projection $G^1/B_k\to G^1/G_k$
identifies $\{t(g)>T_0\}$ with the complement of
$\mathcal F_{T_0}$: every such class has a unique rational cusp
above $T_0$. The cusp coordinates and inherited measure are
\begin{equation}\label{eq:finite-part-cusp-coordinates}
 g=u a(\underline\lambda)a(t_1,t_2)n(c),\qquad
 dg^1=\lambda^{-2}\frac{d\lambda}{\lambda}\,
               du\,d^\times t_1d^\times t_2dc,
\end{equation}
where the transverse quotient has finite measure and is represented
in fixed compact sets.

Proposition~\ref{prop:eis-singular-theta-asymptotic} gives
\begin{equation}\label{eq:finite-part-theta-expansion}
 \Theta_{S,\Phi}(z,g)=z^{-2}\lambda^4A_\Phi(u)+r_{\Phi,z}(g),
 \qquad r_{\Phi,z}(g)=O_{M,z}(\lambda^{-M})
\end{equation}
for every $M$, locally uniformly in $z>0$. Define
\begin{equation}\label{eq:finite-part-integrand}
\begin{aligned}
 H_{\Phi,\tau}(z,g)&=E(\tau,g^{-1})\Theta_{S,\Phi}(z,g),\\
 \mathcal L_{\Phi,\tau}(z,g)&=z^{-2}A_\Phi(u)
       \sum_{\epsilon=\pm1}b_\epsilon(\tau;u,t_1,t_2)
                              \lambda^{5+2\epsilon\tau}.
\end{aligned}
\end{equation}
The exact remainder identity is
\begin{equation}\label{eq:finite-part-remainder-identity}
 \mathcal R:=H_{\Phi,\tau}-\mathcal L_{\Phi,\tau}
 =E(\tau,g^{-1})r_{\Phi,z}(g)
       +z^{-2}\lambda^4A_\Phi(u)E_{\mathrm n}(\tau,g^{-1}).
\end{equation}
The Eisenstein series has moderate growth and its nonconstant part
has rapid cusp decay. It follows that $\mathcal R$ decreases
faster than every power of $\lambda$, so its integral over
$\lambda>T_0$ is absolutely convergent.

Since $\operatorname{vol}(\A/k)=1$, integration over the complete
transverse quotient, followed by \eqref{eq:eis-transverse-coefficients}, gives
\begin{equation}\label{eq:finite-part-leading-shell}
 \int_{T_0<\lambda\le T}\mathcal L_{\Phi,\tau}(z,g)\,dg^1
                   =B_T(z,\tau)-B_{T_0}(z,\tau).
\end{equation}
All complex powers of positive real numbers use the real logarithm.
Splitting at $T_0$ yields the exact finite-cutoff identity
\begin{equation}\label{eq:finite-part-finite-cutoff-identity}
 \mathcal A_T-B_T
 =\int_{\mathcal F_{T_0}}H_{\Phi,\tau}\,dg^1-B_{T_0}
                    +\int_{T_0<\lambda\le T}\mathcal R\,dg^1.
\end{equation}
Taking the limit proves
\begin{equation}\label{eq:finite-part-convergent-representation}
 R^{\mathrm{reg}}_{\Phi,\tau}(z)
 =\int_{\mathcal F_{T_0}}H_{\Phi,\tau}(z,g)\,dg^1-B_{T_0}(z,\tau)
                       +\int_{\lambda>T_0}\mathcal R(z,\tau,g)\,dg^1.
\end{equation}
In particular, this expression is independent of the auxiliary $T_0$.

We give the uniformity argument because it also justifies spectral
continuation. Let $K\subset\operatorname{int}K_1$ be compact subsets
of an open set of regular parameters on which $q_+q_-\ne0$,
and restrict $z$ to $[z_0,z_1]\subset(0,\infty)$.
The standard Fourier estimates for the Eisenstein family give,
for every $A>0$ and some fixed $B$,
\begin{equation}\label{eq:uniform-fp-eisenstein-growth}
 \sup_{\tau\in K_1}|E(\tau,g^{-1})|\ll\lambda^B,
 \qquad
 \sup_{\tau\in K_1}|E_{\mathrm n}(\tau,g^{-1})|\ll_A\lambda^{-A}.
\end{equation}
Here one may take $B=1+2\Gamma$ for any
$\Gamma>\sup_{K_1}|\Re\tau|$.
The local Whittaker estimates and their summation giving the second
bound are proved in Lemma~\ref{lem:75-cusp-decay}; they are independent
of cusp regularization. The finite Fourier support confines the
rational frequencies to a fixed fractional ideal, and
Lemma~\ref{lem:uniform-lattice} makes the resulting rapid-decay
bound uniform. Removable factorization poles are treated by
Lemma~\ref{lem:normalization-removal}, applied to the holomorphic
nonconstant family.

The theta estimate is uniform for $z\in[z_0,z_1]$, by the compact
parameter assertions in the preceding lemmas. Thus
\begin{equation}\label{eq:uniform-fp-remainder-bound}
 \sup_{\tau\in K_1,\ z\in[z_0,z_1]}|\mathcal R(z,\tau,g)|
                   \ll\lambda^{B-M}+\lambda^{4-A}.
\end{equation}
Choosing $M>B+2$ and $A>6$ gives a common majorant
$O(\lambda^{-2})$, integrable against the cusp measure
$\lambda^{-2}d\lambda/\lambda$. Its tail tends to zero uniformly.
This proves local uniform convergence and holomorphy of the last
integral in \eqref{eq:finite-part-convergent-representation}.
For explicit differentiation, choose $\delta>0$ such that every
circle of radius $\delta$ about $K$ lies in $K_1$. Cauchy's estimate
bounds the $j$th spectral derivative of $\mathcal R$ by
$j!\delta^{-j}$ times the same integrable majorant. Hence
\begin{equation}\label{eq:uniform-fp-differentiation}
 \partial_\tau^j\int_{\lambda>T_0}\mathcal R\,dg^1
       =\int_{\lambda>T_0}\partial_\tau^j\mathcal R\,dg^1.
\end{equation}
The compact integral is holomorphic by local uniform holomorphy of
the smooth Eisenstein family, and $B_{T_0}$ is holomorphic away from
$q_+q_-=0$.

Finally, fix any spectral parameter $\tau_*$ and a small disk $D$
around it. Meromorphy of the smooth Eisenstein family supplies a
nonzero holomorphic function $d(\tau)$, independent of $g$, such that
\begin{equation}\label{eq:uniform-fp-cleared-family}
 E^{[d]}=dE,\qquad b_\epsilon^{[d]}=db_\epsilon,
 \qquad D_\epsilon^{[d]}=dD_\epsilon
\end{equation}
are holomorphic on $D$. For
$\GL_2$ the intertwining and Fourier formulas are given in
\cite[\S4.1.7]{MV}. On a fixed finite-dimensional compact-picture
space their finitely many coefficient functions have a common
local denominator, while the Fourier remainder is locally normally
convergent with its smooth moderate-growth seminorms. Choose the boundary of $D$ to avoid
poles. The uniform boundary estimates and the maximum modulus
principle give the same moderate-growth and rapid-decay bounds
for the cleared families throughout $D$.
Consequently
\begin{equation}\label{eq:uniform-fp-holomorphic-numerator}
 \mathcal H_d(z,\tau)=
 \int_{\mathcal F_{T_0}}E^{[d]}\Theta_{S,\Phi}\,dg^1
                        +\int_{\lambda>T_0}d(\tau)\mathcal R\,dg^1
\end{equation}
is holomorphic there. Away from the zeros of $d q_+q_-$,
\begin{equation}\label{eq:uniform-fp-meromorphic-expression}
 R^{\mathrm{reg}}_{\Phi,\tau}(z)
 =\frac1{d(\tau)}\left(
 \mathcal H_d(z,\tau)-z^{-2}
 \sum_{\epsilon=\pm1}D_\epsilon^{[d]}(\Phi,\tau)
                         \frac{T_0^{q_\epsilon}}{q_\epsilon}\right).
\end{equation}
Multiplication by $d q_+q_-$ makes the right side holomorphic on
$D$. This proves the claimed local meromorphy, including at the
resonant exponents $q_\epsilon=0$. When $(\eta_1,\eta_2)\ne(1,1)$, both
$D_\epsilon$ vanish identically, so the fixed-threshold formula
has no $q_\epsilon$ denominators and proves the stronger assertion.
\end{proof}

\begin{lemma}\label{lem:eis-zero-orbit}
For $T>1$ sufficiently large, put
\begin{equation}\label{eq:zero-orbit-truncated-average}
 \mathcal E_T(\tau)=\int_{\mathcal F_T}E(\tau,g^{-1})\,dg^1.
\end{equation}
Then, as an identity of meromorphic functions in $\tau$,
\begin{equation}\label{eq:zero-orbit-average-formula}
 \mathcal E_T(\tau)=\sum_{\epsilon=\pm1}
     e_\epsilon(\tau)
          \frac{T^{-1+2\epsilon\tau}}{-1+2\epsilon\tau}.
\end{equation}
In particular this average is identically zero unless $\eta_1=\eta_2=1$. At a regular parameter
$\tau\ne\pm1/2$, the zero orbit has cutoff finite part
\begin{equation}\label{eq:zero-orbit-finite-part}
 \operatorname{FP}_{T\to\infty}
               \bigl(\Phi(0)\mathcal E_T(\tau)\bigr)=0.
\end{equation}
\end{lemma}
\begin{proof}
For $\Re\tau>1/2$ truncate the inducing section by height:
\begin{equation}\label{eq:zero-orbit-truncated-function}
 \mathcal P_T(\tau,g)
  =\sum_{\gamma\in G_k/B_k}
     \mathbf1_{\{t(g\gamma)\le T\}}
                         f_\tau((g\gamma)^{-1}).
\end{equation}
The absolute value of the section is bounded by a constant times
$t(g\gamma)^{1+2\Re\tau}$. Thus the absolute unfolded integral
is bounded by a constant times
\begin{equation}\label{eq:zero-orbit-absolute-unfolding}
 C^2\int_0^T\lambda^{2\Re\tau-1}\frac{d\lambda}{\lambda}<\infty.
\end{equation}
Ordinary unfolding consequently gives
\begin{equation}\label{eq:zero-orbit-complex-unfolding}
 \int_{G^1/G_k}\mathcal P_T(\tau,g)\,dg^1
   =e_+(\tau)\frac{T^{2\tau-1}}{2\tau-1}.
\end{equation}
On $\mathcal F_T$ the sum equals $E(\tau,g^{-1})$.
On the complementary high cusp exactly the inducing summand for
that cusp is removed, so
\begin{equation}\label{eq:zero-orbit-on-high-cusp}
 \mathcal P_T(\tau,g)
 =E(\tau,g^{-1})-\lambda^{p_+}b_+(\tau;u,t_1,t_2).
\end{equation}
Its unipotent average is the reflected branch. Integration of this
branch over the complete transverse quotient and then in height
gives
\begin{equation}\label{eq:zero-orbit-high-cusp-integral}
 \int_{(G^1/G_k)\setminus\mathcal F_T}\mathcal P_T(\tau,g)\,dg^1
   =e_-(\tau)\frac{T^{-1-2\tau}}{1+2\tau}.
\end{equation}
Subtracting proves \eqref{eq:zero-orbit-average-formula} initially;
compact integration on $\mathcal F_T$ after local pole-clearing
continues the identity meromorphically. Character orthogonality
proves the asserted vanishing. For $\eta_1=\eta_2=1$,
compact averaging commutes with $M_\tau$ and projects onto the
one-dimensional spherical section, so
$e_-(\tau)=m(\tau)e_+(\tau)$; no such scalar replacement is
valid before this compact projection.

Away from the zero exponents $-1\pm2\tau=0$, the displayed
average has no constant term in its cutoff expansion, proving
\eqref{eq:zero-orbit-finite-part}. At an exceptional parameter the
complete expression must be specialized before its cutoff finite
part is taken. At every finite cutoff the zero orbit is retained
in the exact identity
\begin{equation}\label{eq:zero-orbit-exact-accounting}
 \int_{\mathcal F_T}E(\tau,g^{-1})\Theta_{S,\Phi}(z,g)\,dg^1
 =\int_{\mathcal F_T}E(\tau,g^{-1})P f_z(g)\,dg^1
                         +\Phi(0)\mathcal E_T(\tau).
\end{equation}
\end{proof}

\begin{proposition}\label{prop:eis-finite-unfolding}
For $T>1$, define
\begin{equation}\label{eq:eis-unfolded-cutoff}
 U_T(z,\tau)=\int_{\substack{G^1/B_k\\t(g)\le T}}
                            E(\tau,g^{-1})f_z(g)\,dg^1.
\end{equation}
At generic spectral parameters this integral is absolutely convergent,
and
\begin{equation}\label{ex:unfold-limit}
 R^{\rm reg}_{\Phi,\tau}(z)=
 \lim_{T\to\infty}\left\{U_T(z,\tau)-B_T(z,\tau)
                              +\Phi(0)\mathcal E_T(\tau)\right\}.
\end{equation}
The limit is locally uniform in $z>0$ and in generic $\tau$.
\end{proposition}
\begin{proof}
The function $f_z$ and height $t$ are right $B_k$-invariant. Put
\begin{equation}\label{eq:prop72-truncated-function}
 f_{z,T}(g)=\mathbf1_{\{t(g)\le T\}}f_z(g),\qquad
 P f_{z,T}(g)=\sum_{\gamma\in G_k/B_k}f_{z,T}(g\gamma).
\end{equation}
The latter series converges absolutely at every fixed $g$, since
its absolute sum is bounded by the Schwartz lattice sum on
$S_k\setminus\{0\}$.

We first justify integration on $G^1/B_k$, including its endpoint
$t(g)=\lambda\to0$. Let $f_z^{\mathrm{abs}}$ denote the function
formed with $|\Phi|$. Lemma~\ref{lem:eis-small-function} gives
\begin{equation}\label{eq:prop72-small-function}
 f_z^{\mathrm{abs}}(g)\ll_{M,z}\lambda^M
                         \qquad(0<\lambda\le1)
\end{equation}
for every $M$. The moderate-growth estimate on $G_\A$, before passage to
the automorphic quotient, gives $|E(\tau,g^{-1})|\ll\|g^{-1}\|^b$.
For $0<\lambda\le1$, matrix and inverse-matrix norms in these
coordinates are $O(\lambda^{-b_0})$ for a fixed $b_0$, because
the transverse variables are compact. Increasing $b$ therefore
gives $|E(\tau,g^{-1})|\ll\lambda^{-b}$ here.
Choosing $M>b+2$ yields
\begin{equation}\label{eq:prop72-small-integrability}
 \int_{\substack{G^1/B_k\\t(g)\le1}}
 |E(\tau,g^{-1})|f_z^{\mathrm{abs}}(g)\,dg^1
 \ll\int_0^1\lambda^{M-b-2}\frac{d\lambda}{\lambda}<\infty.
\end{equation}
On $1\le\lambda\le T$, boundedness on the compact parameter set
and locally uniform convergence of the Schwartz sum give the same
absolute integrability. Since
$E(\tau,(g\gamma)^{-1})=E(\tau,g^{-1})$, Tonelli and quotient
unfolding now imply
\begin{equation}\label{eq:prop72-unfolding}
 U_T(z,\tau)=\int_{G^1/G_k}E(\tau,g^{-1})P f_{z,T}(g)\,dg^1.
\end{equation}
The same equality with absolute values provides the justification
for interchanging the series and integral.

On $\mathcal F_T$, the truncated Poincar\'e sum equals $P f_z$.
On its complementary cusp there is a unique rational height greater
than $T$; choosing its representative with $t(g)=\lambda>T$,
only the identity summand is deleted. Consequently
\begin{equation}\label{eq:prop72-on-high-cusp}
 P f_{z,T}(g)=P f_z(g)-f_z(g)\qquad(\lambda>T).
\end{equation}
The projection of this high-cusp region in $G^1/B_k$ is injective
and preserves the inherited quotient measure. Using
$P f_z=\Theta_{S,\Phi}-\Phi(0)$ on $\mathcal F_T$ therefore gives
\begin{equation}\label{ex:exact-cutoff}
 U_T(z,\tau)=\mathcal A_T(z,\tau)-\Phi(0)\mathcal E_T(\tau)
                                      +D_T(z,\tau),
\end{equation}
where
\begin{equation}\label{eq:prop72-zero-and-error}
 D_T(z,\tau)=\int_{\substack{G^1/B_k\\t(g)>T}}
 E(\tau,g^{-1})\bigl(P f_z(g)-f_z(g)\bigr)\,dg^1.
\end{equation}
This is an exact finite-cutoff identity, including the zero vector.

By Lemma~\ref{lem:eis-other-cusps}, the absolute function difference
is $O_{M,z}(\lambda^{-M})$ for every $M$. With
$|E(\tau,g^{-1})|\ll\lambda^B$ in the high cusp, it follows that
\begin{equation}\label{eq:prop72-error-estimate}
 |D_T(z,\tau)|\ll\int_T^\infty\lambda^{B-M-2}\frac{d\lambda}{\lambda}
            =\frac{T^{B-M-2}}{M+2-B}\qquad(M+2>B).
\end{equation}
Thus $D_T=O_N(T^{-N})$ for every $N$. All constants are uniform
when $z$ and $\tau$ range over compact sets as in the statement.
Subtracting $B_T$ in \eqref{ex:exact-cutoff} and applying
Proposition~\ref{prop:eis-finite-part-exists} proves
\eqref{ex:unfold-limit} and its asserted uniformity.
\end{proof}

\section{The two constant-term contributions}\label{subsec:eis-constant}
We evaluate both nonzero singular orbits with the same height cutoff.
The compact-type dependence enters the test function. The two
inducing characters give different selection rules for the two orbits;
nontriviality of the central character alone does not imply vanishing.

For $\epsilon\in\{+1,-1\}$ write
\begin{equation}\label{eq:73-parameters}
 \begin{gathered}
 (\eta_{+,1},\eta_{+,2})=(\eta_1,\eta_2),\qquad
 (\eta_{-,1},\eta_{-,2})=(\eta_2,\eta_1),\\
 p_\epsilon=1+2\epsilon\tau,\qquad
 u_\epsilon=(2-p_\epsilon)/3,\qquad v_\epsilon=-p_\epsilon-1.
 \end{gathered}
\end{equation}
The functions $c_\epsilon(\tau,u)$ are the two compact restrictions
of the section and its intertwining transform. The spectral parameter
is generic, as above. For a character $\xi$ of
$\mathcal C_k=\A^1/k^\times$, write $\delta(\xi)=1$ when $\xi$ is trivial,
and $\delta(\xi)=0$ otherwise. Put
\begin{equation}\label{eq:73-full-selection}
 \Delta=\delta(\eta_1)\delta(\eta_2),\qquad \vol(\mathcal C_k)=C.
\end{equation}
Because all inducing characters are norm-normalized, triviality on
$\mathcal C_k$ is equivalent to triviality as a Hecke character.

\begin{lemma}\label{lem:73-functions}
Put
\begin{equation}\label{eq:73-function-splitting}
 f_{1,z}(g)=\sum_{\beta\in k^\times}
       \Phi(d(\underline z)g(0,0,0,\beta)),\qquad
 f_{2,z}(g)=\sum_{\substack{\alpha\in k^\times\\\beta\in k}}
       \Phi(d(\underline z)g(0,0,\alpha,\beta)).
\end{equation}
These are right $B_k$-invariant and sum to the nonzero singular
function. Define
\begin{equation}\label{eq:73-one-dimensional-tests}
 \begin{aligned}
 \Phi_\epsilon(X)&=\int_U c_\epsilon(\tau,u)\Phi(uX)\,du,\\
 F_\epsilon(t)&=\Phi_\epsilon(0,0,0,t),&
 f_\epsilon(t)&=\int_\A\Phi_\epsilon(0,0,t,v)\,dv,\\
 B_\epsilon&=\int_\A F_\epsilon(t)\,dt,&
 A_\epsilon&=\int_\A f_\epsilon(t)\,dt.
 \end{aligned}
\end{equation}
They are Schwartz--Bruhat functions or scalars, as appropriate,
and are locally meromorphic in $\tau$. They satisfy
\begin{equation}\label{eq:73-test-values}
 F_\epsilon(0)=\Phi(0)\int_Uc_\epsilon(\tau,u)\,du,
 \qquad f_\epsilon(0)=B_\epsilon,
 \qquad \widehat f_\epsilon(0)=A_\epsilon.
\end{equation}
Moreover the coefficient of the common cusp subtraction is
\begin{equation}\label{eq:73-D-selection}
 D_\epsilon(\Phi,E_\tau)
   =C^2\Delta A_\epsilon.
\end{equation}
\end{lemma}
\begin{proof}
The two sums partition the plane $x_1=x_2=0$ minus its origin,
and $B_k$ preserves each subset. Compact averaging against a smooth
function preserves the Schwartz space: every archimedean seminorm
of the average is bounded by a constant times a finite sum of
seminorms of $\Phi$, uniformly for the compact group action. At
finite places fixed support and invariance level are sent into
finitely many such data. Restriction and additive integration are
continuous on these spaces. The same estimates apply after a
local scalar has cleared the poles of $c_\epsilon$; hence they
also prove the meromorphic-family assertion. The identities at zero
are definitions. Finally, in the leading cusp coefficient the
norm-one variables occur only through
$\eta_{\epsilon,1}(t_1)^{-1}\eta_{\epsilon,2}(t_2)^{-1}$.
Their Haar integrals are $C\delta(\eta_{\epsilon,1})$ and
$C\delta(\eta_{\epsilon,2})$. This proves
\eqref{eq:73-D-selection}.
\end{proof}

\begin{lemma}\label{lem:73-tate}
For $h\in\cS(\A)$ and a character $\xi$ trivial on the positive
norm section, set
\begin{equation}\label{eq:73-theta-fiber}
 \vartheta_{h,\xi}(r)=\int_{\mathcal C_k}\xi(t)
      \sum_{\alpha\in k^\times}h(\underline r t\alpha)\,d^\times t.
\end{equation}
This function decreases faster than every power at infinity, and
\begin{equation}\label{ex:tate-poisson}
 \vartheta_{h,\xi}(r)
 =r^{-1}\vartheta_{\widehat h,\xi^{-1}}(r^{-1})
   +C\delta(\xi)(r^{-1}\widehat h(0)-h(0)).
\end{equation}
The integral
\begin{equation}\label{eq:73-tate-initial}
 Z(h,\xi,w)=\int_{\A^\times}h(t)\xi(t)|t|_\A^w\,d^\times t
          =\int_0^\infty\vartheta_{h,\xi}(r)r^w\frac{dr}{r}
\end{equation}
converges absolutely for $\Re w>1$ and continues meromorphically as
\begin{equation}\label{ex:tate-cont}
 \begin{aligned}
 Z(h,\xi,w)={}&\int_1^\infty\vartheta_{h,\xi}(r)r^w\frac{dr}{r}
 +\int_1^\infty\vartheta_{\widehat h,\xi^{-1}}(r)r^{1-w}\frac{dr}{r}\\
 &+C\delta(\xi)\left(\frac{\widehat h(0)}{w-1}
                                      -\frac{h(0)}w\right).
 \end{aligned}
\end{equation}
For $w\ne0,1$ we have 
\begin{equation}\label{eq:73-tate-exact-cutoff}
 \begin{aligned}
 \int_\delta^\infty\vartheta_{h,\xi}(r)r^w\frac{dr}{r}
  ={}&Z(h,\xi,w)
    -C\delta(\xi)\frac{\widehat h(0)\delta^{w-1}}{w-1}
    +C\delta(\xi)\frac{h(0)\delta^w}{w}\\
    &-\mathcal R_{h,\xi}(\delta,w),\\
 \mathcal R_{h,\xi}(\delta,w)
  :={}&\int_0^\delta
       \vartheta_{\widehat h,\xi^{-1}}(r^{-1})r^{w-1}\frac{dr}{r}.
 \end{aligned}
\end{equation}
The remainder is $O(\delta^N)$ for every $N$, uniformly on compact
sets of $w$ and bounded subsets of the Schwartz space. Thus
$Z(h,\xi,w)$ is entire in $w$ when $\xi\ne1$; when $\xi=1$
its only possible poles are at $0,1$, of residues $-Ch(0)$ and
$C\widehat h(0)$.
\end{lemma}
\begin{proof}
This is Tate's one-dimensional Poisson argument \cite{Tate},
with the character and measure normalizations retained explicitly.
Choose compact representatives for $\mathcal C_k$. The lattice estimate of
\S\ref{sec:radial}, uniformly over these representatives, bounds
the nonzero sum by $O_L(r^{-L})$ as $r\to\infty$, and the full
sum by $O(1+r^{-1})$ for $r\le1$. Since $|\xi|=1$, these
bounds justify integration and all quotient unfoldings on absolute
values. Scalar Poisson summation applied to $h(\underline r t\,\cdot)$
and the change $t\mapsto t^{-1}$ give
\eqref{ex:tate-poisson}; the zero terms have integral
$\int_{\mathcal C_k}\xi=C\delta(\xi)$. Norm disintegration proves
\eqref{eq:73-tate-initial}. Splitting at $r=1$ and applying
\eqref{ex:tate-poisson} on $(0,1)$ proves \eqref{ex:tate-cont}.
Both upper integrals are entire, since rapid decay absorbs every
fixed power and logarithm of $r$. Integrating the same identity
on $(0,\delta)$ and subtracting from the full integral gives
\eqref{eq:73-tate-exact-cutoff}. The remainder is entire in $w$;
if $a$ is a lower bound for $\Re w$ on a compact set, rapid decay
bounds its integrand by $O(r^{L+a-1})\,dr/r$. Choose
$L+a-1>N$ and integrate. All estimates involve continuous
Schwartz seminorms, proving the claimed uniformity. Continuation
then proves the identities at all stated arguments.
\end{proof}

\begin{proposition}
Let $U_{1,T}^\epsilon$ denote the contribution of the $\epsilon$
constant-term branch against $f_{1,z}$ on $G^1/B_k$ with
$\lambda\le T$. Then
\begin{equation}\label{ex:triple-cutoff}
 U_{1,T}^\epsilon
 =\frac C3\delta(\eta_{\epsilon,1}^2\eta_{\epsilon,2})
 z^{-u_\epsilon}\int_{zT^{-3}}^\infty
       \vartheta_{F_\epsilon,\eta_{\epsilon,1}}(r)
                       r^{u_\epsilon}\frac{dr}{r}.
\end{equation}
Consequently,
\begin{equation}\label{ex:triple-asymptotic}
 \begin{aligned}
 U_{1,T}^\epsilon={}&a_\epsilon z^{-u_\epsilon}+O_{N,z,\tau}(T^{-N})\\
 &+C^2\Delta\left(
 z^{-1}B_\epsilon\frac{T^{p_\epsilon+1}}{p_\epsilon+1}
 -F_\epsilon(0)\frac{T^{p_\epsilon-2}}{p_\epsilon-2}\right),\\
 a_\epsilon=a_\epsilon(\Phi,E_\tau)
 :={}&\frac C3\delta(\eta_{\epsilon,1}^2\eta_{\epsilon,2})
           Z(F_\epsilon,\eta_{\epsilon,1},u_\epsilon).
 \end{aligned}
\end{equation}
In particular,
\begin{equation}\label{eq:eis-r1}
 R^{\mathrm{reg}}_{1,\Phi,E_\tau}(z)
                      =\sum_\epsilon a_\epsilon z^{-u_\epsilon}.
\end{equation}
\end{proposition}
\begin{proof}
Compact integration gives the test function $F_\epsilon$ and the triple-root
argument $\underline z\underline\lambda^{-3}t_2^2t_1^{-1}\beta$.
The unipotent quotient has mass one. On $\mathcal C_k^2$ put
\begin{equation}\label{eq:73-triple-torus-change}
 s=t_2,\qquad t=t_2^2t_1^{-1},\qquad(t_1,t_2)=(s^2t^{-1},s).
\end{equation}
This group automorphism preserves Haar measure, while its character
weight becomes
\begin{equation}
 \eta_{\epsilon,1}(t_1)^{-1}\eta_{\epsilon,2}(t_2)^{-1}
 = (\eta_{\epsilon,1}^2\eta_{\epsilon,2})(s)^{-1}
                                      \eta_{\epsilon,1}(t).
\end{equation}
The $s$ integral is
$C\delta(\eta_{\epsilon,1}^2\eta_{\epsilon,2})$.
This proves the selection rule. The other quotient
integral and rational sum give
$\vartheta_{F_\epsilon,\eta_{\epsilon,1}}(z\lambda^{-3})$.
The height density is $\lambda^{p_\epsilon-2}d\lambda/\lambda$;
putting $r=z\lambda^{-3}$ gives \eqref{ex:triple-cutoff} because
$\lambda^{p_\epsilon-2}=z^{-u_\epsilon}r^{u_\epsilon}$ and
$d\lambda/\lambda=-dr/(3r)$. Apply Lemma~\ref{lem:73-tate}.
Its endpoint terms require both
$\eta_{\epsilon,1}^2\eta_{\epsilon,2}=1$ and
$\eta_{\epsilon,1}=1$, hence precisely $\Delta=1$. The identities
$-3(u_\epsilon-1)=p_\epsilon+1$ and
$-3u_\epsilon=p_\epsilon-2$ give precisely the displayed
coefficients and signs. The rapid Tate remainder gives arbitrary
order in $T$, uniformly on compact $z$-intervals and compact
sets of generic spectral parameters.
\end{proof}

\begin{proposition}
With the analogous notation for $f_{2,z}$,
\begin{equation}\label{ex:double-cutoff}
 U_{2,T}^{\epsilon,\mathrm{ct}}
 =C\delta(\eta_{\epsilon,1})z^{p_\epsilon}
       \int_{zT^{-1}}^\infty
       \vartheta_{f_\epsilon,\eta_{\epsilon,2}^{-1}}(r)
                                    r^{v_\epsilon}\frac{dr}{r}.
\end{equation}
It follows that
\begin{equation}\label{ex:double-asymptotic}
 \begin{aligned}
 U_{2,T}^{\epsilon,\mathrm{ct}}={}&d_\epsilon z^{p_\epsilon}+O_{N,z,\tau}(T^{-N})\\
 &+C^2\Delta\left(z^{-2}A_\epsilon
                \frac{T^{p_\epsilon+2}}{p_\epsilon+2}
            -z^{-1}B_\epsilon
                \frac{T^{p_\epsilon+1}}{p_\epsilon+1}\right),\\
 d_\epsilon=d_\epsilon(\Phi,E_\tau)
 :={}&C\delta(\eta_{\epsilon,1})
             Z(f_\epsilon,\eta_{\epsilon,2}^{-1},v_\epsilon),
 \end{aligned}
\end{equation}
and hence
\begin{equation}\label{eq:eis-r2ct}
 R^{\mathrm{reg},\mathrm{ct}}_{2,\Phi,E_\tau}(z)
                                  =\sum_\epsilon d_\epsilon z^{p_\epsilon}.
\end{equation}
\end{proposition}
\begin{proof}
For fixed $\alpha\in k^\times$ the two nonzero coordinates are
\begin{equation}\label{eq:73-double-coordinates}
 x=\underline z\underline\lambda^{-1}t_2\alpha,
 \qquad q(c\alpha+\beta),\qquad
 q=\underline z\underline\lambda^{-3}t_2^2t_1^{-1}.
\end{equation}
Multiplication by $\alpha$ preserves the quotient measure on
$\A/k$, by the product formula. Additive unfolding therefore gives
\begin{equation}\label{eq:73-double-additive-unfolding}
 \int_{\A/k}\sum_{\beta\in k}
 \Phi_\epsilon(0,0,x,q(c\alpha+\beta))\,dc
 =|q|_\A^{-1}f_\epsilon(x)=z^{-1}\lambda^3f_\epsilon(x).
\end{equation}
Tonelli applied to absolute values justifies the unfolding. The
remaining test is independent of $t_1$, so its weighted integral
is $C\delta(\eta_{\epsilon,1})$. Integration in $t_2$ and
summation in $\alpha$ give the weighted theta function with
character $\eta_{\epsilon,2}^{-1}$. Thus the remaining height
integral is
\[
 C\delta(\eta_{\epsilon,1})z^{-1}
 \int_0^T\lambda^{p_\epsilon+1}
 \vartheta_{f_\epsilon,\eta_{\epsilon,2}^{-1}}
                   (z\lambda^{-1})\frac{d\lambda}{\lambda}.
\]
Set $r=z\lambda^{-1}$. Since
$z^{-1}\lambda^{p_\epsilon+1}=z^{p_\epsilon}r^{v_\epsilon}$,
this proves \eqref{ex:double-cutoff}. Apply
Lemma~\ref{lem:73-tate}, using
$\widehat f_\epsilon(0)=A_\epsilon$ and
$f_\epsilon(0)=B_\epsilon$. The endpoint terms require
$\eta_{\epsilon,1}=\eta_{\epsilon,2}=1$, hence $\Delta=1$. The identities
$v_\epsilon-1=-p_\epsilon-2$ and
$v_\epsilon=-p_\epsilon-1$ give \eqref{ex:double-asymptotic},
including its uniform error bound.
\end{proof}

\begin{proposition}
\label{prop:73-recombination}
The nonconstant integral
\begin{equation}\label{eq:73-nonconstant-integral}
 R^{\mathrm{nc}}_{\Phi,E_\tau}(z)
    =\int_{G^1/B_k}E_{\tau,\mathrm n}(g^{-1})f_{2,z}(g)\,dg^1
\end{equation}
converges absolutely. The full singular pairing is
\begin{equation}\label{eq:eis-r-decomposition}
 R^{\mathrm{reg}}_{\Phi,E_\tau}(z)
 =\sum_{\epsilon=\pm1}
   \left(a_\epsilon(\Phi,E_\tau)z^{-u_\epsilon}
                 +d_\epsilon(\Phi,E_\tau)z^{p_\epsilon}\right)
             +R^{\mathrm{nc}}_{\Phi,E_\tau}(z).
\end{equation}
All assertions are locally uniform in $z>0$ and in generic spectral
parameters.
\end{proposition}
\begin{proof}
At small height the small-height function bound absorbs the moderate
growth of the Eisenstein series and its constant term. At large
height the nonconstant part is rapidly decreasing, while the absolute
function is $O_z(\lambda^4)$ by the two-dimensional lattice lemma.
These statements give an integrable majorant after multiplication
by $\lambda^{-2}d\lambda/\lambda$, and prove the first assertion.
The nonconstant term against $f_{1,z}$ is zero because that function
is independent of the unipotent variable, whose nonconstant Fourier
average is zero. The zero-orbit contribution has the cutoff expansion
established in \S\ref{subsec:eis-finite-part}; its terms
$C^2\Delta F_\epsilon(0)T^{p_\epsilon-2}/(p_\epsilon-2)$
cancel exactly the corresponding terms of
\eqref{ex:triple-asymptotic}. The two terms involving $B_\epsilon$
in the triple- and double-root expansions also cancel. Their only
remaining primitive is
\[
 z^{-2}\sum_\epsilon D_\epsilon(\Phi,E_\tau)
                          \frac{T^{p_\epsilon+2}}{p_\epsilon+2},
\]
which is exactly the common cusp subtraction. The finite-cutoff
unfolding theorem from \S\ref{subsec:eis-finite-part} now proves
\eqref{eq:eis-r-decomposition}. Uniform convergence follows from the
uniform function bounds, rapid Fourier decay, and weighted Tate remainder
estimate already proved; no independent choices of orbit finite
parts are made.
\end{proof}
\section{The nonconstant terms and Whittaker functions}
\label{subsec:eis-nonconstant}

The nonconstant part of the Eisenstein series gives an ordinary,
absolutely convergent integral at each fixed radius.  This section
identifies its full radial Mellin transform with a finite sum of
two-variable Whittaker integrals.  The inducing characters are
$\eta_1,\eta_2$; they may be ramified, and the section may
have any fixed finite set of compact types. Its central character
is $\omega=\eta_1\eta_2$.  Dependence on the fixed inducing character and section
is suppressed in $E(\tau,h)$ and the associated notation.

Write $E_{\mathrm c}$ for the constant term of the left automorphic Eisenstein
series and set
\begin{equation}\label{eq:74-constant-projection}
 E_{\mathrm c}(\tau,h)=\int_{\A/k}E(\tau,n(b)h)\,db,
 \qquad E_{\rm n}(\tau,h)=E(\tau,h)-E_{\mathrm c}(\tau,h).
\end{equation}
Uniform assertions are understood on compact subsets of the
regular spectral locus.  Near a pole, they apply
after multiplication by a local holomorphic scalar clearing the
poles of these families.  No exclusion arising from the separate
Tate integrals of the preceding section is needed here.

Put
\begin{equation}\label{eq:74-two-parts}
 \begin{aligned}
 f_{1,z}(g)&=\sum_{\beta\in k^\times}
                  \Phi(\underline z g(0,0,0,\beta)),\\
 f_{2,z}(g)&=\sum_{\substack{\alpha\in k^\times\\\beta\in k}}
                  \Phi(\underline z g(0,0,\alpha,\beta)).
 \end{aligned}
\end{equation}
Thus $f_z=f_{1,z}+f_{2,z}$.  Both functions are right
$B_k$-invariant.

\begin{proposition}
\label{prop:74-fixed-radius}
For every $z>0$ and every regular $\tau$, the integral
\begin{equation}\label{ex:Rnc}
 R^{\rm nc}_{\Phi,\tau}(z)
 =\int_{G^1/B_k}E_{\rm n}(\tau,g^{-1})f_{2,z}(g)\,dg^1
\end{equation}
converges absolutely, including after replacing the summands of
$f_{2,z}$ by their absolute values.  Convergence is locally uniform
in $z>0$ and in regular $\tau$.  The integral is holomorphic in
$\tau$ on the regular locus and defines a locally meromorphic
family in $\tau$.

The triple-root contribution vanishes:
\begin{equation}\label{eq:74-triple-vanishes}
 \int_{G^1/B_k}E_{\rm n}(\tau,g^{-1})f_{1,z}(g)\,dg^1=0.
\end{equation}
If $x_2=(0,0,1,0)$ and
$H_2(k)=\{a(q,1):q\in k^\times\}$, then
\begin{equation}\label{ex:fixed-radius-unfold}
 R^{\rm nc}_{\Phi,\tau}(z)
 =\int_{G^1/H_2(k)}E_{\rm n}(\tau,g^{-1})
                         \Phi(\underline z gx_2)\,dg^1,
\end{equation}
again with absolute convergence.
\end{proposition}

\begin{proof}
The constant term $E_{\mathrm c}$ is left $B_k$-invariant.  Indeed, left
translation by $n(q)$, $q\in k$, translates its integration variable,
whereas left translation by a rational diagonal matrix rescales that
variable by an element of $k^\times$.  Such rescaling preserves the
Haar measure of $\A/k$ by the product formula.  Left $G_k$-invariance
of $E$ then proves the assertion.  Consequently $E_{\rm n}$ is left
$B_k$-invariant, so the integrand in \eqref{ex:Rnc} is well defined.

Use the right-Iwasawa coordinates
\begin{equation}\label{eq:74-fixed-coordinates}
 g=u a(\underline\lambda)a(t_1,t_2)n(c),\qquad
 dg^1=du\,\lambda^{-2}\frac{d\lambda}{\lambda}
                         d^\times t_1d^\times t_2dc,
\end{equation}
with representatives of $t_i\in\A^1/k^\times$ and $c\in\A/k$
in fixed compact sets.  Denote the full absolute function by
\begin{equation}
 S_z(g)=\sum_{(\alpha,\beta)\in k^2\setminus\{(0,0)\}}
             |\Phi(\underline z g(0,0,\alpha,\beta))|.
\end{equation}
For $z$ in a compact subset $J$ of $(0,\infty)$, the small-height
estimate \eqref{ex:function-small} gives
\begin{equation}\label{eq:74-function-endpoints}
 S_z(g)\ll_{J,N}\lambda^N\quad(0<\lambda\le1),
 \qquad
 S_z(g)\ll_J\lambda^4\quad(\lambda\ge1).
\end{equation}
Here and below the constants may depend on $\Phi$.
Both statements follow from Lemma~\ref{lem:triangular-lattice}
on the repeated-root plane: its archimedean scales are
$(z\lambda^{-1})^{1/n}$ and $(z\lambda^{-3})^{1/n}$,
where $n=[k:\mathbb Q]$.  The nonzero stratum bound gives the first
estimate, while the full counting bound gives
$(1+z^{-1}\lambda)(1+z^{-1}\lambda^3)\ll_J\lambda^4$ for the second.

On a compact spectral set $K$, moderate growth at small height and
rapid decay of the nonconstant term at large height give constants
$B_K$ and, for each $L$, $C_{K,L}$ such that
\begin{equation}\label{eq:74-eisenstein-endpoints}
 |E_{\rm n}(\tau,g^{-1})|
 \le
 \begin{cases}
 C_K\lambda^{-B_K},&0<\lambda\le1,\\
 C_{K,L}\lambda^{-L},&\lambda\ge1.
 \end{cases}
\end{equation}
These are the uniform Eisenstein estimates stated in
Section~\ref{subsec:eis-family}; their Whittaker proof is given in
Section~\ref{subsec:eis-local}.
The transverse variables have finite total measure.  Combining
\eqref{eq:74-function-endpoints} and \eqref{eq:74-eisenstein-endpoints}
with \eqref{eq:74-fixed-coordinates} bounds the absolute integral by
a constant times
\begin{equation}
 \int_0^1\lambda^{N-B_K-2}\frac{d\lambda}{\lambda}
 +\int_1^\infty\lambda^{2-L}\frac{d\lambda}{\lambda}.
\end{equation}
Choose $N>B_K+2$ and $L>2$.  This proves absolute and locally
uniform convergence, both for $f_{2,z}$ and for $f_{1,z}$.
Dominated integration preserves holomorphy in $\tau$.
Applying the same argument to a locally pole-cleared Eisenstein
family proves the meromorphic-family assertion.

For $g$ as in \eqref{eq:74-fixed-coordinates}, the function
$f_{1,z}(g)$ is independent of $c$, since
$n(c)(0,0,0,\beta)=(0,0,0,\beta)$.  Moreover
\begin{equation}
 \int_{\A/k}E_{\rm n}
   (\tau,n(-c)a(t_1,t_2)^{-1}a(\underline\lambda)^{-1}u^{-1})\,dc=0
\end{equation}
by \eqref{eq:74-constant-projection}.  Absolute convergence permits
this integration to be performed first, proving
\eqref{eq:74-triple-vanishes}.

Finally, the action on $x_2$ is
\begin{equation}\label{eq:74-orbit-action}
 n(c)a(A,D)x_2=(0,0,D,cD).
\end{equation}
It follows directly that
\begin{equation}
 B_kx_2=\{(0,0,\alpha,\beta):\alpha\in k^\times,\ \beta\in k\},
 \qquad \operatorname{Stab}_{B_k}(x_2)=H_2(k).
\end{equation}
Thus
\begin{equation}
 f_{2,z}(g)=\sum_{b\in B_k/H_2(k)}\Phi(\underline z gbx_2).
\end{equation}
Inserting this sum in \eqref{ex:Rnc}, absolute convergence permits
ordinary quotient unfolding.  The identity
$E_{\rm n}(\tau,(gb)^{-1})=E_{\rm n}(\tau,g^{-1})$
and the inherited quotient measures give
\eqref{ex:fixed-radius-unfold}.  The orbit parametrization is
bijective, so no additional multiplicity occurs.
\end{proof}

The next lemma fixes the signs in the ordinary automorphic Fourier expansion for the lower unipotent subgroup. 

\begin{lemma}
\label{lem:74-fourier}
Define
\begin{equation}\label{eq:74-whittaker-definition}
 W_\tau(h)=\int_{\A/k}E(\tau,n(b)h)\psi_k(b)\,db.
\end{equation}
Then
\begin{equation}\label{eq:74-whittaker-equivariance}
 W_\tau(n(b)h)=\psi_k(-b)W_\tau(h),
\end{equation}
and
\begin{equation}\label{eq:74-fourier-expansion}
 E_{\rm n}(\tau,h)
   =\sum_{\alpha\in k^\times}
                  W_\tau(a(\alpha^{-1},1)h).
\end{equation}
The latter series converges absolutely and locally uniformly in $h$
and in regular $\tau$.  It has central character $\omega$ and is right $U$-finite.
These assertions hold for every right compact translate of the
inducing section.
\end{lemma}

\begin{proof}
For fixed $h$, expand the smooth function $b\mapsto E(\tau,n(b)h)$
on the compact abelian group $\A/k$.  Its characters are indexed by
$k$.  The zero coefficient is $E_{\mathrm c}(\tau,h)$.
For $\alpha\in k^\times$, rational automorphy and rescaling the additive
variable by $\alpha$ identify the coefficient of $\psi_k(-\alpha b)$
with $W_\tau(a(\alpha^{-1},1)h)$.  Setting $b=0$ gives the
asserted expansion.  Translating the integration variable proves
Whittaker equivariance, with the displayed minus sign.

On a compact set of $h$ and regular $\tau$, finite adelic invariance
restricts the nonzero Fourier coefficients to a fixed fractional
ideal in $k$.  Integration by parts in the archimedean variables
bounds them by $C_N(1+\|\alpha\|)^{-N}$ for every $N$, uniformly
on this compact set.  Lemma~\ref{lem:uniform-lattice} gives absolute
and locally uniform summability.  The same argument applies after
clearing spectral poles.  Right compact finiteness and the central transformation law
are inherited from $E$ under the defining integral.
\end{proof}

To retain the full compact type, use the compact realization of the
induced representation.  Its right $U$-action is independent of $\tau$.
Choose a basis of the finite-dimensional $U$-span of the fixed flat
section.  Its basis sections give Eisenstein series $E_j(\tau,h)$ and
matrix coefficients $c_j(u)$ such that
\begin{equation}\label{eq:74-compact-types}
 E(\tau,hu^{-1})=\sum_{j=1}^d c_j(u)E_j(\tau,h).
\end{equation}
Taking constant terms or Whittaker coefficients commutes with this
finite identity.  Denote the corresponding coefficients by
$E_{{\rm n},j}$ and $W_{\tau,j}$, respectively, with the sign convention
\eqref{eq:74-whittaker-definition}.  The coefficients $c_j$ are smooth
and bounded, and are locally constant in the finite compact variables.
They are independent of the spectral parameter. For a holomorphic
linear combination of flat sections, use a common finite set of
compact types. For each occurring irreducible unitary compact
representation $\sigma$ of dimension $d_\sigma$, define
\[
 P_{\sigma,ij}=d_\sigma\int_U\overline{\sigma_{ij}(v)}R(v)\,dv.
\]
Schur orthogonality identifies these operators with the matrix units
on the $\sigma$ factor, tensored with the identity on its multiplicity
space, and makes them zero on the other isotypic components. Hence
\[
 R(u^{-1})=\sum_{\sigma,i,j}\sigma_{ij}(u^{-1})P_{\sigma,ij}
\]
on the chosen finite set of compact types. Taking $E_j$ to be these
fixed compact convolutions of $E$, with the indices reenumerated, gives
\eqref{eq:74-compact-types} also for the varying section. In particular,
every $E_j$ is holomorphic wherever the full family $E$ is holomorphic,
including when holomorphic section coefficients cancel a spectral
pole. A common pole-clearing scalar for $E$ also works for these
compact projections.  Define
\begin{equation}\label{eq:74-two-averages}
 \Phi_j(X)=\int_U c_j(u)\Phi(uX)\,du,
 \qquad \Phi_j^{\rm abs}(X)=\int_U|c_j(u)|\,|\Phi(uX)|\,du.
\end{equation}
Each $\Phi_j$ is Schwartz--Bruhat, whereas $\Phi_j^{\rm abs}$ need
not be smooth.  The latter dominates $|\Phi_j|$ and is the appropriate
majorant before compact averaging.  Put
\begin{equation}\label{eq:74-positive-marginal}
 P_{\Phi,j}(x)=\int_\A\Phi_j^{\rm abs}(0,0,x,v)\,dv.
\end{equation}
For real $\sigma$ define the nonnegative, initially possibly infinite,
quantities
\begin{equation}\label{eq:74-majorant-factors}
 \begin{aligned}
 \mathcal A_{\Phi,j}(\sigma)
   &=\int_{\A^\times}P_{\Phi,j}(x)|x|_\A^{2\sigma-1}d^\times x,\\
 \mathcal B_{\tau,j}(\sigma)
   &=\int_{\A^\times}|W_{\tau,j}(a(1,t))|
                                  |t|_\A^\sigma d^\times t.
 \end{aligned}
\end{equation}

\begin{lemma}
\label{lem:normalization-removal}
Let $A(\tau,x)$ be holomorphic in $\tau$ on a neighborhood of a
closed disc $D$, and measurable in $x$. Suppose a factorized
expression for $A$ is regular on $\partial D$, and its absolute
integral against a positive measure $d\mu(x)$ is uniformly bounded
there. Then that integral is locally uniformly bounded on $D$.
A fixed nonnegative weight may be included if the corresponding
weighted boundary bound holds. The conclusion applies to an
$x$-derivative if it exists, is holomorphic in $\tau$, and satisfies
the corresponding boundary integral bound. A meromorphic family is
handled after multiplication by a scalar clearing its genuine poles.
\end{lemma}
\begin{proof}
On each smaller concentric disc, Cauchy's formula gives
$|A(\tau,x)|\le C_D\int_{\partial D}|A(\zeta,x)|\,|d\zeta|$.
Integration and Tonelli's theorem transfer the boundary bound to
the smaller disc. Apply the same argument to the weighted function,
to each derivative, or to a pole-cleared family. A finite cover
gives uniformity on compact parameter sets. When a growth bound
involves $|\Re\tau|$, the disc is chosen inside the available strict
margin; an arbitrarily small loss in the exponent absorbs its radius.
\end{proof}

\begin{lemma}
\label{lem:74-majorant}
Let $\tau$ be regular and put $\gamma=|\Re\tau|$.
If $\sigma>\max\{1,\tfrac12+\gamma\}$, then both quantities in
\eqref{eq:74-majorant-factors} are finite.  Moreover
\begin{equation}\label{eq:74-majorant-factorization}
 \begin{aligned}
 &\int_{(\A^\times)^2}P_{\Phi,j}(x)|W_{\tau,j}(a(1,xy))|
       |x|_\A^{3\sigma-1}|y|_\A^\sigma d^\times x d^\times y\\
 &\hspace{25mm}=\mathcal A_{\Phi,j}(\sigma)
                         \mathcal B_{\tau,j}(\sigma).
 \end{aligned}
\end{equation}
The bounds are uniform on compact subsets of this region with regular
$\tau$.  They remain valid after multiplication of the integrand by
any fixed power of $|\log(|x|_\A^3|y|_\A)|$.  They also hold for
locally pole-cleared spectral families.
\end{lemma}

\begin{proof}
Boundedness of $c_j$, compactness of $U$, and the Schwartz bounds give,
for a finite set $S$ containing all infinite places and every sufficiently
large $L$,
\begin{equation}\label{eq:74-marginal-bound}
 P_{\Phi,j}(x)\le C_{j,L}
 \prod_{v\mid\infty}(1+\|x_v\|^2)^{-L}
 \prod_{\substack{v\in S\\v\nmid\infty}}
             \mathbf1_{\varpi_v^{-b_v}\OO_v}(x_v)
 \prod_{v\notin S}\mathbf1_{\OO_v}(x_v),
\end{equation}
where $b_v\ge0$.  Indeed, an archimedean Schwartz seminorm bounds the
test uniformly under the compact action by a product of decaying
factors in the two coordinates of the repeated-root plane.  Integrating
the last coordinate gives the asserted decay in $x$.  At the finite
places, enlarge the support to a product of fractional ideals.  Outside
$S$, the data and compact action are integral, and the additive unit
ball has volume one.  This argument does not require factorization
of the archimedean test among the places.
With $w=2\sigma-1$, the unramified local majorant has integral
\begin{equation}
 \int_{k_v^\times}\mathbf1_{\OO_v}(x)|x|_v^w d^\times x
                  =(1-q_v^{-w})^{-1}.
\end{equation}
The product converges for $w>1$; the finitely many other local factors
converge for $w>0$, using Schwartz decay at archimedean infinity.
This proves finiteness of $\mathcal A_{\Phi,j}$.

For the Whittaker factor, a smooth finite-level, archimedean
$U$-finite section is a finite sum of factorizable sections.  The
Jacquet-integral calculation of the nonzero Fourier coefficient
therefore gives, away from a discrete set of normalization parameters,
a finite sum of products
\begin{equation}\label{eq:74-whittaker-product}
 W_{\tau,j}(a(1,t))
     =\sum_\ell \kappa_{j,\ell}(\tau)
                         \prod_v w_{\tau,j,\ell,v}(t_v).
\end{equation}
This identity is first obtained in an absolute-convergence region and
then continued meromorphically; see \cite[Section~4.1.7]{MV} and the
local construction in \S\ref{subsec:eis-local}.  All but finitely many
local factors are the normalized spherical factors.  At such a place,
write
\begin{equation}\label{eq:74-unramified-bound}
 \begin{gathered}
 \alpha_v=\eta_{2,v}(\varpi_v)q_v^{-\tau},\qquad
 \beta_v=\eta_{1,v}(\varpi_v)q_v^{\tau},\\
 w_{\tau,v}(\varpi_v^m u)
   =q_v^{-m/2}\sum_{a=0}^{m}\alpha_v^a\beta_v^{m-a}
       \quad(m\ge0, u\in\OO_v^\times),\\
 w_{\tau,v}(\varpi_v^m u)=0\quad(m<0).
 \end{gathered}
\end{equation}
This is the $\GL_2$ spherical Whittaker formula
\cite{ShintaniWhittaker,CS}.  Since the inducing characters are unitary,
its absolute Mellin integral is at most
\begin{equation}
 \sum_{m\ge0}(m+1)q_v^{-m(\sigma+1/2-\gamma)}
       =(1-q_v^{-(\sigma+1/2-\gamma)})^{-2}.
\end{equation}
The product of these bounds converges when
$\sigma>1/2+\gamma$.

At each of the finitely many remaining places,
Lemma~\ref{lem:75-whittaker}, for the specified ramified or nonspherical
section, gives for every $\eta>0$
\begin{equation}\label{eq:74-local-zero-bound}
 |w_{\tau,j,\ell,v}(t)|
   \ll |t|_v^{1/2-\gamma-\eta}\qquad(0<|t|_v\le1),
\end{equation}
together with bounded upper support at finite places and arbitrary
decay at archimedean infinity.  Their Mellin integrals converge if
$\sigma+1/2-\gamma-\eta>0$.  Choose $\eta$ inside the strict margin.
This proves finiteness of $\mathcal B_{\tau,j}$ wherever the chosen
factorization is regular.

The full coefficient $W_{\tau,j}$ is holomorphic wherever the
Eisenstein family is holomorphic, since its defining unipotent
integral is over a compact quotient. Lemma~\ref{lem:normalization-removal}
therefore extends these estimates through every removable pole of
\eqref{eq:74-whittaker-product}. Clearing genuine Eisenstein poles
gives the corresponding meromorphic-family bounds.

Finally put $t=xy$.  Haar invariance and
\begin{equation}
 |x|_\A^{3\sigma-1}|y|_\A^\sigma
                   =|x|_\A^{2\sigma-1}|t|_\A^\sigma
\end{equation}
give \eqref{eq:74-majorant-factorization} by Tonelli.  The preceding
bounds make both sides finite.  Uniformity in $\sigma$ follows by
bounding with exponents on the endpoints of a slightly larger compact
interval.  For logarithmic factors use
$|\log r|^m\ll_{m,\delta}r^\delta+r^{-\delta}$, taking $\delta>0$
smaller than each strict convergence margin.
\end{proof}

\begin{proposition}
\label{prop:74-completion}
For a regular parameter $\tau$ put $\gamma=|\Re\tau|$ and define
\begin{equation}\label{eq:74-partial-fourier}
 F_j(x,y)=\int_\A\Phi_j(0,0,x,v)\psi_k(yv)\,dv.
\end{equation}
Then $F_j\in\cS(\A^2)$.  In the common initial region
\begin{equation}\label{ex:M-initial}
 \Re s>2+\gamma,
\end{equation}
one has an identity of absolutely convergent integrals
\begin{equation}\label{eq:eis-nc-completion}
 \begin{aligned}
 M^{\rm nc}_{\Phi,\tau}(s)
  &:=\int_0^\infty z^{2s}R^{\rm nc}_{\Phi,\tau}(z)\frac{dz}{z}\\
  &=\frac12\sum_j\int_{(\A^\times)^2}F_j(x,y)
       W_{\tau,j}(a(1,xy))\\
  &\qquad\times\omega(x)^{-2}\omega(y)^{-1}|x|_\A^{3s-1}|y|_\A^s
       d^\times x d^\times y.
 \end{aligned}
\end{equation}
Before Fourier integration in the last additive coordinate, an
absolute majorant is
$\tfrac12\sum_j\mathcal A_{\Phi,j}(\Re s)
                        \mathcal B_{\tau,j}(\Re s)$.
The identity is jointly holomorphic in $(s,\tau)$ on this region
with $\tau$ regular, and holds as an identity of locally meromorphic
families across spectral poles.
\end{proposition}

\begin{proof}
Compact averaging, restriction to a linear subspace, and partial
Fourier transformation preserve the Schwartz--Bruhat space, including
when the infinite-place variables are coupled. Thus $F_j\in\cS(\A^2)$.
The coordinate calculation is the one underlying \cite[Proposition~3.2]{LTcuspidal}. We supply it here for the nonconstant Eisenstein term. It uses only the orbit action and the Fourier expansion in Lemma~\ref{lem:74-fourier}, together with the stated Whittaker equivariance and central character.
Its absolute convergence follows from the present majorant lemma,
as we now check without assuming any automorphy of $E_{\mathrm n}$
under $G_k$.

Start with the nonnegative integral
\[
 \mathcal I(\sigma,\tau)=\int_{G_\A/H_2(k)}|\det g|^\sigma
             |E_{\mathrm n}(\tau,g^{-1})|\,|\Phi(gx_2)|\,dg.
\]
The integrand is well defined because $E_{\mathrm n}$ is invariant
under the rational diagonal subgroup. Write $g=u n(c)a(A,D)$;
$A$ runs over $\A^\times/k^\times$ and $D$ over $\A^\times$.
Apply the finite compact-type expansion, put $v=cD$, and unfold
the nonnegative Fourier sum by Tonelli, as follows. The additive change has Jacobian $dc=|D|_\A^{-1}dv$, so the absolute determinant weight becomes $|A|_\A^\sigma|D|_\A^{\sigma-1}$. For a frequency $\alpha\in k^\times$, lower Whittaker equivariance gives
\begin{equation}\label{eq:unfold-frequency-term}
W_{\tau,j}\bigl(a((A\alpha)^{-1},D^{-1})n(-v/D)\bigr)
=\psi_k(A\alpha v/D^2)\omega(A\alpha)^{-1}
W_{\tau,j}(a(1,A\alpha/D)).
\end{equation}
Indeed, commuting the unipotent matrix to the left gives $n(-A\alpha v/D^2)$; the remaining scalar matrix is $d((A\alpha)^{-1})$. Tonelli's theorem unfolds the nonnegative frequency sum from $A\in\A^\times/k^\times$ to $A\in\A^\times$. Integration in $v$ is bounded by $P_{\Phi,j}(D)$.
The resulting bound is
\begin{equation}\label{eq:74-full-absolute-bound}
 \begin{aligned}
 \mathcal I(\sigma,\tau)
 &\le\frac12\sum_j\int_{(\A^\times)^2}
 P_{\Phi,j}(x)|W_{\tau,j}(a(1,t))|
 |x|^{2\sigma-1}|t|^\sigma\,d^\times x\,d^\times t\\
 &=\frac12\sum_j\mathcal A_{\Phi,j}(\sigma)
                    \mathcal B_{\tau,j}(\sigma)<\infty.
 \end{aligned}
\end{equation}
Here $t=A/D$, and the central factors have absolute value one.

Radial disintegration and Proposition~\ref{prop:74-fixed-radius}
now give
\begin{equation}\label{eq:74-completed-group-integral}
 M^{\mathrm{nc}}_{\Phi,\tau}(s)=
 \int_{G_\A/H_2(k)}|\det g|^s
                 E_{\mathrm n}(\tau,g^{-1})\Phi(gx_2)\,dg.
\end{equation}
The identity $E_{\mathrm n}(\tau,d(\underline z)^{-1}h^{-1})
=E_{\mathrm n}(\tau,h^{-1})$ uses
$\omega(\underline z)=1$. In the absolutely convergent calculation
one may now remove the absolute values. The change
$x=D$, $y=A/D^2$ preserves multiplicative Haar measure and gives
the weights $|x|^{3s-1}|y|^s\omega(x)^{-2}\omega(y)^{-1}$,
the argument $W_{\tau,j}(a(1,xy))$, and the additive character
$\psi_k(yv)$, by the preceding frequency identity.
Integration in $v$ proves \eqref{eq:eis-nc-completion}.

The uniform majorant, with slightly enlarged strict margins,
permits differentiation in $s$ and dominated integration in $\tau$.
Cauchy's formula supplies the corresponding spectral derivative
bounds. This proves joint holomorphy at regular parameters;
locally clearing genuine spectral poles proves the meromorphic
family assertion.
\end{proof}

The majorant actually proves the larger initial region
$\Re s>\max\{1,\tfrac12+|\Re\tau|\}$.
The common sufficient region \eqref{ex:M-initial} is retained in the
continuation theorem.  Equation~\eqref{eq:eis-nc-completion} is an
identity of ordinary integrals before any meromorphic continuation
in $s$.  The next section continues these integrals by local analysis and Euler factorization.

\section{Local continuation}
\label{subsec:eis-local}

We retain the two norm-normalized unitary inducing characters
$\eta_1,\eta_2$ and put $\omega=\eta_1\eta_2$. The compact-type
expansion gives finitely many pairs $(F_j,W_{j,\tau})$. The local
continuation is uniform in the spectral parameter, includes
reducibility points, and allows tests coupling different places.

The two characters of the local standard $L$-factor are
\begin{equation}\label{eq:75-inducing-characters}
 \mu_{1,v,\tau}=\eta_{2,v}|\cdot|_v^\tau,
 \qquad \mu_{2,v,\tau}=\eta_{1,v}|\cdot|_v^{-\tau}.
\end{equation}
Conjugation by $w_0$ converts our lower induced model into this
upper induced model and carries $a(1,t)$ to $a(t,1)$. Their
product is $\omega_v$. Write $\pi_\tau$ for this induced
representation, with its standard $L$-factor understood as the
product of the two character factors at a reducibility point too.

For a local Whittaker family $w_{\tau,v}$ put
\begin{equation}\label{eq:75-local-integral}
 \mathcal T_v(f;s,\tau)=
 \int_{(k_v^\times)^2}f(x,y)w_{\tau,v}(xy)
       \omega_v(x)^{-2}\omega_v(y)^{-1}|x|_v^{3s-1}|y|_v^s\,d^\times x\,d^\times y.
\end{equation}
Here $w_{\tau,v}(t)$ denotes the restriction to $a(1,t)$ of a Jacquet
integral of a fixed finite compact-type family of sections. Meromorphic
coefficients of a section are cleared first. 

The following two analytic statements are recalled from the companion
preprint \cite{LTcuspidal}. Their statements are included because both
are used below for test functions coupling different infinite places.
The first also requires a spectral-family extension, which is proved
here. Its Taylor subtraction is the method of Gel'fand--Shilov
\cite[Chapter~I, Section~3.2]{GS}.

A Schwartz function on the closed positive quadrant is smooth up to
its boundary and has finite seminorms
\[
 q_{L,M}(H)=\max_{i+j\le M}\sup_{r,s\ge0}
 (1+r+s)^L|\partial_r^i\partial_s^jH(r,s)|
 \qquad(L,M\ge0).
\]
Use the analogous convention on the closed half-line.

\begin{lemma}
\label{lem:sector}
Let $h\in C^\infty(\R_{>0})$ have a power bound
$h(t)=O(t^{-B})$ at zero, for some $B\ge0$, and decrease faster
than every power at infinity. Suppose that
\[
 \mathfrak M_h(z)=\int_0^\infty h(t)t^z\frac{dt}{t}
\]
continues meromorphically to $\C$. Fix $e>0$ and
$\lambda_1,\lambda_2\in\R$, and put
\begin{equation}\label{eq:ab}
 a(w)=e(3w-1)+i\lambda_1,\qquad b(w)=ew+i\lambda_2.
\end{equation}
For every Schwartz function $H$ on the closed positive quadrant,
\begin{equation}\label{eq:sector-integral}
 B_H(w)=\int_0^\infty\int_0^\infty
 H(r,s)h(rs)r^{a(w)}s^{b(w)}\frac{dr}{r}\frac{ds}{s}
\end{equation}
is absolutely convergent in a sufficiently far right half-plane
and has meromorphic continuation to $\C$. Locally there is a
nonzero holomorphic clearing factor $p(w)$, independent of $H$,
such that $p(w)B_H(w)$ is holomorphic and, on each compact set $Q$
in that neighborhood,
\begin{equation}\label{eq:sector-uniform-bound}
 \sup_{w\in Q}|p(w)B_H(w)|\le C_Qq_Q(H)
\end{equation}
for a continuous Schwartz seminorm $q_Q$.

The same conclusion holds jointly in $(w,\tau)\in\C\times D$
for $h=h_\tau$, where $D\subset\C$ is open, provided that
$h_\tau(t)$ is holomorphic in $\tau$ for every $t>0$, its power
bound at zero and rapid decay at infinity are uniform on compact
subsets of $D$, and $\mathfrak M_{h_\tau}(z)$ is jointly
meromorphic on $\C\times D$. The clearing factor is then a
nonzero holomorphic function $p(w,\tau)$ independent of $H$,
and the seminorm bound is uniform on compact parameter sets.
\end{lemma}
\begin{proof}
The fixed-function assertion, including the test-independent
clearing factor and the seminorm bound, is
\cite[Lemma~4.3]{LTcuspidal}. We explain why that proof is uniform
in the additional parameter; pointwise continuation alone would
not suffice.

Fix a relatively compact spectral neighborhood $D_0\Subset D$.
The hypotheses give a common exponent $B$ and common constants
there. Choose a bounded connected neighborhood in the $w$-plane
containing the target parameter and meeting a right half-plane
where both $\Re a(w)>B$ and $\Re b(w)>B$. All Taylor orders below
are chosen uniformly on its closure. This ensures agreement with
the original integral on a nonempty open set before continuation
to the target; constructions on overlapping neighborhoods then
agree by the identity theorem. Choose $\kappa\in C_c^\infty([0,1))$ equal to one near
zero. The part of \eqref{eq:sector-integral} multiplied by
$1-\kappa(rs)$ is jointly holomorphic: on its support $rs\ge\delta>0$,
so $r^{-1}\le\delta^{-1}s$ and $s^{-1}\le\delta^{-1}r$,
and Schwartz decay absorbs every parameter power and logarithm.

For the remaining part, split into $r\le s$ and $s\le r$.
On the first sector put $t=rs$ and $R=s$. Its kernel is
\[
 K_H(w,t)=t^{a(w)}\int_{\sqrt t}^\infty
 H(t/R,R)R^{b(w)-a(w)}\frac{dR}{R}.
\]
Put $H_m(R)=\partial_1^mH(0,R)/m!$,
$\beta_{mn}=\partial_1^m\partial_2^nH(0,0)/(m!n!)$, and
$c_m=b-a-m$. Taylor subtraction, exactly as in the cited proof,
gives
\begin{equation}\label{eq:sector-kernel-expansion}
 \begin{split}
 K_H(w,t)={}&\sum_{m<N}\mathcal M(H_m;c_m)t^{a+m}
 -\sum_{\substack{m<N\\n<J}}
 \frac{\beta_{mn}}{c_m+n}t^{(a+b+m+n)/2}
 +\mathcal R_H(w,t),
 \end{split}
\end{equation}
where $\mathcal M(H_m;c)$ is its ordinary scalar Mellin
continuation. On any bounded $w$-neighborhood, choose first $N$
so that
\[
 \Re a+N>B+1,\qquad (\Re(a+b)+N)/2>B+1,
\]
and then $J$ so that, for every $m<N$,
\[
 \Re c_m+J>1,\qquad (\Re(a+b)+m+J)/2>B+1.
\]
The integral remainders in Taylor's formula then give, uniformly,
\[
 |\mathcal R_H(w,t)|\le Cq(H)t^{B+1}(1+|\log t|),
 \qquad 0<t<1.
\]
These choices and the seminorm $q$ are independent of $\tau$ and
$H$. The remainder is holomorphic in $w$, so its integral against
$\kappa(t)h_\tau(t)dt/t$ is jointly holomorphic and seminorm
bounded, by dominated integration. Holomorphy in the separate
parameters and local boundedness also give joint holomorphy.

The truncated Mellin function
$\mathfrak m(z,\tau)=\int_0^1\kappa(t)h_\tau(t)t^zdt/t$
differs from $\mathfrak M_{h_\tau}(z)$ by a jointly entire
function in $z$, holomorphic in $\tau$. Integrating the finite
terms of \eqref{eq:sector-kernel-expansion} therefore gives
\[
 \sum_{m<N}\mathcal M(H_m;c_m)\mathfrak m(a+m,\tau)
 -\sum_{\substack{m<N\\n<J}}
 \frac{\beta_{mn}}{c_m+n}
 \mathfrak m\bigl((a+b+m+n)/2,\tau\bigr).
\]
Each scalar Mellin factor is jointly meromorphic: the two
substitutions in its first variable have nonzero derivative in
$w$. Their finitely many germs, together with the Taylor
 denominators, have a common nonzero holomorphic denominator
independent of $H$. The scalar Taylor-subtraction estimates in
\cite[Lemma~4.3]{LTcuspidal} bound the cleared
$\mathcal M(H_m;c_m)$ by Schwartz seminorms; the same is true of
all finite jets $\beta_{mn}$. This proves the required joint bound
on the first sector. Interchanging $r,s$ and $a,b$ handles the
second. The formulas agree with the convergent integral in a
right half-plane and patch by uniqueness of meromorphic
continuation. Near spectral poles one first clears those poles
and applies the same argument to the holomorphic family.
\end{proof}

\begin{lemma}
\label{lem:tensor-meromorphic}
Let $S$ be a finite set of local fields and let $T_v(\eta)$,
$v\in S$, be meromorphic distribution families on
$\cS(k_v^{d_v})$, with $\eta$ in an open subset of $\C^d$.
Suppose that locally each family has a nonzero holomorphic
clearing factor independent of the test function, and that the
cleared evaluations are holomorphic and satisfy uniform
Schwartz-seminorm bounds on compact parameter sets. At finite
places the bounds are required on each fixed support and
translation level. Then
\[
 F\longmapsto\left(\bigotimes_{v\in S}T_v(\eta)\right)(F)
\]
is a meromorphic distribution on
$\cS(\prod_{v\in S}k_v^{d_v})$, with the same type of uniform
bounds and a test-independent local clearing factor. The
statement applies to arbitrary tests, including tests that are
not decomposable across the places.
\end{lemma}
\begin{proof}
This is \cite[Lemma~5.1]{LTcuspidal}, whose parameter $\eta$
is already allowed to range over $\C^d$; apply it here with
$\eta=(s,\tau)$. The proof uses the Hermite-expansion method for
Schwartz distributions of Simon \cite{Simon}: after clearing the
finitely many local poles, distribution evaluations of Hermite
functions have uniform polynomial growth in their indices,
whereas the Hermite coefficients of a Schwartz test decrease
faster than every power. The resulting series converges normally
in the parameters. Fixed finite support and translation level
supply only a finite-dimensional factor. Thus no extra
factorization assumption on the test function is introduced.
\end{proof}

\begin{lemma}
\label{lem:75-whittaker}
For a holomorphic family of sections in a fixed finite-dimensional
compact-picture space, its Jacquet integral is holomorphic in $\tau$ at
every nonzero Kirillov argument. If $\Omega$ is compact and
$\Gamma=\sup_{\tau\in\Omega}|\Re\tau|$, then, for every $\eta>0$,
\begin{equation}\label{ex:W-estimates}
 |w_{\tau,v}(t)|\ll_{\Omega,v,\eta}|t|_v^{1/2-\Gamma-\eta}
                       \qquad(0<|t|_v\leq1).
\end{equation}
At a finite place its support is bounded above in $|t|_v$, uniformly
on this space. At an infinite place it decreases faster than every
power as $|t|_v\longrightarrow\infty$, uniformly on $\Omega$.

At an unramified place, meromorphically renormalize the spherical
Jacquet family and continue through zeros of the normalizing scalar.
The resulting family has the normalized values
\begin{equation}\label{ex:spherical-values}
 \begin{gathered}
 w_{\tau,v}(\varpi^\ell u)=
 \begin{cases}
 q^{-\ell/2}h_\ell(\alpha,\beta),&\ell\geq0,\quad u\in\OO_v^\times,\\
 0,&\ell<0,
 \end{cases}
 \qquad h_\ell(\alpha,\beta)=\sum_{i=0}^\ell\alpha^i\beta^{\ell-i},\\
 \alpha=\eta_{2,v}(\varpi)q^{-\tau},\quad
 \beta=\eta_{1,v}(\varpi)q^\tau.
 \end{gathered}
\end{equation}
In particular $\alpha\beta=\omega_v(\varpi)$ and
$w_{\tau,v}(1)=1$.

At an infinite place, every angular component $h_{\tau,\eta}(r)$ of
$w_{\tau,v}(r\theta)$ has a scalar Mellin transform
\begin{equation}\label{ex:Bessel-Mellin}
 \mathfrak M_{\eta}(z,\tau)
       =\int_0^\infty h_{\tau,\eta}(r)r^z\frac{dr}{r}
\end{equation}
which continues jointly meromorphically in $(z,\tau)$.
\end{lemma}

\begin{proof}
We give the parameter argument explicitly, since a Mellin theorem for one
representation alone would not prove joint continuation. Use the upper
model with inducing characters $\mu_1,\mu_2$ from
\eqref{eq:75-inducing-characters}. If $b_\tau$ is its section and
$p_\tau(u)=b_\tau(w_0n^+(u))$, the Jacquet integral takes the form
\begin{equation}\label{eq:75-jacquet-line}
 w_{\tau,v}(t)=c_v\mu_{2,v,\tau}(t)|t|_v^{1/2}
          \int_{k_v}p_\tau(u)\psi_v^{\prime}(tu)\,du,
\end{equation}
Outside its half-plane of absolute convergence the integral in
\eqref{eq:75-jacquet-line} is an oscillatory Fourier integral.
At an infinite place it is defined by a smooth partition into a
compact part and dyadic annuli, followed by integration by parts
on the annuli. At a finite place it is the stabilized sum over
compact annuli: for a nonzero Fourier argument all sufficiently
large annuli vanish. The estimates below prove local uniform
convergence in $\tau$ and independence of the smooth partition,
by agreement with the ordinary Jacquet integral in its convergence
half-plane. Here the fixed nontrivial character $\psi_v^{\prime}$
and constant $c_v$ only reflect measure and Weyl-element conventions. This is
\cite[\S3.2.1, equation (3.10)]{MV}. On the chart at infinity of
$\mathbb P^1(k_v)$ the induction law gives
\begin{equation}\label{eq:75-line-chart}
 p_\tau(u)=c(\tau)(\mu_{2,v,\tau}/\mu_{1,v,\tau})(u)
                    |u|_v^{-1}b_\tau^{\infty}(u^{-1}).
\end{equation}
The function $b_\tau^{\infty}$ is smooth near zero, holomorphic in
$\tau$, and locally constant there at a finite place. The harmless
scalar $c(\tau)$ can instead be included in this function.

At a finite place, local constancy makes
$b_\tau^{\infty}(u^{-1})$ constant once $|u|_v$ exceeds a fixed bound.
For a nonzero $t$, the sufficiently large annuli in the oscillatory
integral \eqref{eq:75-jacquet-line} vanish: the multiplicative character
is constant on a fixed relative open subgroup, whereas the additive
character has nonzero frequency there. The cutoff is uniform in $\tau$
because its unit characters and invariance level are fixed. The
remaining finite integrals are holomorphic in $\tau$. Right invariance
under a fixed subgroup $n(\mathfrak p^a)$ gives
$w_{\tau,v}(t)=\psi_v(-tb)w_{\tau,v}(t)$ for all $b\in\mathfrak p^a$;
hence the upper support bound.

For completeness the same argument gives the precise germ at zero.
Extend the homogeneous expression
$(\mu_2/\mu_1)(u)|u|_v^{-1}$ meromorphically by the one-dimensional
Tate subtraction. The difference from $p_\tau$ is supported in a fixed
compact set. Its Fourier transform is constant for sufficiently small
$t$, whereas the homogeneous term transforms into a meromorphic scalar
times $(\mu_1/\mu_2)(t)$. Thus, away from coalescing characters,
\begin{equation}\label{eq:75-finite-germ}
 w_{\tau,v}(t)=|t|_v^{1/2}
       \bigl(A_1(\tau)\mu_{1,v,\tau}(t)
                    +A_2(\tau)\mu_{2,v,\tau}(t)\bigr)
                         \qquad(0<|t|_v\leq q^{-\ell_0}).
\end{equation}
The coefficients are meromorphic and $\ell_0$ is uniform. At coincidence
their sum, rather than its two summands, is continued; this produces at
most a linear factor in the valuation. Equivalently, after division by
$q^{-\ell/2}$, the sequence on each unit coset satisfies the order-two
recurrence with roots $\mu_1(\varpi),\mu_2(\varpi)$, also at coincident
roots. Holomorphic initial values and this recurrence give a bound
$C(1+\ell)q^{\Gamma\ell}$ on compact spectral sets, proving
\eqref{ex:W-estimates}. This also verifies directly, for these families,
the rank-one germ needed for the finite-place continuation below.

At an infinite place let $e=[k_v:\R]$. Formula
\eqref{eq:75-line-chart} and differentiation give symbol bounds of
order $-e(1+2\Re\tau)$ for $p_\tau$ in ordinary Euclidean radius.
They are uniform, with all derivatives, on compact spectral sets.
Split its Fourier integral into a fixed compact part and smooth dyadic
annuli of radius $R\geq1$. An annular term, before oscillation, has
size $O(R^{-2e\Re\tau})$; integration by parts $N$ times adds
$(1+R|t|_{\rm eucl})^{-N}$. Summing first over
$R\leq |t|_{\rm eucl}^{-1}$ and then over larger $R$ gives
\begin{equation}
 |\widehat p_\tau(t)|\ll_{\Omega,\eta}
 |t|_v^{-\eta}\max\{1,|t|_v^{2\Re\tau}\}
                              \qquad(0<|t|_v\leq1).
\end{equation}
The possible logarithm at $\Re\tau=0$ is absorbed by $|t|_v^{-\eta}$.
Multiplication by the prefactor in \eqref{eq:75-jacquet-line} proves
\eqref{ex:W-estimates}. Derivatives of sufficiently high order of
$p_\tau$ are integrable, so repeated integration by parts gives
arbitrary decay for large $t$. The same annular series converges
locally uniformly in $\tau$ for every $t\ne0$, after sufficiently
many integrations by parts, proving its holomorphy. These arguments
apply to real and complex places and do not assume a spherical vector.

To prove the scalar Mellin assertion, use the Fourier--Mellin identity
\cite[\S3.2.1, equation (3.11)]{MV}. Up to a meromorphic local Tate
gamma factor, the Mellin transform
\begin{equation}
 \int_{k_v^\times}w_{\tau,v}(t)\xi(t)|t|_v^{z-1/2}d^\times t
\end{equation}
equals
\begin{equation}\label{eq:75-section-mellin}
 \int_{k_v}p_\tau(u)(\mu_{2,v,\tau}\xi)^{-1}(u)|u|_v^{-z}\,du.
\end{equation}
There is a nonempty joint convergence region near $\Re\tau=0$ in
which this identity has its ordinary meaning. Split
\eqref{eq:75-section-mellin} between the two charts $u$ and $u^{-1}$
of $\mathbb P^1(k_v)$. On each chart its singular endpoint is a smooth
function times a character and a power whose exponent is affine in
$(z,\tau)$. Taylor subtraction at zero gives a jointly meromorphic
continuation: on any compact parameter set finitely many Taylor terms
suffice, the remainder is normally convergent, and the subtracted
monomials have explicit affine denominators. The coefficients and
remainders depend holomorphically on $\tau$. This proves joint
meromorphy, not merely meromorphy in $z$ for a fixed representation.
Projection onto any of the finitely many angular characters and the
substitution $z\mapsto\tfrac12+z/e$ give
\eqref{ex:Bessel-Mellin}. Finally \eqref{ex:spherical-values} is the
$\GL_2$ Casselman--Shalika formula \cite{CS,ShintaniWhittaker}; its
finite sum also proves the asserted normalization at coincident
Satake parameters.
\end{proof}

There is no need to divide a ramified Jacquet function by its value at
$1$: that value can vanish, notably at a reducibility point. Choose a
finite set $\Sigma$ containing all infinite places and all ramified or
nonstandard data. For each global Fourier coefficient define its
\emph{restricted coefficient}
\begin{equation}\label{eq:75-restricted-coefficient}
 W_{\Sigma,j,\tau}(t_\Sigma)
 =W_{j,\tau}\bigl(a(1,(t_\Sigma,1_{\Sigma^c}))\bigr).
\end{equation}
The Jacquet factorization, initially in a convergence region and then
meromorphically, reads
\begin{equation}\label{eq:75-global-factorization}
 W_{j,\tau}(a(1,t))
   =W_{\Sigma,j,\tau}(t_\Sigma)
          \prod_{v\notin\Sigma}w_{\tau,v}(t_v).
\end{equation}
It follows from the usual unfolded Fourier coefficient of an
Eisenstein series; see \cite[\S4.1.7, equation (4.7)]{MV}.
The factor on $\Sigma$ incorporates all global scalar normalizations.
It is holomorphic wherever the corresponding global Eisenstein family
is holomorphic, directly from \eqref{eq:75-restricted-coefficient}.
Locally in $\tau$ it is a finite sum of tensor products of the local
families of Lemma~\ref{lem:75-whittaker}, with meromorphic coefficients.

Lemma~\ref{lem:normalization-removal}, applied to the holomorphic
restricted coefficient, extends the local estimates through removable
factorization poles. The spectral disc may be chosen small enough
that its change in $\sup|\Re\tau|$ is absorbed by the arbitrary
$\eta$ in \eqref{ex:W-estimates}. Genuine spectral poles are
cleared first.

The following lemma is well-known but we include a proof for the convenience of the reader. 

\begin{lemma}
\label{lem:75-cusp-decay}
On compact transverse coordinates of
$g=u a(\underline\lambda)a(t_1,t_2)n(c)$, for every $N>0$,
\begin{equation}\label{eq:75-cusp-decay}
 E_{\rm n}(\tau,g^{-1})\ll_{\Omega,N}\lambda^{-N}
                                  \qquad(\lambda\geq1),
\end{equation}
uniformly on compact sets $\Omega$ of regular spectral parameters.
The estimate holds after locally clearing spectral poles as well.
\end{lemma}

\begin{proof}
Use the finite compact-type Fourier expansion of
\S\ref{subsec:eis-nonconstant}. The finite-place support of
Lemma~\ref{lem:75-whittaker} restricts the nonzero frequency to a fixed
fractional ideal $\mathfrak a$. Its archimedean Kirillov argument has
the form $\underline\lambda^2t\alpha$, with $t$ in a compact subset of
the norm-one ideles. All finite-place factors are bounded by a fixed
power of $1+\|\alpha\|$. Any inverse archimedean powers have the same
property: $|N_{k/\mathbb Q}\alpha|$ is bounded below for
$\alpha\in\mathfrak a\setminus\{0\}$, and the product formula controls
each inverse embedding by the remaining embeddings. At least one
archimedean argument has ordinary size comparable to
$\lambda^{2/n}\|\alpha\|$. Its arbitrary rapid decay absorbs all those
powers, giving, for every $L>0$,
\begin{equation}\label{eq:75-fourier-lattice-majorant}
 |\text{Fourier summand}|\leq
       C_{\Omega,L}(1+\lambda^{2/n}\|\alpha\|)^{-L}.
\end{equation}
Lemma~\ref{lem:uniform-lattice} gives the desired estimate upon
increasing $L$. The preceding restricted-coefficient argument gives
uniformity through removable normalization singularities.
\end{proof}

\begin{lemma}
\label{lem:75-local-holomorphy}
For every place, \eqref{eq:75-local-integral} is absolutely convergent
and jointly holomorphic for a holomorphic Whittaker family on
\begin{equation}\label{ex:local-tube}
 \Re(3s-1)+\tfrac12-|\Re\tau|>0,
 \qquad \Re s+\tfrac12-|\Re\tau|>0.
\end{equation}
It is bounded on compact subsets by a continuous Schwartz seminorm
of $f$. At finite places continuity is understood on each fixed
support and additive-invariance level. For a meromorphic family these
assertions hold after clearing its spectral poles.
\end{lemma}

\begin{proof}
On $|xy|_v\leq1$,
Lemma~\ref{lem:75-whittaker} bounds the integrand by
\begin{equation}
 C|f(x,y)|\,
 |x|_v^{\Re(3s-1)+1/2-\Gamma-\eta}
 |y|_v^{\Re s+1/2-\Gamma-\eta}.
\end{equation}
On a sufficiently small compact parameter box, choose $\eta$ less
than both strict margins in \eqref{ex:local-tube}. Both exponents are
positive. This gives integrability at the axes and their intersection;
Schwartz decay, or compact support at a finite place, handles the
other endpoints. On $|xy|_v\geq1$ every inverse power of one variable
is bounded by a positive power of the other. Schwartz decay and the
large-argument Whittaker bound therefore suffice. These estimates
are uniform and allow all powers of logarithms by slightly shrinking
the strict margins. Dominated holomorphic integration proves the
assertions.
\end{proof}

\begin{proposition}
\label{prop:75-finite}
At a finite place $v$, put $q=q_v$ and
$\lambda_i(\tau)=\mu_{i,v,\tau}(\varpi)$. A denominator independent
of the Schwartz test is
\begin{equation}\label{ex:finite-denominators}
 \Delta_v(s,\tau)=\prod_{i=1}^2
 (1-\lambda_i(\tau)\omega_v(\varpi)^{-2}q^{-3s+1/2})
 (1-\lambda_i(\tau)\omega_v(\varpi)^{-1}q^{-s-1/2}).
\end{equation}
For a holomorphic Jacquet family, $\Delta_v\mathcal T_v$ is a
holomorphic continuous functional on each fixed support and invariance
level, with locally uniform bounds. A meromorphic section adds only
its spectral pole-clearing factor. For a flat section the integral is
rational in $q^{-s},q^\tau$; for a holomorphic linear combination of
flat sections its coefficients have that same holomorphic dependence.
\end{proposition}

\begin{proof}
The finite germ in \eqref{eq:75-finite-germ} implies that, for all
$\ell\geq\ell_0$ and every unit $u$, the holomorphic sequence
$v_\ell(u,\tau)=q^{\ell/2}w_{\tau,v}(\varpi^\ell u)$ satisfies
\begin{equation}\label{eq:75-recurrence}
 v_{\ell+2}=(\lambda_1+\lambda_2)v_{\ell+1}
                                  -\lambda_1\lambda_2v_\ell.
\end{equation}
This identity extends through coincident roots by holomorphy. Its
initial values lie in a fixed finite-dimensional space of locally
constant functions on $\OO_v^\times$.

Decompose the integral into $x=\varpi^m u$, $y=\varpi^n v$. Its test
has a lower support bound $m,n\geq m_0$. Choose $N$ so that it is
constant in either coordinate of valuation at least $N$, and
$N+m_0\geq\ell_0$. Split into the finite rectangle $m,n<N$, two
strips, and the quadrant $m,n\geq N$. Introduce the variables
\begin{equation}\label{eq:75-finite-variables}
 \begin{gathered}
 X=\omega_v(\varpi)^{-2}q^{-3s+1/2},\qquad
 Y=\omega_v(\varpi)^{-1}q^{-s-1/2},\\
 D(T)=(1-\lambda_1T)(1-\lambda_2T).
 \end{gathered}
\end{equation}
The factors $\omega_v(u)^{-2}\omega_v(v)^{-1}$ are retained
in the finite unit integrals. No extension of these characters across
a coordinate axis is made. Each tail sequence, including these unit integrals,
satisfies \eqref{eq:75-recurrence}. For any such sequence $b_\ell$,
\begin{equation}\label{eq:75-shifted-generating}
 S_a(T)=\sum_{j\geq0}b_{a+j}T^j
       =\frac{b_a+(b_{a+1}-(\lambda_1+\lambda_2)b_a)T}{D(T)}.
\end{equation}
Thus the strips have denominators $D(X)$ and $D(Y)$. The quadrant
has generating series
\begin{equation}\label{eq:75-quadrant-generating}
 \sum_{i,j\geq0}b_{a+i+j}X^iY^j
       =\frac{XS_a(X)-YS_a(Y)}{X-Y}.
\end{equation}
After taking the common denominator $D(X)D(Y)$, the numerator is
antisymmetric in $X,Y$ and divisible by $X-Y$. No extra denominator
occurs there, or at $\lambda_1=\lambda_2$. The finite rectangle is
holomorphic by Lemma~\ref{lem:75-whittaker}, proving the claim.
On a fixed support and invariance level there are finitely many
additive-coset basis functions. Bounding their holomorphic coefficient
functions on a compact parameter set proves the asserted uniform
continuity. The denominator does not depend on the chosen level.
\end{proof}

\begin{proposition}
\label{prop:75-infinite}
For $K=\R$ or $\C$, \eqref{eq:75-local-integral} is a jointly
meromorphic continuous functional on $\cS(K^2)$. Locally on a compact
parameter polydisc $\mathcal P$ one nonzero holomorphic function
$D_v(s,\tau)$, independent of the test, clears its poles and satisfies
\begin{equation}\label{eq:75-infinite-seminorm}
 \sup_{(s,\tau)\in\mathcal P}
       |D_v(s,\tau)\mathcal T_v(f;s,\tau)|
                            \leq C_\mathcal P p_\mathcal P(f)
\end{equation}
for a continuous Schwartz seminorm $p_\mathcal P$.
\end{proposition}

\begin{proof}
The fixed compact-type space has finitely many characters on the
compact diagonal subgroup. Write
$w_{\tau,v}(r\theta)=\sum_\eta\eta(\theta)h_{\tau,\eta}(r)$ with $e=[K:\R]$. After angular integration the local functional is the
finite sum of $e^2$ times
\begin{equation}\label{eq:75-arch-radial}
 \int_0^\infty\int_0^\infty
 H_\eta(r,t)h_{\tau,\eta}(rt)
       r^{e(3s-1)-2i\nu_v}t^{es-i\nu_v}\frac{dr}{r}\frac{dt}{t},
\end{equation}
where
\begin{equation}
 H_\eta(r,t)=\int_{\Omega_K^2}f(r\theta,t\nu)
       \eta(\theta\nu)\omega_v(\theta)^{-2}\omega_v(\nu)^{-1}
                                                   d\theta\,d\nu.
\end{equation}
Here $\omega_v(r)=r^{i\nu_v}$ for $r>0$; the fixed real
number $\nu_v$ is a local parameter, even though the global
character is norm-normalized. Each $H_\eta$ is smooth and rapidly
decreasing on the closed quadrant, with seminorms
bounded by those of $f$. Lemma~\ref{lem:75-whittaker} gives all the
family hypotheses of Lemma~\ref{lem:sector}. It gives the asserted meromorphic continuation and seminorm bound
for each angular component. Only finitely many scalar Mellin
functions and Taylor denominators occur on a sufficiently small
parameter polydisc. Their common local holomorphic denominator,
and a spectral denominator if needed, is independent of $f$.
Taking the product over the finite angular set proves the result.
Larger compact sets are covered by finitely many such neighborhoods;
the assertion of local clearing is understood on those neighborhoods.
The exact subtraction formulas agree with the original integral
on \eqref{ex:local-tube}.
\end{proof}

\begin{lemma}
\label{lem:75-coupled}
For any finite $\Sigma$ containing the infinite places, the tensor
product of the local functionals is a jointly meromorphic continuous
functional on $\cS(\prod_{v\in\Sigma}k_v^2)$. Near each parameter
point one pole-clearing factor works for all tests. At finite places
continuity is on each fixed support and invariance level. On
\eqref{ex:local-tube} it agrees with the absolutely convergent
integral, including for nonfactorizable tests.
\end{lemma}

\begin{proof}
Apply Lemma~\ref{lem:tensor-meromorphic} to
Propositions~\ref{prop:75-finite} and \ref{prop:75-infinite}.
Their test-independent denominators and seminorm bounds are exactly
its hypotheses. In \eqref{ex:local-tube}, bound a coupled test by a
product of Schwartz decay weights and finite-place compact support
indicators. Lemma~\ref{lem:75-local-holomorphy} proves absolute
convergence and continuity of the ordinary integral. Fubini gives
agreement on product tests, and their density gives agreement on
all tests. This also proves compatibility of overlapping
meromorphic constructions.
\end{proof}

\begin{proposition}
\label{prop:75-euler}
Choose $\Sigma$ as above, large enough that all $F_j$ are finite sums
of $F_{\Sigma,j}\otimes\bigotimes_{v\notin\Sigma}\one_{\OO_v^2}$.
Then the completed nonconstant contribution continues jointly
meromorphically to $\C^2$ and equals
\begin{equation}\label{eq:eis-euler}
 M^{\rm nc}_{\Phi,\tau}(s)=\frac12
 \frac{L^\Sigma(3s-\tfrac12,\pi_\tau\otimes\omega^{-2})
       L^\Sigma(s+\tfrac12,\pi_\tau\otimes\omega^{-1})}
      {L^\Sigma(4s,\omega^{-2})}\,
                         H_{\Phi,\Sigma}(s,\tau).
\end{equation}
For fixed regular $\tau$, the factor $H_{\Phi,\Sigma}$ is holomorphic
in
\begin{equation}\label{ex:H-tube}
 \Re s>\max\{(\tfrac12+|\Re\tau|)/3,\ |\Re\tau|-\tfrac12\}.
\end{equation}
Equivalently the numerator is the product of the four character
$L$-functions obtained by inserting
\eqref{eq:75-inducing-characters} into the two standard factors.
\end{proposition}

\begin{proof}
At a good place the normalized spherical values
\eqref{ex:spherical-values} give a shell series with
$\alpha=\eta_{2,v}(\varpi)q^{-\tau}$ and
$\beta=\eta_{1,v}(\varpi)q^\tau$. Put $c=\alpha\beta=\omega_v(\varpi)$, $X=c^{-2}q^{-3s+1/2}$ and $Y=c^{-1}q^{-s-1/2}$. The required identity, also used in \cite[Proposition~5.3]{LTcuspidal} and \cite[Lemma~6]{Hough}, is
\begin{equation}\label{eq:generating}
\sum_{m,n\ge0}h_{m+n}(\alpha,\beta)X^mY^n
=\frac{1-cXY}{(1-\alpha X)(1-\beta X)(1-\alpha Y)(1-\beta Y)}.
\end{equation}
It follows by setting $H(T)=((1-\alpha T)(1-\beta T))^{-1}$ and evaluating $(XH(X)-YH(Y))/(X-Y)$. The apparent singularity at $X=Y$ is removable, and the polynomial identity also applies when $\alpha=\beta$. Since $cXY=c^{-2}q^{-4s}$, it yields
\begin{equation}\label{eq:75-unramified-euler-factor}
 \mathcal T_v(\one_{\OO_v^2};s,\tau)=
 \frac{L_v(3s-\tfrac12,\pi_\tau\otimes\omega^{-2})
       L_v(s+\tfrac12,\pi_\tau\otimes\omega^{-1})}
      {L_v(4s,\omega^{-2})}.
\end{equation}

In the initial half-plane of absolute convergence, define the
remaining factor by
\begin{equation}\label{eq:75-H-definition}
 \begin{split}
 H_{\Phi,\Sigma}(s,\tau)=\sum_j
 \int_{(k_\Sigma^\times)^2}
 &F_{\Sigma,j}(x,y)W_{\Sigma,j,\tau}(xy)
\\
 &\hspace{12mm}\cdot\omega_\Sigma(x)^{-2}\omega_\Sigma(y)^{-1}
                    |x|_\Sigma^{3s-1}|y|_\Sigma^s
                                      d^\times x\,d^\times y,
 \end{split}
\end{equation}
where $k_\Sigma=\prod_{v\in\Sigma}k_v$. A repeated index $j$
also includes the finite decomposition of the Schwartz test.
Use the restricted coefficients \eqref{eq:75-restricted-coefficient},
so no division by a possibly vanishing ramified Whittaker value
occurs. The Jacquet factorization, Lemma~\ref{lem:75-coupled}, and
its finite meromorphic linear combinations continue
\eqref{eq:75-H-definition} jointly meromorphically. At a regular
spectral parameter the full restricted coefficients are holomorphic
and satisfy the uniform local bounds proved above. Consequently
\eqref{eq:75-H-definition} itself is normally convergent and
holomorphic in $s$ in \eqref{ex:H-tube}. This argument handles
removable singularities of any auxiliary tensor factorization.

The initial global absolute convergence permits multiplication of
\eqref{eq:75-unramified-euler-factor}, giving \eqref{eq:eis-euler}.
The meromorphic continuation of Hecke $L$-functions and of
$H_{\Phi,\Sigma}$ gives the global continuation. Agreement in the
initial half-plane proves independence of all decompositions and of
$\Sigma$. The region \eqref{ex:H-tube} concerns the finitely many
remaining local distributions; it neither removes numerator poles
nor excludes poles coming from zeros of $L^\Sigma(4s,\omega^{-2})$.
\end{proof}

\begin{proposition}
\label{prop:75-tail}
For every regular $\tau$,
\begin{equation}\label{eq:eis-tail}
 Q_{\Phi,\tau}(s)=\int_1^\infty z^{2s}
                    R^{\rm nc}_{\Phi,\tau}(z)\frac{dz}{z}
\end{equation}
is absolutely convergent for all $s$ and entire in $s$. It is jointly
holomorphic on the regular spectral set and locally meromorphic in
$\tau$ across spectral poles. The lower radial integral continues as
\begin{equation}\label{eq:75-lower-mellin}
 \int_0^1z^{2s}R^{\rm nc}_{\Phi,\tau}(z)\frac{dz}{z}
            =M^{\rm nc}_{\Phi,\tau}(s)-Q_{\Phi,\tau}(s),
\end{equation}
initially as an equality of convergent integrals.
\end{proposition}

\begin{proof}
For the finite compact expansion of the preceding section put
\begin{equation}\label{eq:75-positive-marginal}
 P_{\Phi,j}(x)=\int_\A\int_U
       |c_j(u)|\,|\Phi(u(0,0,x,v))|\,du\,dv,
\end{equation}
and
\begin{equation}\label{eq:75-majorant-factors}
 \mathcal A_{\Phi,j}(L)=\int_{\A^\times}
           P_{\Phi,j}(x)|x|_\A^{2L-1}d^\times x,
 \qquad
 \mathcal B_{\tau,j}(L)=\int_{\A^\times}
           |W_{j,\tau}(a(1,t))||t|_\A^L d^\times t.
\end{equation}
The majorant lemma of \S\ref{subsec:eis-nonconstant} proves their
convergence for $L>\max\{1,\tfrac12+|\Re\tau|\}$. Before Fourier integration in the additive coordinate the change
$x=D$, $t=A/D$ gives
\begin{equation}\label{eq:75-upper-radius-condition}
 r:=|x|_\A^2|t|_\A=|AD|_\A=|\det g|_\A=z^2.
\end{equation}
Thus the same Tonelli argument, restricted to $r\geq1$, bounds the
absolute upper tail by
\begin{equation}\label{eq:75-upper-majorant}
 \frac12\sum_j\int_{r\geq1}
    P_{\Phi,j}(x)|W_{j,\tau}(a(1,t))|
          r^{\Re s}|x|_\A^{-1}d^\times x\,d^\times t.
\end{equation}
On a compact set of $s$ and regular spectral parameters choose
$L>\max\{2+\sup|\Re\tau|,\sup\Re s\}$. Since $r\geq1$, this
is at most
$\tfrac12\sum_j\mathcal A_{\Phi,j}(L)\mathcal B_{\tau,j}(L)$,
uniformly. Any factor $(\log r)^m$ is absorbed by the strict gap
$L-\sup\Re s$, so the integral and its $s$-derivatives converge
normally. Dominated integration gives joint holomorphy, and
Lemma~\ref{lem:normalization-removal} handles removable factorization poles;
at a genuine spectral pole clear the pole first. Splitting the
absolutely convergent complete radial integral at $1$ proves
\eqref{eq:75-lower-mellin}, and \eqref{eq:eis-euler} continues its
right side meromorphically.
\end{proof}
\section{Singular Mellin transforms and the continuation theorem}
\label{subsec:eis-continuation}
The orbit formulas now give the singular Mellin transforms. Poisson
summation is applied first on a compact truncation at a fixed radius;
only after the cusp subtractions cancel is the radius integrated.
This order of operations is the adelic Poisson method of
Shintani and Wright \cite{Shintani,Wright}, with the finite-part
interpretation supplied by the preceding sections.

Write $E_\tau=E(\tau,\eta_1,\eta_2,f)$ and retain the coefficients
$a_\epsilon(\Phi,E_\tau)$ and $d_\epsilon(\Phi,E_\tau)$ from
\S\ref{subsec:eis-constant}. In particular,
\begin{equation}\label{eq:76-orbit-coefficients}
 \begin{aligned}
 a_\epsilon(\Phi,E_\tau)&=\tfrac C3\delta(\eta_{\epsilon,1}^2\eta_{\epsilon,2})
                             Z(F_\epsilon,\eta_{\epsilon,1},u_\epsilon),\\
 d_\epsilon(\Phi,E_\tau)&=C\delta(\eta_{\epsilon,1})
                             Z(f_\epsilon,\eta_{\epsilon,2}^{-1},-p_\epsilon-1).
 \end{aligned}
\end{equation}
Here and below $Z(h,w)$ may also denote the untwisted Tate integral
$Z(h,1,w)$. The generic-parameter assumption is
\begin{equation}\label{eq:76-parameter-conditions}
 E_\tau\text{ and its section data are regular},\qquad
 \tau\notin\{\pm\tfrac12,\pm1,\pm\tfrac32\}.
\end{equation}
The subsequent section removes this restriction by specialization
of the complete formula. The fixed-radius decomposition is
\begin{equation}\label{eq:76-fixed-radius-decomposition}
 R^{\mathrm{reg}}_{\Phi,E_\tau}(z)
 =\sum_{\epsilon=\pm1}
   \bigl(a_\epsilon(\Phi,E_\tau)z^{-u_\epsilon}
                 +d_\epsilon(\Phi,E_\tau)z^{p_\epsilon}\bigr)
       +R^{\mathrm{nc}}_{\Phi,E_\tau}(z).
\end{equation}

\begin{proposition}
\label{prop:76-singular-mellin}
Write $E_\tau^\iota(h)=E_\tau(h^\iota)$; its inducing data are
computed in Lemma~\ref{lem:76-fourier-plane} below. Full family indices
on $Q$ and $M^{\mathrm{nc}}$ specify which Eisenstein series is used.
Let $\gamma=|\Re\tau|$. The integrals
\begin{equation}\label{eq:76-IJ-definition}
 \begin{aligned}
 I^{\mathrm{reg}}(t,\Phi,E_\tau)
   &=\int_0^1z^{2t}R^{\mathrm{reg}}_{\Phi,E_\tau^\iota}(z^{-1})\frac{dz}{z},\\
 J^{\mathrm{reg}}(s,\Phi,E_\tau)
   &=\int_0^1z^{2s}R^{\mathrm{reg}}_{\Phi,E_\tau}(z)\frac{dz}{z}
 \end{aligned}
\end{equation}
converge absolutely for $\Re t>1/2+\gamma$ and $\Re s>2+\gamma$,
respectively. Define
\begin{equation}\label{eq:eis-PI}
 H(t,\Phi,E_\tau)=\sum_{\epsilon=\pm1}
 \left(\frac{a_\epsilon(\Phi,E_\tau)}{2t+u_\epsilon}
       +\frac{d_\epsilon(\Phi,E_\tau)}{2t-p_\epsilon}\right)
\end{equation}
and
\begin{equation}\label{eq:eis-PJ}
 P_J(s,\Phi,E_\tau)=\sum_{\epsilon=\pm1}
 \left(\frac{a_\epsilon(\Phi,E_\tau)}{2s-u_\epsilon}
       +\frac{d_\epsilon(\Phi,E_\tau)}{2s+p_\epsilon}\right)
                   =-H(-s,\Phi,E_\tau).
\end{equation}
Then
\begin{equation}\label{eq:eis-I}
 I^{\mathrm{reg}}(t,\Phi,E_\tau)
                  =H(t,\Phi,E_\tau^\iota)+Q_{\Phi,E_\tau^\iota}(-t),
\end{equation}
\begin{equation}\label{eq:eis-J}
 J^{\mathrm{reg}}(s,\Phi,E_\tau)
  =P_J(s,\Phi,E_\tau)+M^{\mathrm{nc}}_{\Phi,\tau}(s)
                                        -Q_{\Phi,\tau}(s).
\end{equation}
They are jointly meromorphic in their Mellin and spectral parameters.
For a fixed generic $\tau$, the only possible poles of
$I^{\mathrm{reg}}$ are
\begin{equation}\label{eq:76-I-poles}
 t=-u_\epsilon/2=\frac{-1+2\epsilon\tau}{6},
 \qquad t=p_\epsilon/2=\frac12+\epsilon\tau,
\end{equation}
with respective residues $a_\epsilon(\Phi,E_\tau^\iota)/2$ and
$d_\epsilon(\Phi,E_\tau^\iota)/2$ at
distinct locations. If each of
$\eta_1,\eta_2,\eta_1^2\eta_2,\eta_1\eta_2^2$ is nontrivial,
all these coefficients vanish and $I^{\mathrm{reg}}$ is entire.
\end{proposition}
\begin{proof}
The nonconstant upper tail
$Q_{\Phi,\tau}(s)=\int_1^\infty z^{2s}
R^{\mathrm{nc}}_{\Phi,E_\tau}(z)\,dz/z$ is entire by
\S\ref{subsec:eis-local}. Its complete Mellin transform is
$M^{\mathrm{nc}}$ in a sufficiently far right half-plane. For the
nonconstant part of $I$, apply the same result to $E_\tau^\iota$;
the change $r=z^{-1}$ gives $Q_{\Phi,E_\tau^\iota}(-t)$;
for that of $J$, splitting the complete Mellin transform at $1$
gives $M^{\mathrm{nc}}(s)-Q(s)$. For each constant-term power,
\[
 \int_0^1 z^{2t+u_\epsilon}\frac{dz}{z}
       =\frac1{2t+u_\epsilon},\qquad
 \int_0^1 z^{2t-p_\epsilon}\frac{dz}{z}
       =\frac1{2t-p_\epsilon}.
\]
The stated half-plane for $t$ makes both exponents have positive
real part. Similarly the constant-term contributions to $J$ are
$1/(2s-u_\epsilon)$ and $1/(2s+p_\epsilon)$, converging in
the claimed half-plane. The nonconstant Mellin majorants justify
all splittings and give local uniformity. The continuation of
$M^{\mathrm{nc}}$, the entire upper tail, and the meromorphic
Tate coefficients give joint meromorphy. Each displayed denominator
has derivative $2$ in the Mellin variable. This proves the residue
assertion. At coincident locations the residues add; for a fixed
generic parameter no higher-order pole results merely from adding
coincident simple denominators. Finally the selection rules prove
the asserted vanishing criterion.
\end{proof}

\begin{lemma}
\label{lem:76-fourier-plane}
Use the alternating pairing and its self-dual measure
\begin{equation}\label{eq:76-pairing}
 [x,y]=x_1y_4-\tfrac13x_2y_3+\tfrac13x_3y_2-x_4y_1,
 \qquad \widehat\Phi(y)=\int_{V(\A)}\Phi(x)\psi_k([x,y])\,dx.
\end{equation}
For $g^\iota=(\det g)^{-1}g$ one has
\begin{equation}\label{eq:76-fourier-covariance}
 \widehat{\Phi\circ\rho(g)}(y)
        =|\det g|_\A^{-2}\widehat\Phi(g^\iota y).
\end{equation}
The involution preserves $G^1/G_k$, its measure, and the truncations
$\mathcal F_T$. Define the dual family
\begin{equation}\label{eq:76-E-involution}
 E_\tau^\iota(h)=E_\tau(h^\iota)
        =\omega(\det h)^{-1}E_\tau(h).
\end{equation}
It has central character $\omega^{-1}$ and normalized inducing pair
$(\eta_2^{-1},\eta_1^{-1})$, with the same spectral parameter $\tau$.
It is induced by $f_\tau^\iota(h)=\omega(\det h)^{-1}f_\tau(h)$;
its intertwining section is the same determinant twist of
$M_\tau f_\tau$. For the plane integral one has
\begin{equation}\label{eq:76-plane-average}
 A_{\widehat\Phi}(u)=A_\Phi(u^\iota),\qquad
 D_\epsilon(\widehat\Phi,E_\tau^\iota)=D_\epsilon(\Phi,E_\tau).
\end{equation}
\end{lemma}
\begin{proof}
The action satisfies $[gx,gy]=\det(g)[x,y]$ and
$\det\rho(g)=(\det g)^2$. On a diagonal matrix its four weights
are $A^2/D,A,D,D^2/A$; the two paired products are $AD$ and
the total product is $(AD)^2$. Substitution of the unipotent action
\[
 n(c)x=(x_1,x_2+3cx_1,x_3+2cx_2+3c^2x_1,
                 x_4+cx_3+c^2x_2+c^3x_1)
\]
proves the pairing identity there, and direct substitution proves
it for the Weyl element. These elements generate $G$. A change of
variables in the Fourier integral proves
\eqref{eq:76-fourier-covariance}. The map $\iota$ is an involutive
group automorphism; its positive Haar scaling factor therefore has
square one. Heights are unchanged by scalar multiplication, so
it also preserves each rational cusp height and the truncations.
The central transformation law proves \eqref{eq:76-E-involution}.
Multiplying the induction law by $\omega(ad)^{-1}$ gives the
new pair $(\eta_2^{-1},\eta_1^{-1})$. Since $\det w=1$ and
$\det n(v)=1$, the same twist commutes with the intertwining integral
and then with its meromorphic continuation.

For the plane identity first work at one place. Relative to the
product scalar self-dual measure, the self-dual measure of the
pairing is $|3|_v^{-1}dx_1\,dx_2\,dx_3\,dx_4$.
If $P(x_1,x_2)=\int\Phi(x_1,x_2,x_3,x_4)\,dx_3\,dx_4$,
then on the plane
\[
 \widehat\Phi(0,0,y_3,y_4)=|3|_v^{-1}\widehat P(y_4,-y_3/3).
\]
The change of variables has additive measure factor $|3|_v$;
Fourier inversion at zero gives
\[
 \int\widehat\Phi(0,0,y_3,y_4)\,dy_3\,dy_4
 =\int\Phi(0,0,x_3,x_4)\,dx_3\,dx_4.
\]
Tensor products and continuity give the adelic identity.
Apply it to $\Phi\circ\rho(u^\iota)$ and use $|\det u|_\A=1$
to obtain the first equality in \eqref{eq:76-plane-average}.
For the dual compact branches one has
$c_\epsilon^\iota(\tau,u^\iota)=c_\epsilon(\tau,u)$, by the
central law or directly by composing twice with $\iota$.
If the inducing pair is nontrivial, both transverse coefficients
in the second equality of \eqref{eq:76-plane-average} vanish.
For the trivial pair, the equality follows by changing $u$ to
$u^\iota$ in the compact integral. The involution preserves Haar
measure on $U$. This proves the weighted plane identity in every case.
\end{proof}

\begin{proposition}
\label{prop:76-fixed-poisson}
The ordinary nonsingular fixed-radius pairing
\begin{equation}\label{eq:76-N-definition}
 N_{\Phi,E_\tau}(z)=\int_{G^1/G_k}E_\tau(g^{-1})
                                  \Theta_{V',\Phi}(z,g)\,dg^1
\end{equation}
is absolutely convergent for each $z>0$. For generic $\tau$ it
satisfies
\begin{equation}\label{ex:Poisson-fiber}
 \begin{split}
 N_{\Phi,E_\tau}(z)={}z^{-4}N_{\widehat\Phi,E_\tau^\iota}(z^{-1})
 +z^{-4}R^{\mathrm{reg}}_{\widehat\Phi,E_\tau^\iota}(z^{-1})
                      -R^{\mathrm{reg}}_{\Phi,E_\tau}(z).
 \end{split}
\end{equation}
\end{proposition}
\begin{proof}
A nonsingular rational vector has either first coordinate nonzero,
or first coordinate zero and second coordinate nonzero. On a high
cusp these coordinates are expanded by $\lambda^3$ and $\lambda$,
respectively, in adelic norm. The triangular lattice estimate,
applied to these two first-nonzero-coordinate strata, gives arbitrary
cusp decay for fixed $z$; the polynomial Eisenstein bound is
therefore integrable. This also proves uniformity on compact
$z$-intervals and compact regular spectral sets.

Poisson summation at a point of $\mathcal F_T$ reads
\[
 \Theta_{V',\Phi}(z,g)+\Theta_{S,\Phi}(z,g)
 =z^{-4}\bigl(\Theta_{V',\widehat\Phi}(z^{-1},g^\iota)
                    +\Theta_{S,\widehat\Phi}(z^{-1},g^\iota)\bigr).
\]
All sums and integrals here converge on the compact truncation.
Integrate against $E_\tau(g^{-1})$ and change $g$ to $g^\iota$
on the right. By Lemma~\ref{lem:76-fourier-plane} the domain is unchanged
and the transformed weight is $(E_\tau^\iota)(g^{-1})$. The two singular subtraction terms agree:
\begin{equation}\label{eq:76-equal-subtractions}
 z^{-4}B_T(z^{-1},\widehat\Phi,E_\tau^\iota)=B_T(z,\Phi,E_\tau),
\end{equation}
since each is
$z^{-2}\sum_\epsilon D_\epsilon(\Phi,E_\tau)
T^{p_\epsilon+2}/(p_\epsilon+2)$.
Subtract them before sending $T$ to infinity. The nonsingular
terms converge ordinarily by the first paragraph; the singular
terms converge by the common finite-part theorem. This proves
\eqref{ex:Poisson-fiber}.
\end{proof}

\begin{proposition}\label{prop:76-Z-split}
In the initial half-plane $\Re s>5/2+\gamma$,
\begin{equation}\label{ex:Z-split}
 \begin{split}
 Z(s,\Phi,E_\tau)={}&Z_+(s,\Phi,E_\tau)
                         +Z_+(2-s,\widehat\Phi,E_\tau^\iota)\\
 &+I^{\mathrm{reg}}(s-2,\widehat\Phi,E_\tau)
                         -J^{\mathrm{reg}}(s,\Phi,E_\tau).
 \end{split}
\end{equation}
The upper radial zeta terms are entire in their first argument.
\end{proposition}
\begin{proof}
Absolute convergence of the original integral and entireness of
its upper part were proved in \S\ref{subsec:eis-family}.
Multiply \eqref{ex:Poisson-fiber} by $z^{2s}$ and integrate on
$(0,1)$. The two singular integrals converge absolutely because
$\Re(s-2)>1/2+\gamma$ and $\Re s>2+\gamma$.
The dual nonsingular term is the upper integral after $r=z^{-1}$:
\[
 \int_0^1 z^{2s-4}N_{\widehat\Phi,E_\tau^\iota}(z^{-1})\frac{dz}{z}
 =\int_1^\infty r^{2(2-s)}N_{\widehat\Phi,E_\tau^\iota}(r)\frac{dr}{r}.
\]
Adding the original upper part proves the formula. All operations
are on the already defined fixed-radius finite parts; no radial
integral is interchanged with a divergent cusp integral.
\end{proof}

\begin{theorem}
\label{thm:eisenstein}
For a generic $\tau$ the ordinary nonsingular zeta integral has the
jointly meromorphic continuation
\begin{equation}\label{eq:eis-Z}
 \begin{split}
 Z(s,\Phi,E_\tau)={}&Z_+(s,\Phi,E_\tau)
                  +Z_+(2-s,\widehat\Phi,E_\tau^\iota)\\
 &+H(s-2,\widehat\Phi,E_\tau^\iota)-P_J(s,\Phi,E_\tau)\\
 &+Q_{\widehat\Phi,E_\tau^\iota}(2-s)+Q_{\Phi,E_\tau}(s)
                              -M^{\mathrm{nc}}_{\Phi,E_\tau}(s).
 \end{split}
\end{equation}
The notation $Q_{\Phi,E_\tau}$ and $M^{\mathrm{nc}}_{\Phi,E_\tau}$
records the full family; when it is fixed these are the earlier
$Q_{\Phi,\tau}$ and $M^{\mathrm{nc}}_{\Phi,\tau}$.
For $\gamma=|\Re\tau|$, the difference
$Z(s,\Phi,E_\tau)-H(s-2,\widehat\Phi,E_\tau^\iota)$ is holomorphic
in $\Re s>1/2+\gamma$. The only possible poles of $Z$ in that half-plane are
\begin{equation}\label{eq:76-Z-poles}
 s=\frac{11}{6}+\frac{\epsilon\tau}{3},\qquad
 s=\frac52+\epsilon\tau,\qquad \epsilon\in\{+1,-1\},
\end{equation}
provided they lie in the indicated half-plane. At distinct locations
their residues are, respectively,
\begin{equation}\label{eq:76-Z-residues}
 \tfrac12a_\epsilon(\widehat\Phi,E_\tau^\iota),\qquad
 \tfrac12d_\epsilon(\widehat\Phi,E_\tau^\iota).
\end{equation}
In the plus branch these residues can be nonzero only when
$\eta_1\eta_2^2=1$ and $\eta_2=1$, respectively; in the minus
branch the conditions are $\eta_1^2\eta_2=1$ and $\eta_1=1$.

Suppose, more strongly, that
\begin{equation}\label{eq:eis-character-nonresonance}
 \eta_1\ne1,\qquad\eta_2\ne1,\qquad
 \eta_1^2\eta_2\ne1,\qquad\eta_1\eta_2^2\ne1.
\end{equation}
Then $I^{\mathrm{reg}}(t,\Phi,E_\tau)$ is entire in $t$, and both
$J^{\mathrm{reg}}(s,\Phi,E_\tau)$ and $Z(s,\Phi,E_\tau)$ are
holomorphic on
\begin{equation}\label{eq:eis-nontrivial-character-region}
 \Re s>\beta(\tau),\qquad
 \beta(\tau)=\max\left\{\frac14,
                      \frac{1/2+|\Re\tau|}{3},
                              |\Re\tau|-\frac12\right\}.
\end{equation}
In particular this is $\Re s>1/4$ on the unitary spectral axis.
For an arbitrary unitary norm component, replace $s$ by
$s-i\kappa/2$ as in Lemma~\ref{lem:central-norm-normalization}.
\end{theorem}
\begin{proof}
Insert \eqref{eq:eis-I} and \eqref{eq:eis-J} into the dual-family
Poisson identity \eqref{ex:Z-split}. All terms are jointly meromorphic
and agree with the ordinary integral in its initial domain. The two
upper zeta integrals and two upper nonconstant tails are entire in $s$.
The poles of $P_J$ have real parts at most $(1+2\gamma)/6$ and
$-1/2+\gamma$, both less than $1/2+\gamma$.
The four numerator factors of \eqref{eq:eis-euler} are
\begin{equation}\label{eq:76-four-Hecke-factors}
 \begin{gathered}
 L^\Sigma(3s-\tfrac12+\tau,(\eta_1^2\eta_2)^{-1}),\qquad
 L^\Sigma(3s-\tfrac12-\tau,(\eta_1\eta_2^2)^{-1}),\\
 L^\Sigma(s+\tfrac12+\tau,\eta_1^{-1}),\qquad
 L^\Sigma(s+\tfrac12-\tau,\eta_2^{-1}).
 \end{gathered}
\end{equation}
Because the characters are norm-normalized, a numerator pole can
occur only for a trivial displayed character, at argument $1$.
All such poles have real part at most $1/2+\gamma$.
The denominator $L^\Sigma(4s,\omega^{-2})$ has no zeros in
$\Re s>1/4$. Its possible pole at $4s=1$ produces a zero,
not a pole, in the quotient. The remaining local factor is holomorphic
when $3\Re s-1/2>\gamma$ and $\Re s+1/2>\gamma$, both of
which hold in $\Re s>1/2+\gamma$. This proves the assertion
about the difference. Substituting $t=s-2$ into
\eqref{eq:76-I-poles} gives the locations and residues. Applying
the orbit selection rules to the dual pair
$(\eta_2^{-1},\eta_1^{-1})$ gives the stated branchwise conditions.

Under \eqref{eq:eis-character-nonresonance}, the coefficients of
$H$ and $P_J$ vanish for both the original and dual families, and
all four factors in \eqref{eq:76-four-Hecke-factors} are entire.
Intersecting the zero-free denominator region with the two strict
local inequalities gives exactly \eqref{eq:eis-nontrivial-character-region}.
The formulas for $I,J,Z$ now prove the claims. These conclusions use
the meromorphic Euler continuation, without asserting convergence
of the original Whittaker integral in the enlarged domain.
\end{proof}

\begin{corollary}
\label{cor:76-singular-vanishing}
Suppose $\Phi|_{S(\A)}=0$. At a generic parameter the only possible
poles of $Z(s,\Phi,E_\tau)$ on the whole plane are the locations in
\eqref{eq:76-Z-poles}, with the branchwise selection rules of
Theorem~\ref{thm:eisenstein}. Under
\eqref{eq:eis-character-nonresonance}, this zeta function is entire.
\end{corollary}
\begin{proof}
Every translate of a singular vector is singular, so the repeated-root
function for $\Phi$ vanishes. Its weighted tests, constant coefficients,
nonconstant pairing, $M^{\mathrm{nc}}$, and $Q$ vanish as well.
Formula \eqref{eq:eis-Z} therefore reduces to entire upper terms
and $H(s-2,\widehat\Phi,E_\tau^\iota)$. The latter has only the
listed possible poles and vanishes under the stated character condition.
\end{proof}
\section{Exceptional spectral parameters}
\label{subsec:eis-exceptional}
The exceptional parameters of a particular section include its genuine
Eisenstein poles, and may also include apparent poles introduced by
an orbit formula or a local normalization. This section treats
them through the complete family. A cutoff finite part is the
constant coefficient of a power--logarithm expansion in the height
$T$. A spectral Laurent coefficient is a coefficient in
$\tau-\tau_0$. These are different operations.

\begin{proposition}
\label{prop:77-full-specialization}
Suppose $E_\tau$ is holomorphic near $\tau_0$. If
$q_\epsilon(\tau_0)=3+2\epsilon\tau_0\ne0$ for both signs,
then $R^{\mathrm{reg}}_{\Phi,E_\tau}(z)$ is holomorphic there
and specializes to the direct fixed-parameter cutoff finite part.
The same conclusion holds without the condition on $q_\epsilon$
when $(\eta_1,\eta_2)\ne(1,1)$. All statements are locally uniform in $z>0$.
\end{proposition}
\begin{proof}
The constant term is its compact unipotent average and is therefore
holomorphic whenever $E_\tau$ is. For a fixed large threshold
$T_0$, the finite-part construction writes
\begin{equation}\label{eq:77-fixed-threshold-formula}
 R^{\mathrm{reg}}_{\Phi,E_\tau}(z)
 =\mathcal K_{\Phi}(z,\tau)
   -z^{-2}\sum_{\epsilon=\pm1}D_\epsilon(\Phi,E_\tau)
                     \frac{T_0^{q_\epsilon(\tau)}}{q_\epsilon(\tau)}.
\end{equation}
Here $\mathcal K$ is the sum of the compact integral over
$\mathcal F_{T_0}$ and the integral on the high cusp after
subtracting the full leading integrand. To check holomorphy, write
$\Theta_S=z^{-2}\lambda^4A_\Phi(u)+r_\Theta$, with rapidly
decreasing $r_\Theta$. The remainder integrand is
\[
 E_\tau(g^{-1})r_\Theta(g)
    +z^{-2}\lambda^4A_\Phi(u)E_{\tau,\mathrm n}(g^{-1}).
\]
The first term is rapidly decreasing by moderate growth and the
second by the Fourier estimate. Both bounds are uniform on a slightly
larger compact spectral set; they also bound derivatives by Cauchy's
formula. Thus $\mathcal K$ is holomorphic. When both exponents are
nonzero, the displayed primitives are holomorphic and coincide with
the direct cutoff prescription. When $(\eta_1,\eta_2)\ne(1,1)$, every
$D_\epsilon$ is identically zero by character orthogonality, so
there is no resonant primitive at all. This also proves the claimed
uniformity.
\end{proof}

\begin{proposition}
\label{prop:77-tate-cancellation}
Assume $\eta_1=\eta_2=1$ and $E_\tau$ is regular at $\tau_0=-\epsilon$.
Let $F_\epsilon^0$ and $f_\epsilon^0$ denote the weighted tests
at $\tau_0$, and put $B_\epsilon^0=f_\epsilon^0(0)
=\widehat F_\epsilon^0(0)$. Write
\begin{equation}\label{eq:77-Tate-finite-coefficients}
 Z(F_\epsilon^0,1,w)=\frac{CB_\epsilon^0}{w-1}+\alpha_\epsilon+O(w-1),
 \qquad
 Z(f_\epsilon^0,1,w)=-\frac{CB_\epsilon^0}{w}+\beta_\epsilon+O(w).
\end{equation}
The sum of the corresponding triple- and double-root contributions
has the regular limit
\begin{equation}\label{eq:77-logarithmic-pair}
 z^{-1}(D_{\epsilon,\Phi}+L_{\epsilon,\Phi}\log z),
 \qquad D_{\epsilon,\Phi}=\tfrac C3\alpha_\epsilon+C\beta_\epsilon,
 \quad L_{\epsilon,\Phi}=\tfrac23C^2B_\epsilon^0.
\end{equation}
Its contributions to $I^{\mathrm{reg}}$ and $J^{\mathrm{reg}}$
are respectively
\begin{equation}\label{eq:77-logarithmic-Mellin}
 \frac{D_{\epsilon,\Phi}}{2t+1}
       +\frac{L_{\epsilon,\Phi}}{(2t+1)^2},\qquad
 \frac{D_{\epsilon,\Phi}}{2s-1}
       -\frac{L_{\epsilon,\Phi}}{(2s-1)^2}.
\end{equation}
Thus a regular spectral specialization can produce a double Mellin
pole. These Tate endpoint residues vanish identically unless
$\eta_1=\eta_2=1$; for other inducing pairs there is no such
Tate-pole cancellation to perform.
\end{proposition}
\begin{proof}
Put $\delta=\epsilon\tau+1$. Then
$u_\epsilon=1-2\delta/3$, $p_\epsilon=-1+2\delta$, and
$v_\epsilon=-2\delta$. The common endpoint coefficient
$B_\epsilon(\tau)$ is holomorphic, and the sum of the two
polar parts is exactly
\[
 \frac{C^2B_\epsilon(\tau)}{2\delta}
           z^{-1}(z^{2\delta}-z^{2\delta/3}).
\]
Its limit is $(2/3)C^2B_\epsilon^0z^{-1}\log z$.
The regular parts of the two Tate integrals tend to the displayed
finite coefficients. The derivative of $B_\epsilon$ contributes
nothing because the two residues have canceled before the limit.
This remains valid for a varying nonspherical section: the endpoint
identity holds as an identity of holomorphic families. Finally
$\int_0^1z^a\log z\,dz/z=-a^{-2}$, and substitution of
$z^{-1}$ changes the sign of the logarithm in the $I$ integral.
These facts prove \eqref{eq:77-logarithmic-Mellin}.
\end{proof}

\begin{lemma}
\label{lem:77-moving-exponent}
Let $a,b$ be holomorphic near $\tau_0$, with $b(\tau_0)=0$
and $b'(\tau_0)\ne0$. For fixed $T>0$,
\begin{equation}\label{eq:77-moving-exponent-residue}
 \operatorname*{Res}_{\tau=\tau_0}\frac{a(\tau)T^{b(\tau)}}{b(\tau)}
                           =\frac{a(\tau_0)}{b'(\tau_0)},
\end{equation}
and
\begin{equation}\label{eq:77-moving-exponent-constant}
 \operatorname*{CT}_{\tau=\tau_0}\frac{a(\tau)T^{b(\tau)}}{b(\tau)}
 =a(\tau_0)\log T+\frac{a'(\tau_0)}{b'(\tau_0)}
             -\frac{a(\tau_0)b''(\tau_0)}{2b'(\tau_0)^2}.
\end{equation}
\end{lemma}
\begin{proof}
For $\delta=\tau-\tau_0$ write
$a=a_0+a_1\delta+O(\delta^2)$ and
$b=b_1\delta+b_2\delta^2/2+O(\delta^3)$.
The expansions
$1/b=1/(b_1\delta)-b_2/(2b_1^2)+O(\delta)$ and
$T^b=1+b_1\delta\log T+O(\delta^2)$ give both coefficients
by multiplication.
\end{proof}

\begin{proposition}
\label{prop:77-resonant-comparison}
Suppose $E_\tau$ is regular at $\tau_0=\pm3/2$, and choose
$\epsilon$ with $q_\epsilon(\tau_0)=0$. Let
$R^{\log}_{\Phi,E_{\tau_0}}(z)$ denote the direct cutoff finite
part using $z^{-2}D_\epsilon(\Phi,E_{\tau_0})\log T$ for
that exponent and the ordinary primitive for the other sign. Then
\begin{equation}\label{eq:77-resonant-residue}
 \operatorname*{Res}_{\tau=\tau_0}
 R^{\mathrm{reg}}_{\Phi,E_\tau}(z)
 =-\frac{z^{-2}D_\epsilon(\Phi,E_{\tau_0})}{2\epsilon},
\end{equation}
\begin{equation}\label{eq:77-resonant-constant}
 \operatorname*{CT}_{\tau=\tau_0}
 R^{\mathrm{reg}}_{\Phi,E_\tau}(z)
 =R^{\log}_{\Phi,E_{\tau_0}}(z)
   -\frac{z^{-2}}{2\epsilon}
       \left.\partial_\tau D_\epsilon(\Phi,E_\tau)\right|_{\tau_0}.
\end{equation}
For $(\eta_1,\eta_2)\ne(1,1)$ both corrections are zero.
\end{proposition}
\begin{proof}
In \eqref{eq:77-fixed-threshold-formula} the compact-plus-remainder
term is holomorphic; the only possibly singular primitive is the
negative of $a(\tau)T_0^{b(\tau)}/b(\tau)$, where
$a=z^{-2}D_\epsilon$ and $b=3+2\epsilon\tau$.
Here $b'=2\epsilon$ and $b''=0$. Apply
Lemma~\ref{lem:77-moving-exponent}. At the fixed parameter,
integrating the height density gives $a(\tau_0)\log T$.
Comparison with that direct prescription proves both formulas.
The character selection rule gives the last assertion.
\end{proof}

\begin{lemma}
\label{lem:77-resonant-Poisson}
The spectral residues, and the differences between the two constant
values in \eqref{eq:77-resonant-constant}, cancel between
$I^{\mathrm{reg}}(s-2,\widehat\Phi,E_\tau)$ and
$J^{\mathrm{reg}}(s,\Phi,E_\tau)$ in \eqref{ex:Z-split}.
\end{lemma}
\begin{proof}
Either correction in the fixed-radius pairing is $K_\Phi z^{-2}$,
where $K_{\widehat\Phi,E_\tau^\iota}=K_{\Phi,E_\tau}$
by the weighted plane identity
\eqref{eq:76-plane-average} and its spectral derivative. Its two
radial contributions are respectively
\[
 K_{\widehat\Phi,E_\tau^\iota}\int_0^1z^{2(s-2)+2}\frac{dz}{z}
       =\frac{K_{\widehat\Phi,E_\tau^\iota}}{2s-2},\qquad
 K_{\Phi,E_\tau}\int_0^1z^{2s-2}\frac{dz}{z}=\frac{K_{\Phi,E_\tau}}{2s-2}.
\]
They are equal in a right half-plane and hence meromorphically
everywhere. Spectral extraction and radial integration commute
there by Cauchy's formula on a small spectral circle and the
uniform radial majorants. The opposite signs in Poisson summation
prove cancellation.
\end{proof}

\begin{proposition}
\label{prop:77-constant-specialization}
Let $E^0_\tau$ be the standard spherical family with $\eta_1=\eta_2=1$
and section equal to $1$ on $U$. Put
\begin{equation}\label{eq:77-residue-normalization}
 R_k=\operatorname*{Res}_{w=1}\Lambda_k(w),\qquad
 \rho=\frac{R_k}{2\Lambda_k(2)}.
\end{equation}
Then
\begin{equation}\label{eq:77-E-constant}
 E^0_{-1/2}=1,\qquad
 \operatorname*{Res}_{\tau=1/2}E^0_\tau=\rho.
\end{equation}
For the measures of this paper,
\begin{equation}\label{eq:77-volume}
 V_k:=\vol(G^1/G_k)=\frac{C^2\Lambda_k(2)}{R_k},\qquad
 \vol(\mathcal F_T)=V_k-\frac{C^2}{2T^2}\quad(T\text{ large}).
\end{equation}
Writing $\widetilde\Phi=\int_U\Phi(u\,\cdot)\,du$,
$F_\Phi(t)=\widetilde\Phi(0,0,0,t)$ and
$f_\Phi(t)=\int_\A\widetilde\Phi(0,0,t,v)\,dv$, one has
\begin{equation}\label{eq:77-full-constant-pairing}
 R^{\mathrm{reg}}_{\Phi,E^0_{-1/2}}(z)
 =V_k\Phi(0)+\frac C3z^{-2/3}Z(F_\Phi,1,2/3)
                           +CZ(f_\Phi,1,-1).
\end{equation}
More generally, let $\eta$ be a quadratic Hecke character trivial
on the norm section and take the section whose compact restriction
is $u\mapsto\eta(\det u)$. For the inducing pair $(\eta_1,\eta_2)=(\eta,\eta)$,
\begin{equation}\label{eq:77-quadratic-family}
 E(\tau,\eta,f;h)=\eta(\det h)E^0_\tau(h).
\end{equation}
Its value at $-1/2$ is $\eta\circ\det$ and its residue at $1/2$
is $\rho\,\eta\circ\det$. If $\eta\ne1$, the corresponding
zeta integral twisted by $\eta\circ\det$ is entire.
The constant identities \eqref{eq:77-E-constant} are asserted for
the stated spherical section, not for an arbitrary compact type.
\end{proposition}
\begin{proof}
For the standard spherical section the functional equation and
constant-residue theorem \cite[Section~16]{JL}, with the explicit
residue normalization in \cite[Section~4.1.7, equation~(4.6)]{MV}, give
$E^0_\tau=m(\tau)E^0_{-\tau}$ and the stated residue, where
$m(\tau)=\Lambda_k(2\tau)/\Lambda_k(1+2\tau)$.
The functional equation of $\Lambda_k$ gives, for
$\delta=\tau+1/2$,
\[
 m(-\tfrac12+\delta)
   =-\frac{2\Lambda_k(2)}{R_k}\delta+O(\delta^2),\qquad
 E^0_{1/2-\delta}=-\frac{\rho}{\delta}+O(1).
\]
Their product proves $E^0_{-1/2}=1$. Specialize the finite-cutoff
zero-average identity from \S\ref{subsec:eis-finite-part}:
\[
 \int_{\mathcal F_T}E^0_\tau(g^{-1})\,dg^1
 =C^2\left(\frac{T^{2\tau-1}}{2\tau-1}
       -m(\tau)\frac{T^{-1-2\tau}}{1+2\tau}\right).
\]
At $\tau=-1/2$ its first term is $-C^2/(2T^2)$ and the second
has limit $C^2\Lambda_k(2)/R_k$. This proves the volume formula.
The apparently missing zero-orbit term in the generic decomposition
is recovered by the product
\[
 \lim_{\tau\to-1/2}\frac C3m(\tau)z^{-u_-}
                                      Z(F_\Phi,1,u_-)=V_k\Phi(0),
\]
since $u_-=2\delta/3$ and the Tate residue at zero is
$-C\Phi(0)$. The other minus-branch term tends to zero, and the
plus-branch terms give the two remaining summands in
\eqref{eq:77-full-constant-pairing}. The nonconstant part is zero
at this spherical specialization. Proposition
\ref{prop:77-full-specialization} justifies specialization of the
complete singular pairing, whose cusp exponents are nonzero.

When $\eta^2=1$, multiplication by $\eta(\det h)$ changes the
spherical inducing factor into that of $(\eta,\eta)$; the character
is trivial on rational matrices. This proves
\eqref{eq:77-quadratic-family} directly in the convergent series
and then meromorphically. If $\eta\ne1$, all singular
constant coefficients vanish and the nonconstant part vanishes at
$\tau=-1/2$, because the specialization is a determinant
character. Thus $R^{\mathrm{reg}}=0$ at that parameter, and the
Poisson formula consists of two entire upper zeta integrals. This
proves entireness without any assertion about arbitrary compact-type
values at the reducibility point.
\end{proof}

\begin{remark}\label{rem:rational-volume}
For $k=\mathbb Q$ the measures of \S\ref{subsec:measures} give
$C=1/2$: the compact group
$\{\pm1\}\times\prod_p\mathbb Z_p^\times$ has mass one and
maps onto $\A^1/\mathbb Q^\times$ with kernel $\{\pm1\}$,
which has counting mass two. Since $R_{\mathbb Q}=1$ and
$\Lambda_{\mathbb Q}(2)=\pi/6$, formula \eqref{eq:77-volume}
gives $V_{\mathbb Q}=\pi/24$ and cusp volume $1/(8T^2)$.
This also follows from the classical modular fundamental domain.
After the finite compact factors are removed, $y=\lambda^2$
turns $\lambda^{-2}d\lambda/\lambda$ into $dy/(2y^2)$.
The rational central stabilizer $\{\pm I\}$ and the
determinant-negative component of $\mathrm{GL}_2(\mathbb Z)$
each contribute a factor $1/2$. Thus our quotient measure is
$1/8$ of the hyperbolic measure on the standard
$\mathrm{PSL}_2(\mathbb Z)$ domain, whose area is $\pi/3$.
\end{remark}

The following specialization recovers the quasi-character-twisted
binary-cubic theory of Wright \cite{Wright}; it is a
consistency check on the Eisenstein calculation, rather than a new
determinant-character theorem. The local and arithmetic aspects
are developed in Datskovsky--Wright \cite{DW2,DW3}, and the cubic
character phenomenon also appears in the orbital $L$-functions
of Taniguchi--Thorne \cite{TT1,TT2}.

\begin{proposition}
\label{prop:77-determinant-characters}
Let $\eta$ be any unitary Hecke character trivial on the norm section.
For the inducing pair $(\eta,\eta)$ and compact restriction
$f_\tau(u)=\eta(\det u)$ one has
\begin{equation}\label{eq:77-general-determinant-family}
 E_\tau(h)=\eta(\det h)E^0_\tau(h),\qquad
 E_{-1/2}(h)=\eta(\det h).
\end{equation}
If $\eta\ne1$, put
\begin{equation}\label{eq:77-determinant-coefficient}
 \begin{aligned}
 F_{\Phi,\eta}(t)&=\int_U\eta(\det u)^{-1}\Phi(u(0,0,0,t))\,du,\\
 A_{\Phi,\eta}&=\frac C3\delta(\eta^3)
                          Z(F_{\Phi,\eta},\eta,2/3).
 \end{aligned}
\end{equation}
Then
\begin{equation}\label{eq:77-determinant-pairing}
 R^{\mathrm{reg}}_{\Phi,\eta\circ\det}(z)
                         =A_{\Phi,\eta}z^{-2/3}.
\end{equation}
Consequently the determinant-character zeta integral is entire if
$\eta^3\ne1$. For a nontrivial cubic character its only possible
poles on the whole plane are $s=1/3$ and $s=5/3$, and
\begin{equation}\label{eq:77-determinant-zeta}
 \begin{split}
 Z(s,\Phi,\eta\circ\det)
 ={}&Z_+(s,\Phi,\eta\circ\det)
     +Z_+(2-s,\widehat\Phi,\eta^{-1}\circ\det)\\
 &+\frac{A_{\widehat\Phi,\eta^{-1}}}{2s-10/3}
    -\frac{A_{\Phi,\eta}}{2s-2/3}.
 \end{split}
\end{equation}
In particular the possible pole in $\Re s>1$ is $s=5/3$.
\end{proposition}
\begin{proof}
Multiplication of the spherical section by $\eta\circ\det$ gives
exactly the induction law for $(\eta,\eta)$, and factors through
the rational summation because $\eta$ is trivial on $k^\times$.
This proves the family identity first in its convergence half-plane
and then meromorphically. Its value at $-1/2$ follows from
\eqref{eq:77-E-constant}. The central character is $\eta^2$.

For $\eta\ne1$, both double-root coefficients vanish. The
triple-root condition is $\delta(\eta^3)$ in both branches. The
plus branch at $\tau=-1/2$ has $u_+=2/3$ and compact coefficient
$\eta(\det u)^{-1}$, giving \eqref{eq:77-determinant-coefficient}.
The minus compact coefficient is $m(\tau)\eta(\det u)^{-1}$,
which tends to zero there; its Tate integral has character $\eta$,
and is entire, so its limiting contribution is zero. The full
nonconstant coefficient vanishes because a determinant character is
invariant under the unipotent subgroup. The transverse cusp
subtractions and zero-orbit averages vanish since $\eta\ne1$.
The specialization theorem away from cusp resonance therefore gives
\eqref{eq:77-determinant-pairing}. Apply fixed-radius Poisson
summation to this expression and to the inverse character, and
integrate. The two elementary radial integrals are precisely the
fractions in \eqref{eq:77-determinant-zeta}; the upper integrals
are entire. If $\eta^3\ne1$ both coefficients vanish. This
also recovers the nontrivial quadratic case above.
\end{proof}

\begin{theorem}
\label{thm:77-specialization}
At every regular parameter $\tau_0$, the ordinary nonsingular zeta
integral has a meromorphic continuation in $s$ obtained by combining
the full expression \eqref{eq:eis-Z} and then restricting to
$\tau=\tau_0$. The combined expression has no polar component
$\tau=\tau_0$. Moving Mellin poles may remain or coalesce.

If near a spectral pole
\begin{equation}\label{eq:77-E-Laurent}
 E_\tau(h)=\sum_{j=-r}^{\infty}E_j(h)(\tau-\tau_0)^j,
\end{equation}
then the zeta integral twisted by each $E_j$ converges in a sufficiently
far right half-plane and extends meromorphically in $s$. It is the
corresponding Laurent coefficient of the complete formula
\eqref{eq:eis-Z}, interpreted meromorphically in $s$. This includes
residues and all spectral derivatives at regular parameters.
Under \eqref{eq:eis-character-nonresonance}, the entire assertion
for $I^{\mathrm{reg}}$
and the holomorphy assertion \eqref{eq:eis-nontrivial-character-region}
hold at every regular parameter. They also hold for each spectral
Laurent coefficient, with $\beta(\tau_0)$ in place of $\beta(\tau)$.
\end{theorem}
\begin{proof}
On a small closed disc where $E_\tau$ is holomorphic, choose
$\Gamma>\sup|\Re\tau|$. The nonsingular theta majorants give
uniform absolute convergence of its defining integral for
$\Re s>5/2+\Gamma$. In this common half-plane the zeta integral
is holomorphic in $\tau$. The complete right side of
\eqref{eq:eis-Z} is jointly meromorphic and agrees there with
that integral. At a generic value of $s$, its negative Laurent
coefficients along $\tau=\tau_0$ are meromorphic in $s$ and
vanish on the initial half-plane. The identity theorem makes them
identically zero. Thus there is no vertical polar component, and
restriction to $\tau_0$ is well defined as a meromorphic function
of $s$. At a moving Mellin pole its numerical value need not be
finite; no such assertion is made.

For \eqref{eq:77-E-Laurent}, choose a small circle containing no
other spectral pole. Uniform bounds on this circle and Cauchy's
coefficient formula permit Fubini in a common initial half-plane,
giving exactly
\[
 Z(s,\Phi,E_j)=\int_{G_\A/G_k}|\det g|_\A^s E_j(g^{-1})
                                \sum_{x\in V'_k}\Phi(gx)\,dg.
\]
The same circle bounds $E_j$ by a constant times the uniform
moderate-growth bound and proves absolute convergence. Its coefficient
in the continued formula gives meromorphic continuation. Coefficients
are extracted locally at generic $s$ and then continued; a single
fixed spectral contour need not avoid every moving Mellin pole.

Finally, under \eqref{eq:eis-character-nonresonance}, all orbit
coefficients and cusp subtractions vanish for both $E_\tau$ and
$E_\tau^\iota$. The nonconstant upper tail is
an entire function of $s$ with meromorphic spectral dependence,
which proves the assertion about $I$ and its coefficients. Fix
$s_0$ with $\Re s_0>\beta(\tau_0)$. The inequalities defining
that region remain strict on a small product neighborhood and a
sufficiently small spectral circle. After clearing the genuine
spectral poles, the full Whittaker distribution is holomorphic
there by the local seminorm estimates; its remaining Euler factors
have no pole or zero in their denominators there. Cauchy's coefficient
formula preserves holomorphy in $s$. This proves the assertion for
$J$ and $Z$, including their spectral Laurent coefficients. Removable normalization poles are treated by
Lemma~\ref{lem:normalization-removal}.
\end{proof}

\appendix
\section{Comparison with smoothed Eisenstein series}\label{subsec:eis-smoothed}
In this appendix write $E^0(v,h)$ for the spherical Eisenstein
series induced from the trivial character pair, with section $1$
on $U$. This auxiliary family is distinct from the general twisting
family $E(\tau,h)$ of Section~\ref{sec:eisenstein}.
Shintani's smoothing method \cite[Lemma~2.9]{Shintani}, developed in
the adelic binary-cubic setting by Wright \cite{Wright} and
Yukie \cite{Yukie}, uses the family
\begin{equation}\label{eq:smoothed-definition}
 \mathscr E(\psi,w;g^{-1})=\frac1{2\pi i}
       \int_{\Re v=c}\frac{\psi(v)E^0(v,g^{-1})}{w-v}\,dv,
       \qquad \frac12<c<\Re w,
\end{equation}
where $\psi$ is entire and decreases faster than every inverse power on
closed vertical strips.  Put
\begin{equation}\label{eq:smoothed-residue-normalization}
 \gamma_k=\frac12\operatorname*{Res}_{u=1}\Lambda_k(u),\qquad
 \rho=\operatorname*{Res}_{\tau=1/2}E^0(\tau,g^{-1})
      =\frac{\gamma_k}{\Lambda_k(2)}.
\end{equation}
The factor $1/2$ in $\gamma_k$ is the Jacobian of $u=2\tau$.
We first record the estimate needed to use this residue under an
integral.

\begin{lemma}\label{lem:smoothed-cusp}
The family $\mathscr E(\psi,w;g^{-1})$ continues meromorphically to a
neighborhood of $w=1/2$, with constant residue $\rho\psi(1/2)$.
On a Siegel set of height $\lambda>1$ it has the expansion
\begin{equation}\label{eq:smoothed-cusp-expansion}
 \mathscr E(\psi,w;g^{-1})
   =\psi(w)m(w)\lambda^{1-2w}+\mathcal R(w,g),
\end{equation}
where $\mathcal R$ is holomorphic near $1/2$ and is rapidly decreasing
in $\lambda$, locally uniformly in $w$. The same decomposition holds
locally uniformly on compact subsets of $\Re w>1/2$, with a
holomorphic rapidly decreasing remainder there. In particular, if
$F\in L^1(G_\A^1/G_k)$, then
\begin{equation}\label{eq:smoothed-integrable-residue}
 \lim_{\substack{w\to1/2\\\Re w>1/2}}
 (w-1/2)\int_{G_\A^1/G_k}F(g)\mathscr E(\psi,w;g^{-1})\,dg^1
      =\rho\psi(1/2)\int_{G_\A^1/G_k}F(g)\,dg^1.
\end{equation}
\end{lemma}
\begin{proof}
Shifting the contour in \eqref{eq:smoothed-definition} to
$0<c'<1/2$ crosses only the constant residue at $v=1/2$, and gives
\begin{equation}
 \mathscr E(\psi,w;g^{-1})
 =\frac{\rho\psi(1/2)}{w-1/2}
       +\frac1{2\pi i}\int_{\Re v=c'}
                    \frac{\psi(v)E^0(v,g^{-1})}{w-v}\,dv.
\end{equation}
The second term is holomorphic near $w=1/2$ for fixed $g$.
To obtain a uniform cusp estimate, return to the original contour
$\Re v=c>1/2$ and first work in $\Re w>c$. Split
$E^0(v,g^{-1})$ there into its two constant-term powers and its
nonconstant part. The following shifts start from that original
line, not from the line $c'<1/2$ in the preceding continuation
formula; in particular they introduce no unrecorded crossing of
the scattering pole at $v=1/2$.  In the integral of
$\lambda^{1+2v}$ move the $v$-line to $\Re v=-N$.
The integrand has no pole to the left of the original contour, and
the resulting term is $O(\lambda^{1-2N})$.
For the scattering power $m(v)\lambda^{1-2v}$, move the contour to
$\Re v=R>\max(c,\Re w)$.  There are no scattering poles in
$\Re v>1/2$.  The pole of the kernel at $v=w$ gives
\begin{equation}
 \psi(w)m(w)\lambda^{1-2w}
  +\frac1{2\pi i}\int_{\Re v=R}
             \frac{\psi(v)m(v)\lambda^{1-2v}}{w-v}\,dv,
\end{equation}
and the contour integral is $O(\lambda^{1-2R})$.
The same shift for the nonconstant part gives
$\psi(w)E^0_{\mathrm n}(w,g^{-1})$ and a rapidly decreasing contour
integral.  These are holomorphic near $w=1/2$, since the residue of
$E^0(w,g^{-1})$ there is constant.  Choosing $N$ and $R$ as large as
necessary proves \eqref{eq:smoothed-cusp-expansion} initially.
For any compact subset of $\Re w>1/2$, choose the original line
with $1/2<c<\inf\Re w$ and the right-hand line with
$R>\sup\Re w$. This proves the stated locally uniform expansion
throughout that half-plane. Its right-hand side, with $R$ fixed
large, is meromorphic near $w=1/2$ and agrees with the already
continued smoothed family; the residue of the nonconstant part
is zero there, so its remainder is holomorphic.

For clarity, the contour shifts require bounds only on the strips
just described.  The inducing-power integrand is entire apart from
the kernel pole and is shifted to the left.  The scattering and
nonconstant integrands are shifted to the right, within
$\Re v>1/2$.  On such a right-hand line,
$|E^0(R+it,g^{-1})|\le E^0(R,g^{-1})$ by the defining series.
Group derivatives introduce polynomial factors in $t$; Fourier
expansion then gives rapid cusp decay for the nonconstant part with
polynomial dependence on $t$.  On the strip $c'\le\Re v\le c$,
removing the pole at $1/2$ gives polynomial vertical bounds by the
Fourier expansion and the gamma-factor estimates.  The denominator
$\Lambda_k(1+2v)$ has no zeros there.  The rapid vertical decay of
$\psi$ therefore makes the horizontal contour integrals tend to zero.

For $\Re w>1/2$ sufficiently close to $1/2$,
\eqref{eq:smoothed-cusp-expansion} implies
\begin{equation}
 |(w-1/2)\mathscr E(\psi,w;g^{-1})|
 \le C\lambda^{1-2\Re w}+C_M\lambda^{-M}\le C'.
\end{equation}
The compact part of the quotient satisfies the same bound by the
computed residue and local meromorphy.  Dominated convergence proves
\eqref{eq:smoothed-integrable-residue}.
\end{proof}

For $F=E(\tau)\Theta_S$, the leading cusp powers generally prevent
absolute integrability.  The next proposition identifies the residue
of its meromorphically continued smoothed pairing with the common
finite part of Section~\ref{subsec:eis-finite-part}.

\begin{proposition}\label{prop:smoothed-finite-part}
Fix $z>0$, let $\tau$ be a generic regular spectral parameter, and put
$\gamma=|\Re\tau|$.  The integral
\begin{equation}\label{eq:smoothed-pairing}
 \mathcal A_{\Phi,\tau}(w)
   =\int_{G_\A^1/G_k}E(\tau,g^{-1})\Theta_{S,\Phi}(z,g)
                            \mathscr E(\psi,w;g^{-1})\,dg^1
\end{equation}
converges absolutely for $\Re w>2+\gamma$.  It continues
meromorphically to a neighborhood of $w=1/2$ and satisfies
\begin{equation}\label{eq:eis-smoothed-residue}
 \operatorname*{Res}_{w=1/2}\mathcal A_{\Phi,\tau}(w)
                  =\rho\psi(1/2)R^{\rm reg}_{\Phi,\tau}(z).
\end{equation}
\end{proposition}
\begin{proof}
Choose a cusp cutoff $L>1$, and on $\lambda>L$ put
\begin{equation}
 \begin{split}
 F_z(g)&=E(\tau,g^{-1})\Theta_{S,\Phi}(z,g),\\
 B_z(g)&=z^{-2}A_\Phi(u)
                 \sum_{\epsilon=\pm1}
                    b_\epsilon(\tau;u,t_1,t_2)
                                      \lambda^{5+2\epsilon\tau}.
 \end{split}
\end{equation}
The singular-theta asymptotic of Section~\ref{subsec:eis-finite-part}
and rapid decay of the nonconstant Eisenstein term imply that
$F_z-B_z$ is rapidly decreasing in the cusp.
Lemma~\ref{lem:smoothed-cusp}, with its contour lines chosen far
enough to the right and left, shows that the leading cusp contribution
to \eqref{eq:smoothed-pairing} is
\begin{equation}\label{eq:smoothed-leading-integral}
 -\psi(w)m(w)z^{-2}
       \sum_{\epsilon=\pm1}D_\epsilon(\Phi,\tau)
       \frac{L^{4+2\epsilon\tau-2w}}{4+2\epsilon\tau-2w}.
\end{equation}
Indeed the cusp measure contains $\lambda^{-2}d\lambda/\lambda$,
and $\int_L^\infty\lambda^a d\lambda/\lambda=-L^a/a$ for
$\Re a<0$.  This also proves absolute convergence for
$\Re w>2+\gamma$.

The remaining cusp contributions are
\begin{equation}
 \psi(w)m(w)\int_{\lambda>L}
       (F_z-B_z)\lambda^{1-2w}\,dg^1
           +\int_{\lambda>L}F_z(g)\mathcal R(w,g)\,dg^1.
\end{equation}
Both converge locally uniformly near $w=1/2$ after clearing the
simple pole of $m(w)$.  Their residues are respectively
$\rho\psi(1/2)\int_{\lambda>L}(F_z-B_z)dg^1$ and zero.
The compact part contributes
$\rho\psi(1/2)\int_{\mathcal F_L}F_zdg^1$.
Since $3+2\epsilon\tau\ne0$ at a generic parameter, taking the
residue in \eqref{eq:smoothed-leading-integral} and adding the
remaining terms yields
\begin{equation}
 \rho\psi(1/2)\left\{
 \int_{\mathcal F_L}F_zdg^1+
 \int_{\lambda>L}(F_z-B_z)dg^1
 -z^{-2}\sum_{\epsilon=\pm1}D_\epsilon(\Phi,\tau)
            \frac{L^{3+2\epsilon\tau}}{3+2\epsilon\tau}
 \right\}.
\end{equation}
The expression in braces is the finite part in
\eqref{eq:eis-reg-subtraction}, proving
\eqref{eq:eis-smoothed-residue}.
\end{proof}

When $D_\epsilon$ is nonzero, at $\tau=\pm3/2$ a moving
cusp pole can meet $w=1/2$. If at least one of $\eta_1,\eta_2$ is nontrivial,
both $D_\epsilon$ vanish identically, and this particular
collision is absent.
The exceptional-parameter prescription of
Section~\ref{subsec:eis-exceptional} must then be applied to the complete
meromorphic expression; the two residues cannot simply be
interchanged.  Thus smoothing and function truncation agree at generic
parameters, and the cutoff unfolding shows that all orbit terms use
the common subtraction required by Poisson summation.
For Eisenstein twists the additional constant terms are the rational
expressions $H$ and $P_J$ in \eqref{eq:eis-I}--\eqref{eq:eis-J}.

\end{document}